\documentclass[1p,fleqn,10pt]{elsarticle}

\usepackage{pgfplots}
\usepackage{geometry}
\usepackage{hyperref}

\usepackage{lmodern}

\usepackage{amsthm}
\usepackage{amsmath}
\newtheorem{assumption}{Assumption}[section]
\newtheorem{definition}{Definition}[section]
\newtheorem{remark}{Remark}[section]
\newtheorem{theorem}{Theorem}[section]
\newtheorem{proposition}{Proposition}[section]
\newtheorem{lemma}{Lemma}[section]
\newtheorem{corollary}{Corollary}[section]

\usepackage{tikz}
\usepackage{tikz-cd}
\usetikzlibrary{patterns.meta,3d}
\usepgflibrary{patterns}
\usepackage{newtxmath}
\usepackage{graphicx}
\usepackage{xcolor}
\usepackage{mathtools}
\usepackage{enumitem}
\usepackage[framemethod=TikZ]{mdframed}

\usepackage{subcaption}
\usepackage{algorithm}
\usepackage{algpseudocode}

\usepackage{siunitx}
\newcommand{\pp}{p}
\newcommand{\kk}{m}
\newcommand{\qq}{q}

\renewcommand{\dots}{\ldots}

\newcommand{\Bezier}{B\'ezier}

\newcommand{\trunc}[2]{\mathop{\rm trunc}\nolimits_{#1}\left(#2\right)}

\newcommand{\interior}[1]{\mathop{\rm int}\nolimits\left(#1\right)}
\newcommand{\closure}[1]{\mathop{\rm clos}\nolimits\left(#1\right)}

\newcommand{\argmin}{\mathop{\rm argmin}}

\newcommand{\NN}{{\mathbb N}}

\newcommand{\RR}{{\mathbb R}}

\newcommand{\supp}[1]{\mathop{\rm supp}\nolimits\left(#1\right)}

\newcommand{\lspan}[1]{\mathop{\rm span}#1}
\newcommand{\act}{\mathop{\rm in}}
\newcommand{\nact}{\mathop{\rm ex}}

\newcommand{\DomainHierarchy}{\Psi} 

\newcommand{\piset}[3]{\widehat{\mathcal{I}}^{#1}_{#2}(#3)}

\newcommand{\peset}[3]{\widehat{\mathcal{E}}^{#1}_{#2}(#3)}

\newcommand{\Bspace}[2]{{\mathbb{B}_{#1}^{#2}}}
\newcommand{\Bbasis}[2]{{\mathcal{B}_{#1}^{#2}}}
\newcommand{\Bspline}[2]{{{B}_{#1}^{#2}}}

\newcommand{\Hbasis}[2]{{\mathcal{H}_{#1}^{#2}}}

\newcommand{\Tspace}[2]{{\mathbb{T}_{#1}^{#2}}}
\newcommand{\Tbasis}[2]{{\mathcal{T}_{#1}^{#2}}}
\newcommand{\Tspline}[2]{{{T}_{#1}^{#2}}}

\newcommand{\Spline}[2]{f_{#1}^{#2}}
\newcommand{\Coeff}[2]{\tau_{#1}^{#2}}

\newcommand{\mesh}[2]{\mathcal{M}_{#1}^{#2}}
\newcommand{\meshPbox}[2]{\mathcal{N}_{#1}^{#2}}

\newcommand{\pset}[2]{E_{#1}^{#2}}
\newcommand{\pbox}[2]{\Omega^{\pset{#1}{#2}}}

\newcommand{\Ebound}[3]{{H}_{#1}^{\mathrm{b}\ifx&#2&\else,\fi#2}\left(#3\right)}
\newcommand{\EboundNew}[3]{{H}^{\mathrm{b}\ifx&#2&\else,\fi#2}\left(#1, #3\right)}
\newcommand{\Tbound}[3]{{H}_{#1}^{\mathrm{t}\ifx&#2&\else,\fi#2}\left(#3\right)}

\renewcommand{\vec}[1]{\boldsymbol{#1}}

\newcommand{\sgn}[1]{\text{sgn}\left(#1\right)}

\newenvironment{assumpBox}
{
	\begin{assumption}
}%
{
	\end{assumption}
}

\title{Macro-element refinement schemes for THB-splines: Applications to B\'ezier projection and structure-preserving discretizations}
\author[label1]{Kevin Dijkstra}
\ead{kwdijkstra@tudelft.nl}
\author[label2]{Carlotta Giannelli}
\ead{carlotta.giannelli@unifi.it}
\author[label1]{Deepesh Toshniwal}
\ead{d.toshniwal@tudelft.nl}

\address[label1]{Delft Institute of Applied Mathematics, Delft University of Technology - The Netherlands}
\address[label2]{Dipartimento di Matematica e Informatica ‘‘U. Dini’’, Università degli Studi di Firenze - Italy}

\pgfplotsset{compat=1.18}

\begin{document}

\begin{abstract}
    This paper introduces a novel adaptive refinement strategy for Isogeometric Analysis (IGA) using Truncated Hierarchical B-splines (THB-splines).
    The strategy is motivated by the fact that certain applications may benefit from adaptive refinement schemes, which lead to a higher degree of structure in the locally-refined mesh than usual, and building this structure a priori can simplify the implementation in those contexts.
    Specifically, we look at two applications: formulation of an $L^2$-stable local projector for THB-splines a la \Bezier{} projection [Dijkstra and Toshniwal (2023)], and adaptive structure-preserving discretizations using THB-splines [Evans et al. (2020), Shepherd and Toshniwal (2024)].
    Previously proposed approaches for these applications require mesh modifications to preserve critical properties of the spline spaces, such as local linear independence or the exactness of the discrete de Rham complexes.
    Instead, we propose a macro-element-based refinement approach based on refining $\vec{q} = q_1\times\cdots\times q_n$ blocks of elements, termed $\vec{q}$-boxes, where the block size $\vec{q}$ is chosen based on the spline degree $\vec{p}$ and the specific application.
    \begin{itemize}
        \item For the \Bezier{} projection for THB-splines, we refine $\vec{p}$-boxes (i.e., $\vec{q} = \vec{p}$).
    We show that THB-splines are locally linearly independent on $\vec{p}$-boxes, which allows for a simple extension of the \Bezier{} projection algorithm to THB-splines. This new formulation significantly improves upon the approach previously proposed by Dijkstra and Toshniwal (2023).
        \item For structure-preserving discretizations, we refine $(\vec{p+1})$-boxes (i.e., $\vec{q} = \vec{p}+\vec{1}$).
        We prove that this choice of $\vec{q}$ ensures that the mesh satisfies the sufficient conditions presented in Shepherd and Toshniwal (2024) for guaranteeing the exactness of the THB-spline de Rham complex a priori and in an arbitrary number of dimensions.
        This is crucial for structure-preserving discretizations, as it eliminates the need for additional mesh modifications to maintain the exactness of the complex during adaptive simulations.
    \end{itemize}
    The effectiveness of the proposed framework is demonstrated through theoretical proofs and numerical experiments, including optimal convergence for adaptive approximation and the simulation of the incompressible Navier-Stokes equations.

\end{abstract}

\maketitle

\section{Introduction}

Isogeometric Analysis (IGA) was introduced to bridge the gap between Computer-Aided Design (CAD) and Finite Element Analysis (FEA) \cite{hughes_isogeometric_2005}. By using the same smooth splines for both geometric description and solution approximation, IGA offers superior approximation power per degree of freedom \cite{beirao_da_veiga_estimates_2011, sande_explicit_2020}.
However, standard tensor-product B-splines are computationally expensive in two or more dimensions.
This limitation can be overcome with locally refinable splines, such as T-splines \cite{buffa_linear_2010,sederberg_t-splines_2003}, LR-splines \cite{bressan_properties_2013,dokken_polynomial_2013,johannessen_isogeometric_2014,patrizi_linear_2020}, S-splines \cite{li_s-splines_2019}, and Hierarchical B-splines (HB-splines) \cite{forsey_hierarchical_1988, kraft_adaptive_1998, vuong_hierarchical_2011}.
Focusing on hierarchical B-splines and their truncated counterpart (THB-splines) \cite{giannelli_thb-splines_2012}, this paper proposes a new approach to adaptive refinement and demonstrates its use in adaptive approximation.

Adaptive methods using splines in IGA have been an active topic of research.
They have been investigated experimentally in \cite{bazilevs_isogeometric_2010,dorfel_adaptive_2010,vuong_hierarchical_2011,beirao_da_veiga_estimates_2011,kuru_goal-adaptive_2014,evans_hierarchical_2015} and in the case of (T)HB-splines, were mathematically investigated in \cite{buffa_adaptive_2016,buffa_adaptive_2017,gantner_adaptive_2017}.
Here, conditions that lead to a convergent adaptive method using THB-splines are presented. 
One of these conditions is that the mesh should be admissible, which, roughly speaking, implies that the mesh is graded, leading to the number of basis functions being uniformly bounded on each mesh element and influencing the sparsity of the system.

In the case of general adaptive methods, they consist of a \textit{solve} routine,  an \textit{estimate}, \textit{mark} and \textit{refine} routine \cite{devore_theory_2009, buffa_adaptive_2016} that, as their name suggests, solve the discrete problem, estimate the local error, mark elements for refinement, and refine the mesh.
However, depending on the problem considered, additional routines may be required.

For example, for time-dependent problems, the old solution at the previous time step is required to solve the subsequent time step.
This poses a challenge as this old solution might be defined over a different mesh, requiring an additional \textit{project} step.
A standard \textit{project} method could consist of a global $L^2$ or $H^1$ projection, but in the worst case, these methods scale cubically with the number of mesh elements or degrees of freedom. 
Compare this to local projectors, which instead scale linearly \cite{lanini_characterization_2024,giust_local_2020, thomas_bezier_2015,speleers_effortless_2016} and are thus more efficient.

Another example where an additional routine is required is in adaptive structure-preserving discretizations -- \textit{cohomology-preserving refinements and coarsenings}.
Structure-preserving discretizations aim to preserve the underlying structure of the continuous problem at the discrete level, such as conservation laws of physical quantities.
Formulation of such methods can be done within the framework of finite element exterior calculus (FEEC) \cite{arnold_finite_2018}, where the focus is on discretizing the Hilbert complex associated with the PDE of interest.
This approach has been successfully utilized within IGA to perform structure-preserving discretizations for various problems using tensor-product B-splines.
This means discretizing the de Rham complex for the incompressible Navier--Stokes equations \cite{evans2013isogeometric} and Maxwell's equations \cite{buffa2011isogeometric}.
However, the extension of these methods to locally-refinable splines is not straightforward, as the discrete spaces must form a valid de Rham complex, which is not guaranteed for general meshes \cite{evans2020hierarchical,shepherd_locally-verifiable_2024}.
The purpose of the additional \textit{cohomology-preserving refinements and coarsenings} routine is to ensure that the adaptively-refined spaces satisfy the necessary conditions for structure-preserving discretizations of the de Rham complex and, thus, lead to a stable discretization of the associated PDEs.
While such routines have been previously proposed for two-dimensional problems \cite{buffa2014isogeometric,johannessen2015divergence,cabanas_construction_2025}, doing so in arbitrary dimensions remains an open problem.

\subsection{Motivation and contributions of this paper}
In this paper, we propose a new, simple-to-implement refinement strategy for THB-splines that can be useful in the two example applications mentioned above.
The strategy is based on:
\begin{itemize}
    \item interpreting the element mesh as a collection of mesh of macro-elements called $\vec{q}$-boxes, which are blocks of $\vec{q} = (q_1,\dots,q_n) \in \NN^n$ mesh elements;
    \item limiting refinement and coarsening coarsening one or more $\vec{q}$-boxes in each adaptive step.
\end{itemize}
A specific version of this approach was proposed in \cite{giannelli_bases_2013} for studying the dimension of hierarchical B-spline spaces in two dimensions.
Clearly, this point of view is a generalization of the standard element-wise refinement, which corresponds to the case $\vec{q} = \vec{1}$.
However, other application-specific choices of $\vec{q}$ are possible, and we will illustrate the advantages of this generalization in the two example applications discussed above.

\paragraph{Application 1: Local projector for THB-splines.}~\\
In \cite{lanini_characterization_2024}, an extension of the B\'ezier projector proposed in \cite{thomas_bezier_2015} to the context of maximally-smooth THB-splines of admissibility class 2 was presented.
The \Bezier{} projector has the advantage of an intuitive and straightforward construction, and an efficient implementation that only utilizes extraction operators \cite{borden2011isogeometric,d2018multi}.
A key ingredient of the projector is that local projections are needed over mesh elements and requires the underlying spline basis to be locally linearly independent.
THB-splines do not possess this property, and \cite{lanini_characterization_2024} circumvented this issue by performing the local projections by partitioning the mesh into a set of non-overloaded macro-elements.
However, the existence of this partition is not guaranteed for general meshes, and mesh modifications that would ensure its existence are not obvious.
The constructive approach presented in \cite{lanini_characterization_2024} was thus limited to the case of two dimensions, and only for maximally smooth THB-splines of admissibility class 2.

In this paper we show that the proposed $\vec{q}$-box refinement strategy can overcome the above challenges: the choice of $\vec{q} = \vec{\pp}$, where $\vec{\pp}$ is the spline degree, enables a constructive formulation of the \Bezier{} projector for THB-splines of any admissibility class in arbitrary dimensions, while preserving its optimal approximation properties.
We test this formulation in several numerical examples, including adaptive approximation.

\paragraph{Application 2: Structure-preserving discretizations.}~\\
Structure-preserving adaptive discretizations using THB-splines have been recently proposed in \cite{cabanas_construction_2025} for two-dimensional problems.
Their approach is based on ensuring that all adaptive refinement steps satisfy the sufficient conditions presented in \cite{shepherd_locally-verifiable_2024} for guaranteeing the construction of a valid THB-spline de Rham complex.
We will show that the proposed $\vec{q}$-box refinement strategy can be used to extend this approach to arbitrary dimensions, and without the need for iterative mesh modifications -- we prove that any refinement or coarsening of $\vec{q}$-boxes with $\vec{q} = \vec{\pp} + \vec{1}$ leads automatic satisfaction of the sufficient conditions presented in \cite{shepherd_locally-verifiable_2024}.
We demonstrate the approach with a simple adaptive structure-preserving time-dependent incompressible Navier-Stokes method.

\subsection{Outline of the paper}
This paper is organised as follows.
In Section \ref{sec:theory}, we introduce the necessary background on THB-splines.
In Section \ref{sec:q-boxes}, we introduce the concept of $\vec{q}$-boxes and the associated refinement strategy.
In Section \ref{sec:non-overloaded-and-bezier-projection}, we present the \Bezier{} projector for THB-splines based on $\vec{p}$-box meshes, and in Section \ref{sec:struc-pres-methods}, we present the structure-preserving discretizations based on $(\vec{p}+\vec{1})$-box meshes.
Finally, in Section \ref{sec:numerical-results}, we present numerical experiments to demonstrate the performance of the proposed approaches.

\section{THB-splines}
\label{sec:theory}
The following section will introduce Truncated Hierarchical B-splines via (tensor) B-splines and the Hierarchical B-splines.

\subsection{Univariate B-splines}
\label{sec:univariate-b-spline}
Consider the unit domain $\Omega = (0,1)$ and a partition of $\Omega$ into $M$ elements by the breakpoints
\begin{equation}
    \widehat{\boldsymbol{\xi}} \,:=\, \left\{0<\widehat{\xi}_1<\dots<\widehat{\xi}_{M-1}<1 \right\}\,.
\end{equation}
We call $\widehat{\boldsymbol{\xi}}$ a breakpoint sequence and define the mesh elements $\Omega^{e} := \left(\hat{\xi}_{e},\hat{\xi}_{e+1}\right)$ such that $\mesh{}{} := \left\{0,1,\dots,M \right\}$ is the set of all mesh elements.
For a given set of breakpoints $\widehat{\boldsymbol{\xi}}$, polynomial degree $p\in\mathbb{N}_0$ and integer $1 \leq \kk \leq \pp +1$, consider the finite increasing sequence of real numbers over $\Omega$:
\begin{eqnarray}
    \label{eq:knot-sequence}
    \nonumber \boldsymbol{\xi}_{\pp,\kk}\,:=\,
    \{\underbrace{0, \dots, 0}_{\pp+1\,\text{times}} 
    < \underbrace{\widehat{\xi}_1= \dots= \widehat{\xi}_1}_{\kk\,\text{times}} 
    < \underbrace{\widehat{\xi}_2= \dots= \widehat{\xi}_2}_{\kk\,\text{times}} 
    <\dots  
    < \underbrace{1, \dots, 1}_{\pp+1\,\text{times}}\}\;.
\end{eqnarray}
This sequence is called a knot sequence, where each number is referred to as a knot.
The knot sequence's first and final knots have a multiplicity (the number of occurrences of a given knot) of $\pp+1$, and every interior knot has a multiplicity of $\kk\leq \pp+1$.
For our purposes, we only consider knot sequences generated from breakpoint sequences by increasing the multiplicity of the breakpoints.
One can also generate a breakpoint sequence from a knot sequence by reducing the multiplicity of the knots to 1.
For a knot sequence $\boldsymbol{\xi}_{\pp,\kk}$, we can define B-splines $\Bspline{j,\pp,\boldsymbol{\xi}_{\pp,\kk}}{}$, $j = 1, \dots, N := |\boldsymbol{\xi}_{\pp,\kk}|-\pp-1$, using the recursion,
\begin{equation}
    \Bspline{j,\pp,\boldsymbol{\xi}_{\pp,\kk}}{}(x)\,:=\, \frac{x - \xi_j}{\xi_{j+\pp} - \xi_j}\Bspline{j,\pp-1,\boldsymbol{\xi}_{\pp,\kk}}{}(x)\,+\,\frac{\xi_{j+\pp+1}-x}{\xi_{j+\pp+1}-\xi_{j+1}}\Bspline{j+1,\pp-1,\boldsymbol{\xi}_{\pp,\kk}}{}(x)\;,\quad j=1,\dots,N\;.
\end{equation}
The fractions are assumed to be zero whenever the denominator is zero. For the base case $\pp=0$, the B-splines are defined by the unit functions:
\begin{equation}
        \Bspline{j,0,\boldsymbol{\xi}_{\pp,\kk}}{}(x):= \begin{cases} 1& \text{if } x\in \left[\xi_j,\xi_{j+1}\right) \text{ and } j \neq N\;\text{ or }x\in \left[\xi_j,\xi_{j+1}\right] \text{ and } j = N\;, \\ 0 & \text{otherwise}\;,\end{cases}
\end{equation}

The functions $\Bspline{j,\pp,\boldsymbol{\xi}_{\pp,\kk}}{}$ are non-negative, form a partition of unity on $\closure{\Omega}$ and are $C^{\pp-\kk}$ smooth.
Lastly, most noteworthy for this paper, the collection of B-splines with support on any given mesh element $\Omega^e$ is linearly independent and can reproduce any polynomial of degree $\pp$ on $\Omega^e$.
The set containing all B-splines is denoted as $\Bbasis{\pp,\boldsymbol{\xi}_{\pp,\kk}}{}$, and its span is called the B-spline space and denoted as $\Bspace{\pp,\boldsymbol{\xi}_{\pp,\kk}}{}:= \lspan{\Bbasis{\pp,\boldsymbol{\xi}_{\pp,\kk}}{}}$.

\subsection{Tensor-product B-splines}
Consider the $n$-dimensional hypercube $\Omega = (0,1)^n$ in $\mathbb{R}^n$ and partition $\Omega$ via $\widehat{\Xi} := \left\{\widehat{\boldsymbol{\xi}}^1,\dots,\dots,\widehat{\boldsymbol{\xi}}^n\right\}$.
Here, each $\widehat{\boldsymbol{\xi}}^i$ partitions $\Omega$ along the $i$-th dimension.
Let the vector $\vec{\pp} := (\pp^1,\dots,\pp^n)\in \mathbb{N}^n$ denote the polynomial degrees per dimension, the vector $\vec{\kk} := (\kk^1,\dots,\kk^n)\in \mathbb{N}^n$ denote the interior knot multiplicities per dimension and let $\Xi_{\vec{\pp},\vec{\kk}} = (\boldsymbol{\xi}_{\pp^1,\kk^1}^1,\dots,\boldsymbol{\xi}_{\pp^n,\kk^n}^n)$ collect the knot sequences in each dimension, where each $\boldsymbol{\xi}_{\pp^i,\kk^i}^i$ is generated from the breakpoint sequence $\widehat{\boldsymbol{\xi}}^i$.
The tensor-product B-splines, $\Bspline{\vec{j},\vec{\pp},\Xi_{\vec{\pp},\vec{\kk}}}{}:= \bigotimes\limits_{i=1}^n\Bspline{j^i,\pp^i,\boldsymbol{\xi}^i_{\pp^i,\kk^i}}{}$, are then defined as tensor-products of univariate B-splines over the knot sequences $\boldsymbol{\xi}_{\pp^i,\kk^i}^i$, $i=1,\dots,n$.
To simplify notation, the spline degree $\pp$ and interior knot multiplicity $\kk$ are assumed to be fixed and are henceforth omitted. Hence, tensor-product B-splines will be denoted as $\Bspline{\vec{j},\Xi}{}$, the set of all tensor-product B-splines as $\Bbasis{\Xi}{}$, and the spline space spanned by them as $\Bspace{\Xi}{}$.
The mesh elements are denoted as $\Omega^{\vec{e}} = \bigtimes_{i=1}^n \left(\widehat{\xi}_{e^i}^i,\widehat{\xi}_{e^i+1}^i\right)$ for the index vector $\vec{e}$.
Collect these vectors $\vec{e}$ in $\mesh{\Xi}{}$ to define the set of mesh elements associated with $\Xi$.

\subsection{Truncated Hierarchical B-splines}
Consider a nested sequence of tensor-product B-spline spaces over the domain $\Omega = (0,1)^n$:
\begin{eqnarray}\label{eq:NestedBsplineSpaces}
	\Bspace{\boldsymbol{\Xi}_1}{} \subset \Bspace{\boldsymbol{\Xi}_2}{} \subset \dots \subset \Bspace{\boldsymbol{\Xi}_L}{}\;,
\end{eqnarray}
for appropriately chosen tensor-product level-$\ell$ knot sequences $\boldsymbol{\Xi}_\ell = (\boldsymbol{\xi}^1_\ell,\dots,\boldsymbol{\xi}^n_\ell)$ and associated breakpoint sequences $\widehat{\boldsymbol{\Xi}}_\ell$, $\ell = 1, \dots, L$.
By fixing the knot sequences at each level $\ell$, we simplify the notation of the B-spline spaces to $\Bspace{\ell}{} := \Bspace{\boldsymbol{\Xi}_\ell}{}$.
The corresponding level-$\ell$ B-spline basis will be denoted as $\Bbasis{\ell}{} := \Bbasis{\boldsymbol{\Xi}_\ell}{}$, the level-$\ell$ B-splines by $\Bspline{\vec{j},\ell}{} := \Bspline{\vec{j},\boldsymbol{\Xi}_\ell}{}$.
The level-$\ell$ mesh elements are denoted by $\vec{e}_\ell := [\vec{e}^T~ \ell]^T$ and collected in $\vec{e}_\ell\in \mesh{\ell}{}$ for all $\vec{e}\in \mesh{\Xi_\ell}{}$.
Then, for a vector $\vec{t}\in \mathbb{Z}^n$, define the addition $\vec{e}_l + \vec{t} := (\vec{e}+\vec{t})_\ell$ as the translated element of level $\ell$,
\begin{eqnarray}
    \Omega^{\vec{e}_\ell+\vec{t}} := \bigtimes_{i=1}^n \left(\widehat{\xi}_{e^i+t^i,\ell}^i,\widehat{\xi}_{e^i+t^i+1,\ell}^i\right)\;.
\end{eqnarray}

To ensure \eqref{eq:NestedBsplineSpaces}, we assume that the B-spline space $\Bspace{\ell}{}$ is attained by bisecting the mesh elements of the previous level.
In addition, we assume that at each level, the mesh is quasi-uniform.
\begin{assumpBox}
    \label{ass:bisection}
    The mesh element sizes at level $\ell$ satisfy the quasi-uniformity condition. That is, there exists a constant $\eta \geq 1$ such that for all breakpoints $\widehat{\xi}^i_{e^i,\ell}\in\widehat{\boldsymbol{\xi}}^i_\ell$, the mesh-size ratio is bounded as
    \begin{equation}
        \eta^{-1} \,\leq\, \frac{\widehat{\xi}_{e^i+1,\ell}^i - \widehat{\xi}_{e^i,\ell}^i}{\widehat{\xi}_{e^i,\ell}^i - \widehat{\xi}_{e^i-1,\ell}^i} \,\leq\, \eta\;,\quad i = 1,\dots, n\;,\quad \ell=1,\dots,L\;.
    \end{equation}
    Moreover, for $\ell > 1$, the level-$\ell$ breakpoint sequence in each direction is obtained by bisecting each $\Omega^{\vec{e}_{\ell-1}}$ of the corresponding level-$(\ell-1)$ breakpoint sequence.
\end{assumpBox}
Besides decreasing the element size, nested B-spline spaces can also be generated by increasing the spline degree $\vec{\pp}$ or the knot multiplicities $\vec{\kk}$.
These will not be considered.
\begin{assumpBox}\label{ass:only-h-refinement}
    The interior knot multiplicities $\vec{\kk}$ and spline degree $\vec{\pp}$ are the same for all $\boldsymbol{\Xi}_\ell$, $1\leq \ell\leq L$.
\end{assumpBox}
\begin{remark}
    The requirement on knot muliplicities imposed by Assumption \ref{ass:only-h-refinement} is placed to simplify the theoretical discussions in the different sections of the paper.
    For instance, the proofs in Sections \ref{sec:non-overloaded-and-bezier-projection} and \ref{sec:struc-pres-methods} can be extended to the case where knot multiplicities are different for different levels (assuming that such a choice still leads to a nested sequence of spaces).
    Moreover, the logic of the proofs remains the same, but the notation and the book-keeping become much more cumbersome -- this is why we have placed the above assumption on the knot multiplicities.
    The requirement on the degrees, however, is is mandatory for the current discussion: this choice ensures nestedness of the $\mathbf{q}$-boxes (see Section \ref{sec:q-boxes}) between levels, a key ingredient in our proofs.
\end{remark}

In addition, consider a sequence of nested, closed subsets of $\closure{\Omega}$:
\begin{eqnarray}
    \Omega_L \subseteq \dots \subseteq \Omega_1 := \closure{\Omega}\;,
\end{eqnarray}
where each $\Omega_\ell$ is given by a set of level-$(\ell-1)$ mesh elements
\begin{eqnarray}\label{eq:def-refinement-domain-elements}
    \Omega_\ell := \bigcup_{\vec{e}_{\ell-1}\in I_\ell} \closure{\Omega^{\vec{e}_{\ell-1}}}\;,
\end{eqnarray}
for some subset $I_\ell \subset \mesh{\ell-1}{}$. The collection of those subsets is denoted by:
\begin{eqnarray}
    \DomainHierarchy := \left\{\Omega_1,\dots,\Omega_L\right\}\;,
\end{eqnarray}
and will be referred to as the domain hierarchy on $\Omega$.
For a given level $\ell$, refinement domain $\Omega_k$ partitions the B-spline basis functions $\Bbasis{\ell}{}$ as:
\begin{subequations}\label{eq:active-nonactive-mesh-elements}
\begin{eqnarray}
\Bbasis{\ell}{\act,k} &:=& \left\{\,\Bspline{\vec{i},\ell}{} \in \Bbasis{\ell}{} \,:\, \supp{\Bspline{\vec{i},\ell}{}} \subseteq \Omega_k\,\right\}\;,\\
\Bbasis{\ell}{\nact,k} &:=& \left\{\,\Bspline{\vec{i},\ell}{} \in \Bbasis{\ell}{} \,:\, \supp{\Bspline{\vec{i},\ell}{}} \nsubseteq \Omega_k\,\right\}\,=\,\Bbasis{\ell}{}\backslash \Bbasis{\ell}{\act,k}\;.
\end{eqnarray}
\end{subequations}
Here, $\Bbasis{\ell}{\act,k}$ are the level-$\ell$ B-splines that are contained in $\Omega_k$, from which the HB-splines are defined. 
Note, in this work, we consider the support of a spline to be an open set.
\begin{definition}
    \label{def:HBsplines}
    Given a domain hierarchy $\DomainHierarchy$, the corresponding set of HB-spline basis functions is denoted by $\Hbasis{\DomainHierarchy}{}$ and defined recursively as follows:
    \begin{enumerate}[leftmargin=0.28in]
        \item $\Hbasis{1}{} := \Bbasis{1}{}$ ,
        \item for $\ell = 2,\dots,L$ :
        \[\Hbasis{\ell}{} := \,\Hbasis{\ell}{C}\, \cup \, \Hbasis{\ell}{F}\;,\]
        where
        \begin{eqnarray*}
            \Hbasis{\ell}{C} &:=& \left\{\, \Bspline{\vec{j},k}{} \in \Hbasis{\ell-1}{} \,:\, \supp{\Bspline{\vec{j},k}{}} \nsubseteq \Omega_\ell \, \right\}\;,
            \qquad
            \Hbasis{\ell}{F} := \Bbasis{\ell}{\act,\ell}\;,
        \end{eqnarray*}
        \item $\Hbasis{\DomainHierarchy}{} := \Hbasis{L}{}$ .
    \end{enumerate}
\end{definition}

\begin{figure}[htb]
    \centering
    \begin{subfigure}{0.4\textwidth}
        \begin{tikzpicture}
            \node at (0,0) {\includegraphics[width =\textwidth]{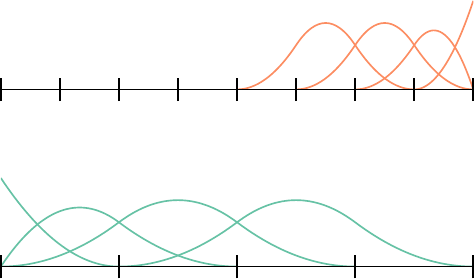}};
        \end{tikzpicture}
        \caption{HB-spline basis}
        \label{fig:HB-THB-comparison-HB}
    \end{subfigure}
    \hfill
    \begin{subfigure}{0.4\textwidth}
        \begin{tikzpicture}
            \node at (0,0) {\includegraphics[width =\textwidth]{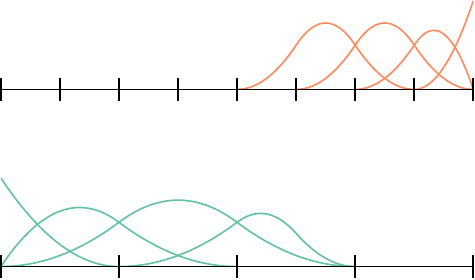}};
        \end{tikzpicture}
        \caption{THB-spline basis}
        \label{fig:HB-THB-comparison-THB}
    \end{subfigure}
    \caption{An HB-spline basis \ref{fig:HB-THB-comparison-HB} and a THB-spline basis \ref{fig:HB-THB-comparison-THB}. The green and orange highlighted basis functions represent the (T)HB-splines at levels 1 and 2, respectively.}
    \label{fig:HB-THB-comparison}
\end{figure}

See Figure \ref{fig:HB-THB-comparison-HB} for an example.
It is well-known that HB-spline basis functions lack the partition of unity property.
Moreover, the number of overlapping basis functions associated with different hierarchical levels increases quickly. This motivates the construction of a different basis,  the THB-spline basis. This basis relies on the following definition.
\begin{definition}
    Given $\ell=2,\dots,L$. Let $\Spline{}{}\in \Bspace{\ell-1}{}$ be represented in the B-spline basis $\Bbasis{\ell}{}$:
    \begin{eqnarray}
        \label{eq:defTruncationSum}
        \Spline{}{} = \;\sum_{\mathclap{\vec{j}:\Bspline{\vec{j},\ell}{}\in\Bbasis{\ell}{}}}\; \Coeff{\vec{j},\ell}{}\Bspline{\vec{j},\ell}{}\;.
    \end{eqnarray}
    The truncation of $\Spline{}{}$ at hierarchical level $\ell$ is defined as the sum of the terms appearing in \eqref{eq:defTruncationSum} corresponding to the B-splines in $\Bbasis{\ell}{\nact,\ell}$:
    \begin{eqnarray}
        \label{eq:def-trunc-operator}
        \trunc{\ell}{\Spline{}{}} := \;\sum_{\mathclap{\vec{j}:\Bspline{\vec{j},\ell}{}\in \Bbasis{\ell}{\nact,\ell}}}\; {\Coeff{\vec{j},\ell}{}\Bspline{\vec{j},\ell}{}}\;.
    \end{eqnarray}
\end{definition}
Since the B-spline spaces are nested, see \eqref{eq:NestedBsplineSpaces}, $\trunc{\ell}{\cdot}$ is well defined for the B-splines of $\Bbasis{1}{},\dots,\Bbasis{\ell}{}$. Exploiting this, the HB-splines are altered to form the THB-splines.
\begin{definition}
    \label{def:THBsplines}
    Given a domain hierarchy $\DomainHierarchy$, the corresponding set of THB-splines basis is denoted by $\Tbasis{\DomainHierarchy}{}$ and defined recursively as follows:
    \begin{enumerate}[leftmargin=0.28in]
        \item $\Tbasis{1}{} := \Bbasis{1}{}$ ,
        \item for $\ell = 2,\dots,L$ :
        \[\Tbasis{\ell}{} \,:=\, \Tbasis{\ell}{C} \,\cup\, \Tbasis{\ell}{F}\;,\]
        where
        \begin{eqnarray*}
            \Tbasis{\ell}{C} &:=& \left\{\,\trunc{\ell}{\Tspline{}{}}\,:\,\Tspline{}{} \in \Tbasis{\ell-1}\;,\;\supp{\Tspline{}{}} \nsubseteq \Omega_\ell \,\right\}\;,
            \qquad
            \Tbasis{\ell}{F} := \Bbasis{\ell}{\act,\ell}\;,
        \end{eqnarray*}
        \item $\Tbasis{\DomainHierarchy}{} := \Tbasis{L}{}$ .
    \end{enumerate}
\end{definition}
See Figure \ref{fig:HB-THB-comparison-THB} for an example.
The $N := |\Tbasis{\DomainHierarchy}{}|$ THB-spline basis functions are  indexed by some ordering $j = 1,\dots,N$ so that the THB-spline basis and space are given by
\begin{eqnarray}
    \Tbasis{\DomainHierarchy}{} := \left\{\Tspline{j}{}\right\}_{j=1}^N\;,\,\Tspace{\DomainHierarchy}{} := \lspan{\Tbasis{\DomainHierarchy}{}}\;.
\end{eqnarray}
Lastly, the set of active level-$\ell$ mesh elements are denoted by $\mesh{\ell}{\act}$ and is defined as:
\begin{equation}
    \mesh{\ell}{\act} :=
    \left\{\,
        {\vec{e}_{\ell}} \in \mesh{\ell}\,:\,\Omega^{\vec{e}_{\ell}} \subset \Omega_\ell \,\land\, \Omega^{\vec{e}_{\ell}} \cap \Omega_{\ell+1}= \emptyset
    \,\right\}\;.
\end{equation}
The collection of all active mesh elements is defined as:
\begin{equation}
    \mesh{\DomainHierarchy}{} := \left\{\,{\vec{e}_{\ell}}\in \mesh{\ell}{\act}\,:\, 1\leq \ell \leq L\,\right\}\;.
\end{equation}
Then, the admissibility class of a THB-spline space is defined as follows \cite{buffa_complexity_2016}.
\begin{definition}\label{def:admiss-class}
    A mesh is of admissible class $c$ if the truncated splines of $\Tbasis{\DomainHierarchy}{}$ which are non-zero on an element $\Omega^{\vec{e}_{\ell}}$ for ${\vec{e}_{\ell}}\in \mesh{\DomainHierarchy}{}$, belong to at most $c$ successive levels.
\end{definition}
\begin{remark}
Each intermediate set $\Tbasis{\ell}{}$ in Definition \ref{def:THBsplines} is again a THB-spline space. These are obtained by limiting the domain hierarchy $\DomainHierarchy$ to the first $\ell$ levels. We denote this domain hierarchy as $\DomainHierarchy_\ell$ and call the associated THB-spline space (respectively, basis) its level-$\ell$ intermediate THB-spline space (respectively, basis).
\end{remark}

\section{Macro-element-based adaptive refinement}\label{sec:q-boxes}
We introduce a refinement strategy in which larger clusters of elements, referred to as macro-elements, are refined collectively. 
For this, let $\vec{\qq}\in \mathbb{N}^n_{+}$ where $\mathbb{N}_{+} := \mathbb{N}\backslash\{0\}$, then, the macro-elements we will consider are defined as:
\begin{definition}\label{def:pbox}
    Given level $\ell$ and $\vec{\qq}\in\mathbb{N}_{+}^n$, assume that the breakpoint sequences $\widehat{\boldsymbol{\xi}}^i_\ell$ contains $\qq^i R^i_\ell+1$ breakpoints {where $\vec{R}_\ell \in \mathbb{N}_{+}^n$}.
    Then, for $r^i \in \{0,\dots,R^i_\ell-1\}$ define the \textbf{$\vec{\qq}$-box} $\pset{\vec{r},\ell}{}\subseteq \mesh{\ell}{}$,
	\begin{eqnarray}\label{eq:pbox-definition}
		\pbox{\vec{r},\ell}{} := \bigtimes_{i=1}^n \left(\widehat{\xi}^i_{r^i\qq^i,\ell},~\widehat{\xi}^i_{(r^i+1)\qq^i,\ell}\right) = \interior{\bigcup_{{\vec{e}_{\ell}} \in \pset{\vec{r},\ell}{}} \closure{\Omega^{\vec{e}_{\ell}}}} \;.
	\end{eqnarray}
    Denote the set of all possible level-$\ell$ $\vec{\qq}$-boxes by $\meshPbox{\ell}{\vec{\qq}} := \{~\pset{\vec{r},\ell}{}~:~\pbox{\vec{r},\ell}{}\subset\Omega~\}$.
\end{definition}
By Assumption \ref{ass:bisection}, each non-empty element is bisected under refinement.
This extends naturally to $\vec{\qq}$-boxes, as Figure \ref{fig:example-pbox-bisections} shows.
\begin{figure}[htb]
    \centering
    \begin{subfigure}{.45\textwidth}
        \centering
        \begin{tikzpicture}
            \def\w{6}
            \def\wp{3}
            \draw (0,1) -- (\w,1);
            \foreach \i in {0,0.25,0.5,0.75,1}{
                \draw (\i*\w,1-0.2) -- (\i*\w,1+0.2);
            }
            \draw (0,0) -- (\w,0);
            \foreach \i in {0,0.5,1}{
                \draw (\i*\w,-0.2) -- (\i*\w,0.2);
            }
            \foreach \i in {0,1}{
                \draw [decorate,decoration={brace,amplitude=5pt,mirror,raise=4ex}] (\i*\wp,0) -- ++ (\wp,0) node[midway,yshift=-3em]{$\pbox{\i,1}{}$};
            }
            \foreach \i in {0,1,2,3}{
                \draw [decorate,decoration={brace,amplitude=5pt,raise=4ex}] (\i*\wp/2,1) -- ++ (\wp/2,0) node[midway,yshift=3em]{$\Omega^{{\i_1}}{}$};
            }
        \end{tikzpicture}
        \caption{$\ell=1$}
    \end{subfigure}
    \hfill
    \begin{subfigure}{.45\textwidth}
        \centering
        \begin{tikzpicture}
            \def\w{6}
            \def\wp{1.5}
            \draw (0,1) -- (\w,1);
            \foreach \i in {0,0.125,0.25,0.375,0.5,0.625,0.75,0.875,1}{
                \draw (\i*\w,1-0.2) -- (\i*\w,1+0.2);
            }
            \draw (0,0) -- (\w,0);
            \foreach \i in {0,0.25,0.5,0.75,1}{
                \draw (\i*\w,-0.2) -- (\i*\w,0.2);
            }
            \foreach \i in {0,1,2,3}{
                \draw [decorate,decoration={brace,amplitude=5pt,mirror,raise=4ex}] (\i*\wp,0) -- ++ (\wp,0) node[midway,yshift=-3em]{$\pbox{\i,2}{}$};
            }
            \foreach \i in {0,1,2,3,4,5,6,7}{
                \draw [decorate,decoration={brace,amplitude=5pt,raise=4ex}] (\i*\wp/2,1) -- ++ (\wp/2,0) node[midway,yshift=3em]{$\Omega^{{\i_2}}{}$};
            }
        \end{tikzpicture}
        \caption{$\ell=2$}
    \end{subfigure}
    \caption{An element mesh $\Omega^{e_{\ell}}$ (top) and a $\vec{q}$-box mesh $\pbox{r,\ell}{}$ (bottom, $q=2$) of levels $\ell=1,2$. Like mesh elements, $\vec{q}$-boxes are bisected when refined.}
    \label{fig:example-pbox-bisections}
\end{figure}
This allows us to restrict refinement with the help of $\vec{\qq}$-boxes.
\begin{assumpBox}\label{ass:p-box-refinement-domain}
The refinement domain $\Omega_{\ell+1}$ consists of level-$\ell$ $\vec{\qq}$-boxes,
\begin{equation*}
    \Omega_{\ell+1} =  \bigcup_{\pset{\vec{r},\ell}{}\in R_\ell} \closure{\pbox{\vec{r},\ell}{}} \;,
\end{equation*}
where $R_\ell \subset \meshPbox{\ell}{\vec{\qq}}$.
\end{assumpBox}
Given a domain hierarchy $\DomainHierarchy$ adhering to Assumption \ref{ass:p-box-refinement-domain}, we define the following $\vec{\qq}$-boxes.
\begin{definition}
    The set of active level-$\ell$ $\vec{\qq}$-boxes, $\vec{\qq} \in \NN^n_{+}$, is defined as
    \begin{eqnarray}
        \meshPbox{\ell}{\vec{\qq},\act}:=\left\{ \,\pset{\vec{r},\ell}{}\in\meshPbox{\ell}{\vec{\qq}}\,:\,\pbox{\vec{r},\ell}{}\subset \Omega_\ell 
        \;\; \land \;\; \pbox{\vec{r},\ell}{}\nsubset \Omega_{\ell+1} \,\right\}\;,
    \end{eqnarray}
    and the collection of all active $\vec{\qq}$-boxes as the $\vec{\qq}$-box mesh
    \begin{eqnarray}
        \meshPbox{\DomainHierarchy}{\vec{\qq}} := \left\{\,\pset{\vec{r},\ell}{}\in\meshPbox{\ell}{\vec{\qq},\act}\,:\,1\leq\ell\leq L\,\right\}\;.
    \end{eqnarray}
    For $\ell>1$, call a $\vec{\qq}$-box $\pset{\vec{r},\ell}{} \in \meshPbox{\ell}{\vec{q}}$ a border $\vec{\qq}$-box if
    \begin{equation}
        \partial \pbox{\vec{r},\ell}{} \cap \left(\partial \Omega_\ell \backslash\partial \Omega \right) \neq \emptyset\;.
    \end{equation}
    Secondly, we call a border $\vec{\qq}$-box of level $\ell$ a well-behaved $\vec{\qq}$-box if it is not contained in any well-behaved $\vec{\qq}$-box of levels $\ell=2,\dots,\ell-1$.  
    Lastly, any active $\vec{\qq}$-box of $\meshPbox{\DomainHierarchy}{\vec{\qq}}$ that is not contained in any well-behaved $\vec{\qq}$-boxes is called a regular $\vec{\qq}$-box.
    Then, the set containing all well-behaved $\vec{\qq}$-boxes and regular $\vec{\qq}$-boxes partitions the mesh elements $\mesh{\DomainHierarchy}{}$.
\end{definition}

\begin{figure}[htb]
    \centering
    \begin{subfigure}{0.3\textwidth}
    \begin{tikzpicture}

    \definecolor{color1}{RGB}{102,194,165}
    \definecolor{color2}{RGB}{252,141,98}
    \definecolor{color3}{RGB}{141,160,203}

    \draw[dashed,gray,step=1cm] (0,0) grid (4,4);
    \draw[dashed,gray,step=0.5cm] (0,0) grid (2,2);
    \draw[dashed,gray,step=0.25cm] (0,2) grid (1,0) grid (2,1);
    \draw[dashed,gray,step=0.125cm] (1.5,0.5) grid (2,1);
    \draw[dashed,gray,step=0.125cm] (0,0) grid (.5,.5);

    \draw[step=2cm] (0,0) grid (4,4);
    \draw[step=1cm] (0,0) grid (2,2);
    \draw[step=0.5cm] (0,2) grid (1,0) grid (2,1);
    \draw[step=0.25cm] (1.5,0.5) grid (2,1);
    \draw[step=0.25cm] (0,0) grid (.5,.5);

    \end{tikzpicture}
    \caption{}
    \label{fig:example-border-and-well-behaved-p-boxes-mesh}
    \end{subfigure}
    \hfill
    \begin{subfigure}{0.3\textwidth}
    \begin{tikzpicture}

    \definecolor{color3}{RGB}{102,194,165}
    \definecolor{color1}{RGB}{0,0,0}
    \definecolor{color2}{RGB}{141,160,203}
    \definecolor{meshcolor}{RGB}{127,127,127}

    \draw[step=2cm,meshcolor] (0,0) grid (4,4);
    \draw[step=1cm,meshcolor] (0,0) grid (2,2);
    \draw[step=0.5cm,meshcolor] (0,2) grid (1,0) grid (2,1);
    \draw[step=0.25cm,meshcolor] (1.5,0.5) grid (2,1);
    \draw[step=0.25cm,meshcolor] (0,0) grid (.5,.5);

    \draw[line width = 0.5mm,color1,opacity = 1] (1,1) rectangle ++ (1,1);

    \draw[line width = 0.5mm,color1,opacity = 1] (0,1.5) rectangle ++ (0.5,0.5);
    \draw[line width = 0.5mm,color1,opacity = 1] (0.5,1.5) rectangle ++ (0.5,0.5);
    \draw[line width = 0.5mm,color1,opacity = 1] (0.5,1) rectangle ++ (0.5,0.5);
    \draw[line width = 0.5mm,color1,opacity = 1] (0.5,0.5) rectangle ++ (0.5,0.5);
    \draw[line width = 0.5mm,color1,opacity = 1] (1,0.5) rectangle ++ (0.5,0.5);
    \draw[line width = 0.5mm,color1,opacity = 1] (1.5,0) rectangle ++ (0.5,0.5);

    \draw[line width = 0.5mm,color1,opacity = 1] (0,0.25) rectangle ++ (0.25,0.25);
    \draw[line width = 0.5mm,color1,opacity = 1] (0.25,0.25) rectangle ++ (0.25,0.25);
    \draw[line width = 0.5mm,color1,opacity = 1] (0.25,0) rectangle ++ (0.25,0.25);
    \draw[line width = 0.5mm,color1,opacity = 1] (1.5,0.5) rectangle ++ (0.25,0.25);
    \draw[line width = 0.5mm,color1,opacity = 1] (1.75,0.5) rectangle ++ (0.25,0.25);
    \draw[line width = 0.5mm,color1,opacity = 1] (1.5,0.75) rectangle ++ (0.25,0.25);
    \draw[line width = 0.5mm,color1,opacity = 1] (1.75,0.75) rectangle ++ (0.25,0.25);

    \end{tikzpicture}
    \caption{}
    \label{fig:example-border-and-well-behaved-p-boxes-border-p-boxes}
    \end{subfigure}
    \hfill
    \begin{subfigure}{0.3\textwidth}
    \begin{tikzpicture}

    \definecolor{color1}{RGB}{102,194,165}
    \definecolor{color2}{RGB}{252,141,98}
    \definecolor{color3}{RGB}{141,160,203}
    \definecolor{color3}{RGB}{0,0,0}
    \definecolor{meshcolor}{RGB}{127,127,127}




    \draw[step=2cm,meshcolor] (0,0) grid (4,4);
    \draw[step=1cm,meshcolor] (0,0) grid (2,2);
    \draw[step=0.5cm,meshcolor] (0,2) grid (1,0) grid (2,1);
    \draw[step=0.25cm,meshcolor] (1.5,0.5) grid (2,1);
    \draw[step=0.25cm,meshcolor] (0,0) grid (.5,.5);

    \draw[line width = 0.5mm,color3,opacity = 1] (1,1) rectangle ++ (1,1);
    \draw[line width = 0.5mm,color3,opacity = 1] (0,1) rectangle ++ (1,1);
    \draw[line width = 0.5mm,color3,opacity = 1] (1,0) rectangle ++ (1,1);

    \draw[line width = 0.5mm,color3,opacity = 1] (0.5,0.5) rectangle ++ (0.5,0.5);

    \draw[line width = 0.5mm,color3,opacity = 1] (0,0.25) rectangle ++ (0.25,0.25);
    \draw[line width = 0.5mm,color3,opacity = 1] (0.25,0.25) rectangle ++ (0.25,0.25);
    \draw[line width = 0.5mm,color3,opacity = 1] (0.25,0) rectangle ++ (0.25,0.25);

    \end{tikzpicture}
    \caption{}
    \label{fig:example-border-and-well-behaved-p-boxes-well-behaved-p-box}
    \end{subfigure}
    \caption{In (a), a $\vec{q}$-box mesh is depicted for $\vec{q}=(2,2)$ and the associated element mesh with dashed lines. 
    In (b) and (c), respectively, the active border and well-behaved $\vec{q}$-boxes are highlighted. Observe that well-behaved $\vec{q}$-boxes can either be active or non-active.}
    \label{fig:example-border-and-well-behaved-p-boxes}
\end{figure}

See Figure \ref{fig:example-border-and-well-behaved-p-boxes} for examples of these $\vec{\qq}$-boxes.
We note that when $\vec{\qq} = \vec{1}$, the resulting $\vec{\qq}$-boxes are single mesh-elements.
Two specific choices of $\vec{\qq}$ will be considered: $\vec{\qq}=\vec{\pp}$ (Section 4) and $\vec{\qq}=\vec{\pp}+\mathbf{1}$ (Section 5), leading to $\vec{\pp}$- and $(\vec{\pp}+\vec{1})$-boxes related to the spline degree $\vec{\pp}$.
We extend Definition \ref{def:admiss-class} to $\vec{\qq}$-boxes.
\begin{definition}
    A $\vec{\qq}$-box mesh $\meshPbox{\DomainHierarchy}{\vec{\qq}}$ is of admissibility class $c$ if the truncated splines of $\Tbasis{\DomainHierarchy}{}$ which are non-zero on an active $\vec{\qq}$-box $\pset{\vec{r},\ell}{}$, belong to at most $c$ successive levels.
\end{definition}
In addition, to preserve an admissibility class upon refinement or coarsening, we extend algorithms 3, 4, and 5 from \cite{carraturo_suitably_2019} to $\vec{\qq}$-boxes, considering a marking strategy based on $\vec{\qq}$-boxes rather than elements.
They extend naturally to $\vec{\qq}$-boxes by altering the relevant definitions of multilevel support extension $S$, refinement neighborhood $N_r$, coarsening neighborhood $N_c$, see equations (3.3)--(3.5) in \cite{carraturo_suitably_2019}, as follows:
\begin{align}\label{eq:multilevel-support-pbox}
    S(\pset{\vec{r},\ell}{}, k) &:= \left\{\pset{\vec{s},k}{}\,:\,\exists\Bspline{\vec{i},k}{}\in \Bbasis{k}{},\supp{\Bspline{\vec{i},l}{}}\cap \pbox{\vec{s},k}{}\neq \emptyset \land  \supp{\Bspline{\vec{i},l}{}}\cap \pbox{\vec{r},\ell}{} \neq \emptyset\, \right\}\;,\nonumber\\
    N_r(\meshPbox{\DomainHierarchy}{\vec{\qq}},\pset{\vec{r},\ell}{},c) &:=\left\{\pset{\vec{t},\ell-c+1}{}\in \meshPbox{\ell-c+1}{\vec{\qq},\act}\,:\,\exists \pset{\vec{s},\ell-c+2}{}\in S(\pset{\vec{r},\ell}{},\ell-c+2),\,\pbox{\vec{s},\ell-c+2}{}\subseteq \pbox{\vec{t},\ell-c+1}{}\,\right\}\;,\\
    N_c(\meshPbox{\DomainHierarchy}{\vec{\qq}},\pset{\vec{r},\ell}{},c) &:=\left\{\pset{\vec{t},\ell+c}{}\in\meshPbox{\ell+c}{\vec{\qq},\act}\,:\,\exists \pset{\vec{s},\ell+1}{}\in \meshPbox{\ell+1}{\vec{\qq},\act}\land \pbox{\vec{s},\ell+1}{}\subset \pbox{\vec{r},\ell}{},\,\pset{\vec{t},\ell+c}{}\in S(\pset{\vec{s},\ell+1}{},\ell+1)\,\right\}\;.\nonumber
\end{align}
Lastly, there are two minor results regarding the admissibility classes for $\vec{\qq}$-boxes.
\begin{proposition}\label{prop:regular-p-box-spline-level}
    Consider THB-splines of degree $\vec{\pp}$ and $\vec{\qq}$-box mesh $\meshPbox{\DomainHierarchy}{\vec{\qq}}$ with $\qq^i\geq \pp^i$ for all $i$. Every level-$\ell$ regular $\vec{\qq}$-box $\pset{\vec{r},\ell}{}\in \meshPbox{\DomainHierarchy}{\vec{\qq}
    }$ only supports THB-splines of level $\ell$.
\end{proposition}
\begin{proposition}\label{prop:sufficient-condition-admiss-class-2}
    For THB-splines of degree $\vec{\pp}$ and $\vec{\qq}$-box mesh $\meshPbox{\DomainHierarchy}{\vec{\qq}}$ with $\qq^i\geq \pp^i$ for all $i$, the following are equivalent
    \begin{itemize}
        \item every well-behaved $\vec{\qq}$-box is (an) active (border $\vec{\qq}$-box),
        \item $\meshPbox{\DomainHierarchy}{\vec{\qq}}$ is of admissibility class $c=2$. 
    \end{itemize}
\end{proposition}
The proofs of Proposition \ref{prop:regular-p-box-spline-level} and Proposition \ref{prop:sufficient-condition-admiss-class-2} are presented in \ref{sec:AppendixA}.

\section{Locally linearly independent macro-elements and applications to B\'ezier projection}
\label{sec:non-overloaded-and-bezier-projection}
The first main result of this paper is that we can use the $\vec{\qq}$-boxes to characterize the local linear independence of THB-splines. 
For this, we consider THB-splines of degree $\vec{\pp}$ and $\vec{\qq}$-boxes with $\vec{\qq} =\vec{\pp}$, which we call $\vec{\pp}$-boxes.
Then, we introduce the notion of overloading:
\begin{definition}\label{def:overloaded}
    For a given set of functions, $\mathcal{S}$ on $\Omega$, we say that $\widetilde{\Omega}\subset\Omega$ is overloaded w.r.t. $\mathcal{S}$ if
    $$ S_{\tilde{\Omega}} = \left\{~f\vert_{\widetilde{\Omega}}~:~f\in\mathcal{S},~f\vert_{\widetilde{\Omega}}\neq 0~ \right\}$$
    is linearly dependent; else, $\widetilde{\Omega}$ is not overloaded.
\end{definition}
We can now state the local linear independence result for $\vec{\pp}$-boxes:
\begin{theorem}\label{thm:non-overloaded-well-behaved-p-box}
	Given a $\vec{\pp}$-box mesh $\meshPbox{\DomainHierarchy}{\vec{\pp}}$, neither the well-behaved $\vec{\pp}$-boxes nor the regular $\vec{\pp}$-boxes are overloaded w.r.t. $\Tbasis{\DomainHierarchy}{}$.
\end{theorem}
The proof of Theorem \ref{thm:non-overloaded-well-behaved-p-box} is split into two parts.
First, we consider the special case of admissibility class $c=2$, for which the proof is much simpler. 
Then, regarding the full proof of Theorem \ref{thm:non-overloaded-well-behaved-p-box}, it can be decomposed into smaller problems, where each problem is of admissibility class $c=2$.
\subsection{The case of admissibility class 2}
When we assume that the mesh is of admissibility class $2$, by Proposition \ref{prop:sufficient-condition-admiss-class-2}, Theorem \ref{thm:non-overloaded-well-behaved-p-box} reduces to the following Lemma.
\begin{lemma}\label{lem:non-overloaded-border-p-box}
    For a $\vec{\pp}$-box mesh $\meshPbox{\DomainHierarchy}{\vec{\pp}}$ of admissibility class 2,
    neither the active border $\vec{\pp}$-boxes nor the regular $\vec{\pp}$-boxes are overloaded w.r.t. $\Tbasis{\DomainHierarchy}{}$.
\end{lemma}
For the proof of Lemma \ref{lem:non-overloaded-border-p-box}, we note that the $\vec{\pp}$-boxes considered can each be subdivided into smaller macro-elements, and most importantly, into the macro-elements introduced in \cite{lanini_characterization_2024}.
In \cite{lanini_characterization_2024}, alternative macro-elements were introduced that are not overloaded. 
As a result, Lemma \ref{lem:non-overloaded-border-p-box} is trivial under the assumptions of \cite{lanini_characterization_2024}, where only maximal smoothness was considered. Since the proof of Lemma \ref{lem:non-overloaded-border-p-box} for the more general case of arbitrary smoothness follows similar arguments, we present it in \ref{sec:appendixB}.

\subsection{The general admissibility class case}
From Proposition \ref{prop:regular-p-box-spline-level}, regular $\vec{\pp}$-boxes only support THB-splines from a single level, which must be B-splines, and are thus never overloaded.
Hence, to prove Theorem \ref{thm:non-overloaded-well-behaved-p-box}, we will show that well-behaved $\vec{\pp}$-boxes are not overloaded.
For this, observe that a well-behaved $\vec{\pp}$-box $\pset{\vec{r},\ell}{}$ of level $\ell$ might consist of active border $\vec{\pp}$-boxes, active regular $\vec{\pp}$-boxes, deactivated border $\vec{\pp}$-boxes, and deactivated regular $\vec{\pp}$-boxes of various levels.
See Figure \ref{fig:motivation-well-behaved-p-box}.
But, if we consider one of these deactivated $\vec{\pp}$-boxes of level $k$, we note that they were active for the intermediate THB-spline space of level $k$.
We will use this observation to prove Theorem \ref{thm:non-overloaded-well-behaved-p-box}, by considering these intermediate THB-spline spaces of levels $\ell,\ell+1,\dots,L$, and showing that the well-behaved border $\vec{\pp}$-box is non-overloaded for each of these levels via induction.
However, we do not need to check all levels $\ell,\ell+1,\dots, L$, in light of the following definition.

\begin{figure}[htb]
    \centering
    \begin{tikzpicture}
        \def\w{3}
        \def\h{2}
        \def\xoffsetPbox{-6}
        \def\xoffsetBasis{6.5}
        \def\yoffsetBasis{-.5}

        \draw (\xoffsetPbox,0) -- ++ (\w,0) --++ (0,\w) --++(-\w,0) --cycle;
        \draw (\xoffsetPbox,+\w/2) -- ++ (\w,0);
        \draw (\xoffsetPbox+\w/2,0) -- ++ (0,\w);
        \draw (\xoffsetPbox+\w/4,0) -- ++ (0,\w);
        \draw (\xoffsetPbox,+\w/4) -- ++ (\w/2,0);
        \draw (\xoffsetPbox,+3*\w/4) -- ++ (\w/2,0);
        
        \draw[fill = gray!30!white] (0,0,0) -- (\w,0,0) -- (\w,0,\w) -- (0,0,\w) -- cycle;

        \draw[fill = gray!30!white] (0,\h,\w/2) -- (\w,\h,\w/2) -- (\w,\h,\w) -- (0,\h,\w) -- cycle;
        \draw (0,\h,0) -- (\w,\h,0) -- (\w,\h,\w) -- (0,\h,\w) -- cycle;
        \draw (\w/2,\h,0) -- (\w/2,\h,\w);
        \draw (0,\h,\w/2) -- (\w,\h,\w/2);

        \draw[fill = gray!30!white] (\w/4,2*\h,0) -- (\w/4,2*\h,3*\w/4) -- (0,2*\h,3*\w/4) -- (0,2*\h,\w) -- (\w/2,2*\h,\w) -- (\w/2,2*\h,0) -- cycle;
        \draw (0,2*\h,0) -- (\w/2,2*\h,0) -- (\w/2,2*\h,\w) -- (0,2*\h,\w) -- cycle;
        \draw (\w/4,2*\h,0) -- (\w/4,2*\h,\w);
        \draw (0,2*\h,\w/4) -- (\w/2,2*\h,\w/4);
        \draw (0,2*\h,\w/2) -- (\w/2,2*\h,\w/2);
        \draw (0,2*\h,3*\w/4) -- (\w/2,2*\h,3*\w/4);

        \node at (\xoffsetBasis,\yoffsetBasis)      {level $\ell\,:\,\Tbasis{\ell}{} := \Bbasis{\ell}{\act,\ell} \cup \trunc{\ell}{\Bbasis{\ell-1}{\nact,\ell}}$};
        \node at (\xoffsetBasis,\yoffsetBasis+\h)   {level $\ell+1\,:\,\Tbasis{\ell+1}{} := \Bbasis{\ell+1}{\act,\ell+1} \cup \trunc{\ell+1}{\Tbasis{\ell}{}}$};
        \node at (\xoffsetBasis,\yoffsetBasis+2*\h) {level $\ell+2\,:\,\Tbasis{\ell+2}{} := \Bbasis{\ell+2}{\act,\ell+2} \cup \trunc{\ell+2}{\Tbasis{\ell+1}{}}$};
    \end{tikzpicture}
    \caption{On the left a well-behaved $\vec{\pp}$-box $\pbox{\ell}{}$ of level $\ell$ is depicted. The middle depicts the active and deactivated $\vec{\pp}$-boxes in $\pbox{\ell}{}$ per level. The border $\vec{\pp}$-boxes are highlighted in grey. On the right, the intermediate THB-spline spaces $\Tbasis{k}{}$ are given for the levels $\ell,\ell+1$ and $\ell+2$.}
    \label{fig:motivation-well-behaved-p-box}
\end{figure}

\begin{definition}\label{def:class-well-behaved-p-box}
Consider a $\vec{\pp}$-box mesh $\meshPbox{\DomainHierarchy}{\vec{\pp}}$ and a well-behaved $\vec{\pp}$-box $\pset{\vec{r},\ell}{}$. Define the depth of $\pset{\vec{r},\ell}{}$ as the smallest $d\geq 0$ such that $\pset{\vec{r},\ell}{}$ contains an active border $\vec{\pp}$-box of level $\ell+d$.
\end{definition}
We will prove Theorem \ref{thm:non-overloaded-well-behaved-p-box} via induction.
We start by showing that $\pbox{\vec{r},\ell}{}$ is not overloaded w.r.t. $\Tbasis{\DomainHierarchy_\ell}{}$, the intermediate THB-spline basis of level $\ell$.
Then, for the induction step, we show that for $\ell < k\leq \ell + d$, $\pbox{\vec{r},\ell}{}$ is not overloaded w.r.t. $\Tbasis{\DomainHierarchy_k}{}$.
Both results rely on Lemma \ref{lem:non-overloaded-border-p-box}.
However, for the induction step, we require some intermediate results.\\

For the induction step, we note that any level-$(k-1)$ THB-spline $\Tspline{}{}\in \Tbasis{\DomainHierarchy_{k-1}}{}$ can be written as a linear combination of level-$(k-1)$ B-splines.
\begin{eqnarray}\label{eq:THB-spline-level-k-B-spline-comb}
    \Tspline{}{} = \quad\sum_{\mathclap{\vec{j}:\Bspline{j,k-1}{}\in \Bbasis{k-1}{}}}\; \Coeff{\vec{j}}{} \Bspline{\vec{j},k-1}{}\;.
\end{eqnarray}
As a result, we start by investigating B-splines before considering THB-splines.
Note that if a level-$(k-1)$ THB-spline has non-zero truncation on $\pset{\vec{r},\ell}{}$, at least one of the level-$(k-1)$ B-splines of \eqref{eq:THB-spline-level-k-B-spline-comb} must be a member of $\Bbasis{k-1}{\nact,k}$.
For these B-splines, we have the following result.

\begin{lemma}\label{lem:sufficient-projection-elements-level-k+1}
Consider a level-$\ell$ well-behaved $\vec{\pp}$-box $\pset{\vec{r},\ell}{}$ of depth $d$ and a B-spline $\Bspline{}{} \in \Bbasis{k-1}{\nact,k}$ for $\ell < k \leq \ell + d$ with support on $\pbox{\vec{r},\ell}{} \cap \Omega_{k}$. Then, $\Bspline{}{}$ has support on a level-$k$ border $\vec{\pp}$-box $\pset{\vec{t},k}{}$ and $\pbox{\vec{t},k}{}\subset \pbox{\vec{r},\ell}{}$.
\end{lemma}

\begin{proof}
The B-spline $\Bspline{}{}\in \Bbasis{k-1}{\nact,k}$ has support inside and outside of $\Omega_k$, so that,
\begin{equation}
    \supp{\Bspline{}{}} \cap  \left( \partial \Omega_k\backslash\partial\Omega \right) \neq \emptyset\;.
\end{equation}
Note that $\Bspline{}{}$ has support on a level-$k$ element $\Omega^{\vec{e}_{k}} \subset \Omega_{k}\cap \pbox{\vec{r},\ell}{}$. 
This element is contained in a level-$(k-1)$ element $\Omega^{\star {e}_{k}}  \subset \Omega_{k}\cap \pbox{\vec{r},\ell}{}$ which is contained in the $\vec{\pp}$-box ${\pbox{\vec{s},k-1}{}} \subseteq \Omega_{k}\cap \pbox{\vec{r},\ell}{}$. 
Note that the support of a level-$(k-1)$ B-spline counts at most $p^i+1$ level-$(k-1)$ mesh elements in each dimension $i=1,\dots,n$.
Then, since $\Bspline{}{}$ has support on $\pset{\vec{s},k-1}{}$, the only other level-$(k-1)$ $\vec{\pp}$-boxes it can have support on are directly adjacent to $\pbox{\vec{s},k-1}{}$. Hence, we have that,
\begin{equation}\label{eq:proof-317-eq1}
    \supp{\Bspline{}{}} \cap \left(  \partial \pbox{\vec{s},k-1}{}  \cap \left( \partial \Omega_k\backslash\partial\Omega \right) \right) \neq \emptyset\;.
\end{equation}
Since \eqref{eq:proof-317-eq1} is non-empty, there is a $\pbox{\vec{t},k}{}\subset \pbox{\vec{s},k-1}{}$ such that
\begin{equation}
    \supp{\Bspline{}{}} \cap \left(  \partial \pbox{\vec{t},k}{}  \cap \left( \partial \Omega_k\backslash\partial\Omega \right) \right) \neq \emptyset\;,
\end{equation}
showing that $\pbox{\vec{t},k}{}\subset\pbox{\vec{r},\ell}{}$ is a border $\vec{\pp}$-box.
\end{proof}
Next, we show that $\trunc{k}{\Bbasis{k-1}{\nact,k}}$ is linearly independent over a well-behaved $\vec{\pp}$-box $\pset{\vec{r},\ell}{}$.
\begin{lemma}\label{lem:linear-independence-B-k-1-nact-k}
Consider a level-$\ell$ well-behaved $\vec{\pp}$-box $\pset{\vec{r},\ell}{}$ of depth $d$  and $\ell < k \leq \ell + d $. Then $\trunc{k}{\Bbasis{k-1}{\nact,k}}$ is linearly independent on $\pset{\vec{r},\ell}{}$.
\end{lemma}
\begin{proof}
    Consider any linear combination $\Spline{}{} = \sum_{\vec{j}} \Coeff{\vec{j}}{} \Bspline{\vec{j}}{}$ with $\Bspline{\vec{j}}{} \in \Bbasis{k-1}{\nact,k}$ such that $\trunc{k}{\Spline{}{}} = 0$. 
    The splines $\Bspline{\vec{j}}{}\in \Bbasis{k-1}{\nact,k}$ which have support on $\pbox{\vec{r},\ell}{}\backslash\Omega_{k}$ are unaffected by $\trunc{k}{\Bspline{\vec{j}}{}}$ on $\pbox{\vec{r},\ell}{}\backslash\Omega_{k}$. 
    Hence, local linear independence of these $\Bspline{\vec{j}}{}$ are preserved, so that $\Coeff{\vec{j}}{} = 0$. 
    By Lemma \ref{lem:sufficient-projection-elements-level-k+1}, each remaining B-spline $\Bspline{\vec{j}}{}$ must have support on some level-$k$ border $\vec{\pp}$-box $\pbox{\vec{s},k}{}\subset\Omega_k$. 
    But then, the claim follows if $\trunc{k}{\Bbasis{k-1}{\nact,k}}$ is linearly independent on $\pset{\vec{s},k}{}$.\\
    
    For this, consider the following domain hierarchy $\DomainHierarchy^\# = \{\Omega, \Omega,\dots,\Omega_k\}$ where $|\DomainHierarchy^\#| = k$.
    By construction, the associated THB-spline space is given by
    \begin{equation}\label{eq:proof-linear-independence-B-k-1-nact-k-basis}
        \trunc{k}{\Bbasis{k-1}{\nact,k}} \cup \Bbasis{k}{\act,k}\;,
    \end{equation}
    and is of admissibility class $c=2$.
    But then, by Lemma \ref{lem:non-overloaded-border-p-box}, every active $\vec{\pp}$-box, and in particular $\pset{\vec{s},k}{}$, is not overloaded w.r.t. \eqref{eq:proof-linear-independence-B-k-1-nact-k-basis}.
\end{proof}
We extend this result to linear combinations of B-splines (e.g., THB-splines).
\begin{corollary}\label{cor:linear-independence-A-subset-B-k-1-nact-k}
    Consider a level-$\ell$ well-behaved $\vec{\pp}$-box $\pset{\vec{r},\ell}{}$ of depth $d$ and let $\mathcal{A} \subset \lspan{\Bbasis{k-1}{\nact,k}}$, $\ell < k \leq \ell + d $ be a linear independent set on $\pset{\vec{r},\ell}{}$. Then $\trunc{k}{\mathcal{A}}$ is linearly independent on $\pset{\vec{r},\ell}{}$.
\end{corollary}
\begin{proof}
    By the linear independence of $\mathcal{A}$, every $\Spline{i}{} \in \mathcal{A}$ can be written as $\Spline{i}{} = \sum_{\vec{j}} Q_{i\vec{j}} B_{\vec{j}}$ for $\Bspline{\vec{j}}{} \in \Bbasis{k-1}{\nact,k}$ where the rows of $Q$ are linearly independent. Then, consider any linear combination of splines of $\trunc{k}{\mathcal{A}}$ that sums to zero:
    \begin{eqnarray}
        \sum_i \Coeff{i}{}\trunc{k}{\Spline{i}{}}\;=\;
        \sum_i \sum_{\vec{j}} \Coeff{i}{} Q_{i\vec{j}} \trunc{k}{\Bspline{\vec{j}}{}} \;=\;0 \;.
    \end{eqnarray}
    By the linear independence of $\trunc{k}{\Bbasis{k-1}{\nact,k}}$ over $\pset{\vec{r},\ell}{}$ shown in Lemma \ref{lem:linear-independence-B-k-1-nact-k}, this equation can only be zero if $\Coeff{i}{}=0$, showing linear independence of $\trunc{k}{\mathcal{A}}$ over $\pset{\vec{r},\ell}{}$.
\end{proof}
   
With these intermediate results, we return to the proof of Theorem \ref{thm:non-overloaded-well-behaved-p-box}.
    
\begin{proof}[Proof of {Theorem \ref{thm:non-overloaded-well-behaved-p-box}}]
Consider any regular $\vec{\pp}$-box $\pset{\vec{r},\ell}{}$, which by Proposition \ref{prop:regular-p-box-spline-level} only supports B-splines and is thus trivially not overloaded w.r.t. $\Tbasis{\DomainHierarchy}{}$.
Hence, consider a well-behaved $\vec{\pp}$-box $\pset{\vec{r},\ell}{}$ of depth $d$.
Since the set of THB-splines of $\Tbasis{\DomainHierarchy}{}$ with non-empty support on $\pset{\vec{r},\ell}{}$ coincide with the level-$(\ell+d)$ intermediate THB-splines of  $\Tbasis{\DomainHierarchy_{\ell+d}}{}$ with non-empty support on $\pset{\vec{r},\ell}{}$, it is sufficient to show that the latter set is linearly independent on $\pset{\vec{r},\ell}{}$. This is proven via induction over the level-$k$ intermediate THB-spline bases for $\ell\leq k \leq \ell + d$.

For intermediate basis $\Tbasis{\DomainHierarchy_\ell}{}$, $\pset{\vec{r},\ell}{}$ is an active border $\vec{\pp}$-box. 
Hence, by Lemma \ref{lem:non-overloaded-border-p-box}, the level-$\ell$ intermediate basis $\Tbasis{\DomainHierarchy_\ell}{}$ is linearly independent over $\pset{\vec{r},\ell}{}$, proving the base case.

For the induction step, assume that the intermediate THB-spline basis $\Tbasis{\DomainHierarchy_{k-1}}{}$ is linearly independent over $\pset{\vec{r},\ell}{}$ for $\ell < k \leq \ell +d $. The level-$k$ intermediate THB-spline space is given by
\begin{equation}
\trunc{k}{\Tbasis{\DomainHierarchy_{k-1}}{}} \cup \Bbasis{k}{\act,k}\;.
\end{equation}
By construction, ${\Tbasis{\DomainHierarchy_{k-1}}{}}$ can be written as
\begin{equation}
		\Tbasis{\DomainHierarchy_{k-1}}{} = \trunc{k-1}{\mathcal{A}} \cup \Bbasis{k-1}{\act,k-1}\;,
	\end{equation}
for some set $\mathcal{A}$ (potentially empty). Note that all splines of $\trunc{k-1}{\mathcal{A}}$ are linear combinations of $\Bbasis{k-1}{\nact,k-1}\subset \Bbasis{k-1}{\nact,k}$. Hence, when refining to level $k$, only splines from $\Bbasis{k-1}{\act,k-1}$ are omitted to form the set that will be truncated,
\begin{equation}
		 \trunc{k-1}{\mathcal{A}} \;\cup\; \Bbasis{k-1}{\act,k-1}\backslash \Bbasis{k-1}{\act,k}\quad \subset \quad \Bspace{k-1}{\nact,k}
\end{equation}
By Corollary \ref{cor:linear-independence-A-subset-B-k-1-nact-k} the truncation of this set is linearly independent over $\pset{\vec{r},\ell}{}$, showing that $\Tbasis{\DomainHierarchy_{k}}{}$ is linearly independent over $\pset{\vec{r},\ell}{}$; showing the induction step.
\end{proof}

\subsection{Application to B\'ezier projection}
We conclude this section with an application of $\vec{\pp}$-boxes for the local B\'ezier projector for THB-splines, introduced in \cite{lanini_characterization_2024}.
The B\'ezier projector for THB-splines consists of two steps.
The target function $f \in L^2(\Omega)$ is initially projected onto a space of discontinuous splines $\mathbb{V}$ via local $L^2$ projections.
Then, the local projections are smoothened to retrieve a THB-spline $\Pi f\in\Tspace{\DomainHierarchy}{}$.
However, to ensure that $\Pi$ is a projection, the local projections must be performed on non-overloaded domains $D\subset \Omega$ w.r.t. $\Tbasis{\DomainHierarchy}{}$.

For any $\vec{\pp}$-box mesh, in light of Theorem \ref{thm:non-overloaded-well-behaved-p-box}, we collect the well-behaved $\vec{\pp}$-boxes and regular $\vec{\pp}$-boxes in the set $\mathcal{W}$, so that for any $\pset{}{}\in\mathcal{W}$, $\pbox{}{}$ is not overloaded.
Note, the regular and well-behaved $\vec{\pp}$-boxes cover the domain $\Omega$.
Then, we define the spaces
\begin{equation}\label{eq:discont-THB-spline-basis}
\mathbb{V}_{\pset{}{}} := \lspan{}\left\{ \,\Tspline{j}{}\vert_{\Omega^{\pset{}{}}}\,:\,\forall\,\Tspline{j}{}\in \Tbasis{\DomainHierarchy}{}\,\right\}\;.
\end{equation}
The initial projections of $f\in L^2(\Omega)$ over a regular/well-behaved $\vec{\pp}$-box ${\pset{}{}}\in\mathcal{W}$ is defined by
\begin{equation}\label{eq:THB-spline-macroelementwise-projection}
    \begin{split}
        \sum_{\mathclap{j\in \piset{\pset{}{}}{}{{\Tbasis{\DomainHierarchy}{}}}}} \Coeff{j}{\pset{}{}} \Tspline{j}{}\vert_{\Omega^{\pset{}{}}} := \argmin_{\Spline{h}{\pset{}{}}\in {\mathbb{V}_{\pset{}{}}}} \left\Vert \Spline{h}{\pset{}{}} - f \right\Vert_{L^2(\Omega^{\pset{}{}})}\;,\;
        \piset{\pset{}{}}{}{{\Tbasis{\DomainHierarchy}{}}} := \left\{\,j\,:\,\Tspline{j}{}\in\Tbasis{\DomainHierarchy}{}\;,\,\Omega^{\pset{}{}}\cap\supp{\Tspline{j}{}}\neq \emptyset\,\right\}\;.
    \end{split}
\end{equation}
These local projections are smoothed to form the final projection $\Pi f := \sum_{j=1}^N \Coeff{j}{} \Tspline{j}{}$, where
\begin{align}
\Coeff{j}{} &:=\sum_{\mathclap{\pset{}{}\in\peset{}{j}{\mathcal{W}}}} \omega_j^{\pset{}{}} \Coeff{j}{\pset{}{}}\, , & \peset{}{j}{\mathcal{W}}&:=\left\{\,\pset{}{}\in \mathcal{W}\,:\,\Omega^{\pset{}{}}\cap\supp{\Tspline{j}{}}\neq\emptyset\,\right\}\;,
\end{align}
with averaging constants
\begin{equation}
\omega^{\pset{}{}}_j := \frac{\int_{\Omega^{\pset{}{}}}\Tspline{j}{}dx}{\int_{\Omega}\Tspline{j}{}dx}\;.
\end{equation}
For any $\Spline{}{}=\sum_j\Coeff{j}{}\Tspline{j}{}\in\Tspace{\DomainHierarchy}{}$, the resulting $\Coeff{j}{E}$ from \eqref{eq:THB-spline-macroelementwise-projection} will coincide with $\Coeff{j}{}$.
These coefficients are also preserved by averaging, from which we conclude that $\Pi$ is a projector.
In \cite{lanini_characterization_2024}, it was shown that this projector attains optimal local convergence rates.
This result depends on the support extension, which represents the domain of dependence of the projection $\Pi f$ over an element $\Omega^e$.
\begin{definition}
For a element ${\vec{e}_{\ell}} \in \mesh{\DomainHierarchy}{}$ that is contained in a $\vec{\pp}$-box $\pset{}{}\in\mathcal{W}$, $\Omega^{\vec{e}_{\ell}}\subset \pbox{}{}$, its support extension
is given as:
    \begin{eqnarray}
        \widetilde{\Omega}^{\vec{e}_{\ell}} := \;
            \bigcup_{\mathclap{
        	\substack{
                \pset{}{'} \in \peset{}{j}{\mathcal{W}}\\
        		j \in \piset{\pset{{\vec{e}_{\ell}}}{}}{}{\Tbasis{\DomainHierarchy}{}}
        		}}}\;
         \closure{\Omega^{\pset{}{'}}}\;.
    \end{eqnarray}
\end{definition}
Let $h^{\vec{e}_{\ell}}$ be the mesh size of element $\Omega^{\vec{e}_{\ell}}$. 
Then, on the element $\Omega^{\vec{e}_{\ell}}$, the following Theorem gives a local error estimate for $\Pi f$.
\begin{theorem}
\label{thm:LocalTHBsplineProjEst}
    For ${\vec{e}_{\ell}}\in \mesh{\DomainHierarchy}{}$, $0\leq k \leq m \leq \min(\vec{\pp})+1$ and $f\in H^m(\widetilde{\Omega}^{{\vec{e}_{\ell}}})$:
    \begin{eqnarray}
        \left\vert f - \Pi f \right\vert_{H^k({\Omega}^{{\vec{e}_{\ell}}})} \leq C (h^{\vec{e}_{\ell}})^{m-k}\left\vert f \right\vert_{H^m(\text{cube}(\widetilde{\Omega}^{{\vec{e}_{\ell}}}))}\;,
    \end{eqnarray}
	where $C$ is a constant independent of the mesh size $h^{\vec{e}_{\ell}}$ and $\text{cube}(\widetilde{\Omega}^{{\vec{e}_{\ell}}})$ is the smallest hyper-cube that contains $\widetilde{\Omega}^{{\vec{e}_{\ell}}}$.
\end{theorem}

\begin{remark}
    We claim that the above \Bezier{} projection on macro-elements is local in the following sense.
    The size of the $L^2$-projection problems in \eqref{eq:THB-spline-macroelementwise-projection} is bounded from above by $C\prod_{i=1}^n(p^i+1)$, where $C$ is a constant dependent on the admissibility class.
    While this bound grows with growing degree, it is still independent of the total number of mesh elements/degrees of freedom.
    Furthermore, the following are some special refinement/mesh configurations.
    \begin{itemize}
        \item In the case of a {coarse, single-level mesh with only one $\vec{p}$-box}, i.e., only $\prod_{i=1}^n(p^i+1)$ elements, the local projection proposed here becomes global.
        However, in this case, one is simply working with B-splines, so the original B\'ezier projector can be used to recover locality -- this coincides with the choice $\vec{q} = \vec{1}$.
        \item In the case of a {single-level mesh with multiple $\vec{p}$-boxes}, one is again working with B-splines so any choice $\vec{q}$ will lead to a valid B\'ezier projector, including $\vec{q} = \vec{1}$.
        In particular, our choice of $\vec{q} = \vec{p}$ still leads to multiple local problems over each $\vec{p}$-box.
        \item In all other cases of {THB-splines with multiple active levels}, locality is recovered in the sense described above.
    \end{itemize}
\end{remark}

\section{Adaptive structure-preserving methods via macro-element refinement}\label{sec:struc-pres-methods}
The second kind of $\vec{\qq}$-box, where $\vec{\qq}$ is chosen to be one larger than $\vec{\pp}$, has practical applications in structure-preserving methods.
Structure-preservation provides a framework for creating robust discretizations of mixed formulations of PDEs such as the Navier-Stokes equations, Maxwell equations, and various other equations.
To ensure that the discretization is stable and convergent, the finite-element spaces must mimic the geometric and topological structure that is present in the continuous setting and encoded in the PDE.
This is particularly important for capturing the essential features of the solution, such as conservation laws and symmetries.
For the PDEs mentioned above, the relevant structure is encoded in the so-called de Rham Hilbert complex. For example, in two dimensions, the complex can be written down as:
\begin{equation}\label{eq:cont-derham}
    H^1(\Omega)\xrightarrow[]{\mathrm{rot}}H(\mathrm{div};\Omega)\xrightarrow[]{\mathrm{div}} L^2(\Omega)\;,
\end{equation}
and in three dimensions as:
\begin{equation}\label{eq:cont-derham-3d}
    H^1(\Omega)\xrightarrow[]{\mathrm{grad}}H(\mathrm{curl};\Omega)\xrightarrow[]{\mathrm{curl}} H(\mathrm{div};\Omega) \xrightarrow{\mathrm{div}} L^2(\Omega)\;.
\end{equation}
For contractible domains, such as the unit squares/cubes studied in this work, this complex is exact.
For example, in two dimensions, this means that $\text{image}(\mathrm{rot}) = \text{kernel}(\mathrm{div})$ for the complex \eqref{eq:cont-derham}; in higher dimensions, the same relations between the images and kernels of successive differential operators would hold.
For a non-exact de Rham complex, the kernel space would be larger, as it contains additional functions known as harmonics.

In order to discretize the PDEs associated to the de Rham complex (e.g., fluid flows, electromagnetics), we need to construct discrete spaces that form a discrete de Rham complex that mimics the structure of the continuous complex.
In particular, we want to preserve the relationships between the images and kernels of the differential operators.
For instance, for discretizing the two-dimensional complex in \eqref{eq:cont-derham} using THB-splines, we choose the following THB-spline spaces:
\begin{align}\label{eq:discrete-spaces-derham-2d}
    \mathbb{V}_h^0 :=\Tspace{\DomainHierarchy,[\pp^0,\pp^1]}{}\;,\;\;
    \mathbb{V}_h^1 :=\Tspace{\DomainHierarchy,[\pp^0,\pp^1-1]}{}\times \Tspace{\DomainHierarchy,[\pp^0-1,\pp^1]}{}\;,\;\;
    \mathbb{V}_h^2 :=\Tspace{\DomainHierarchy,[\pp^0-1,\pp^1-1]}{}\;,
\end{align}
and for the three dimensional complex in \eqref{eq:cont-derham-3d}, we choose the THB-spline spaces:
\begin{align}\label{eq:discrete-spaces-derham-3d}
    \mathbb{V}_h^0 :=\Tspace{\DomainHierarchy,[\pp^0,\pp^1,\pp^2]}{}\;,\;\;
    \mathbb{V}_h^1 :=\Tspace{\DomainHierarchy,[\pp^0-1,\pp^1,\pp^2]}{}\times \Tspace{\DomainHierarchy,[\pp^0,\pp^1-1,\pp^2]}{}\times \Tspace{\DomainHierarchy,[\pp^0,\pp^1,\pp^2-1]}{}\;,\\
    \mathbb{V}_h^2 :=\Tspace{\DomainHierarchy,[\pp^0,\pp^1-1,\pp^2-1]}{}\times \Tspace{\DomainHierarchy,[\pp^0-1,\pp^1,\pp^2-1]}{}\times \Tspace{\DomainHierarchy,[\pp^0-1,\pp^1-1,\pp^2]}{}\;,\;\;
    \mathbb{V}_h^3 :=\Tspace{\DomainHierarchy,[\pp^0-1,\pp^1-1,\pp^2-1]}{}\;.
\end{align}
See \cite{evans2020hierarchical,shepherd_locally-verifiable_2024} for the arbitrary dimensional case.
However, as shown in \cite{evans2020hierarchical,shepherd_locally-verifiable_2024}, the resulting discrete complex is not necessarily exact for arbitrary domain hierarchies.
In this section, we prove our main result which is stated in Theorem \ref{thm:exact-discrete-complex}.

\begin{theorem}\label{thm:exact-discrete-complex}
    Consider a domain hierarchy generated on an $n$-dimensional cube $\Omega \subset \RR^n$ by refinements of ($\vec{\pp}+\vec{1}$)-boxes. Then, the corresponding THB-spline de Rham complex is exact.
\end{theorem}

The proof of this theorem follows from the following results where we verify that the assumptions of \cite{shepherd_locally-verifiable_2024} are satisfied by domain hierarchies generated by ($\vec{\pp}+\vec{1}$)-box refinements.
The assumptions in \cite{shepherd_locally-verifiable_2024} are sufficient conditions for the exactness of the THB-spline de Rham complex.
\cite[Assumptions 1-2]{shepherd_locally-verifiable_2024} are easily verified by the construction of the THB-splines considered in this work, and the following discussion thus focuses only on \cite[Assumption 3]{shepherd_locally-verifiable_2024}.
This result depends on the concept of shortest-chains and grids.
We first define these concepts, and then show that \cite[Assumption 3]{shepherd_locally-verifiable_2024} is satisfied in Lemmas \ref{lem:ass-3a} and \ref{lem:ass-3b} below.

\begin{definition}
    Between two vectors $\vec{s}_1,\vec{s}_2\in\mathbb{N}^n$, we call $\{\vec{r}_i\in\mathbb{N}^n : i = 0, \dots, k\}$ a chain when
    $\vec{r}_0 = \vec{s}_1, \vec{r}_k = \vec{s}_2$, only a single component of $\delta \vec{r}_i := \vec{r}_i - \vec{r}_{i-1}$ is non-zero and has unit magnitude, and $k = \sum_{i=1}^{k} |\delta \vec{r}_i|_1$.
    When $k = \sum_{j=1}^n |s_1^j - s_2^j|$, we call it the shortest chain between $\vec{s}_1$ and $\vec{s}_2$.
\end{definition}
\begin{definition}
 Given bounding vectors $\vec{L},\vec{U}\in\mathbb{N}^n$ where $L^j \leq U^j$ for all $j$, define a grid $\mathcal{G}$ as:
\begin{equation}
    \mathcal{G} := \left\{~\vec{t}\in\mathbb{N}^n~:~L^j\leq t^j\leq U^j\;,\;j = 1, \dots, n~\right\}\;.
\end{equation}
\end{definition}
\begin{lemma}\label{lem:grid-shortest-chain}
    For any two grids $\mathcal{G}_1,\mathcal{G}_2$ that overlap (i.e., $\exists \vec{t} \in \mathcal{G}_1 \cap \mathcal{G}_2$),  there exists a shortest chain $\{\vec{r}_0, \cdots, \vec{r}_k\}$ between any two $\vec{s}_1,\vec{s}_2\in\mathcal{G}_1\cup\mathcal{G}_2$ where $\vec{r}_i\in \mathcal{G}_1\cup\mathcal{G}_2$.
    Moreover, if $\vec{s}_1\in\mathcal{G}_1$ and $\vec{s}_2\in\mathcal{G}_2$, there exists an index $\alpha$ such that $\vec{r}_i \in\mathcal{G}_1$ if $i\leq \alpha$ and $\vec{r}_i\in\mathcal{G}_2$ if $\alpha \leq i$. 
\end{lemma}
\begin{proof}
    Let the grids $\mathcal{G}_i$ be constructed from the bounding vector $\vec{L}_i,\vec{U}_i$, where we assume w.l.o.g. that $\mathcal{G}_1 \neq \mathcal{G}_2$ and $L_1^j < L_2^j$ for all $j$.
    Clearly, if $\vec{s}_1,\vec{s}_2\in\mathcal{G}_i$ for some $i$, a trivial shortest chain exists.
    Hence, we consider $\vec{s}_1\in\mathcal{G}_1$, $\vec{s}_2\in\mathcal{G}_2$.
    Define $\epsilon^j := \sgn{s_2^j - s_1^j}$ and $J := \sum_{j=1}^n J^j$ where:
    \begin{equation*}
        J^j := \begin{dcases*}
            \min{\left(U_1^j-s_1^j,s_2^j-s_1^j\right)}\;, &if $\epsilon^j = 1$\;,\\
            s_1^j - s_2^j\;, & else\;.
        \end{dcases*}
    \end{equation*}
    Then, consider the index vector chain $\{\vec{r}_0, \dots, \vec{r}_J\}$ defined as:
    \begin{equation}
        \vec{r}_i = \begin{cases}
            (s_1^1 + \epsilon^1 i,\; s_1^2,\;\dots,\; s_1^n)\;, &\text{for }0 \leq i \leq J^1\;,\\
            (s_1^1 + \epsilon^1 J^1,\;s_1^2+\epsilon^2(i - J^1),\;\dots,\; s_1^n)\;, &\text{for }J^1 < i \leq J^1 + J^2\;,\\
            ~~\vdots\\
            (s_1^1 + \epsilon^1 J^1,\; s_1^2 + \epsilon^2 J^2,\;\dots,\; s_1^n + \epsilon^n(i - J + J^n))\;, &\text{for }J - J^n < i \leq J\;.
        \end{cases}
    \end{equation}
    If $\vec{r}_J = \vec{s}_2$, then this is the desired shortest chain.
    Else, we note that $\vec{r}_J\in\mathcal{G}_2$ since $L_2^j \leq \min\left( U_1^j, s_2^j \right) = \vec{r}_J^j$.
    Thus, we can trivially extend the chain $\{\vec{r}_0, \dots, \vec{r}_J\}$ to form a shortest chain between $\vec{s}_1$ and $\vec{s}_2$ that passes through $\vec{r}_J$.
    Moreover, $\alpha = J$ satisfies the claim.
\end{proof}
\begin{lemma}[Assumption 3a from \cite{shepherd_locally-verifiable_2024}]\label{lem:ass-3a}
    Consider any two level-$\ell$ B-splines $\Bspline{\vec{s}_1,\ell}{},\Bspline{\vec{s}_2,\ell}{}\in \Bbasis{\ell}{\act,\ell+1}$, so that there exists $\vec{r}',k_0$,
\begin{subequations}\label{eq:line-intersection}
    \begin{align}
    \closure{\supp{\Bspline{\vec{s}_1,\ell}{}}} \cap \closure{\supp{\Bspline{\vec{s}_2,\ell}{}}} \supseteq \times_{k=1}^n I^k\;,\\
    I^k:= \begin{cases}
        \left(\widehat{\xi}^k_{r'^k,\ell+1},\widehat{\xi}^k_{r'^k+1,\ell+1}\right)\;, &k\neq k_0\;,\\
        \left\{ \widehat{\xi}^k_{r'^k,\ell+1}\right\}\;, & k = k_0\;.
    \end{cases}
\end{align}
\end{subequations}
Then, these exist a shortest chain $\{\vec{r}_0,\vec{r}_1,\dots\}$ between $\vec{s}_1,\vec{s}_2$, such that $\Bspline{\vec{r}_i,\ell}{}\in \Bbasis{\ell}{\act,\ell+1}$ for all $i$.
\end{lemma}
As a small remark, \eqref{eq:line-intersection} is a slightly stronger condition than what is required in \cite{shepherd_locally-verifiable_2024}.
Nonetheless, this stronger assumption is easier to prove for ($\vec{\pp}+\vec{1}$)-boxes.

\begin{proof}
Without loss of generality, we assume that $k_0 = n$. 
Then we define
\begin{align*}
    \mathcal{G}_i^{\pset{}{}} &:= \left\{~\vec{t}~:~\pset{\vec{t},\ell}{}\in \meshPbox{\ell}{\vec{\pp}+\vec{1}}~:~\supp{\Bspline{\vec{s}_i,\ell}{}}\cap \pbox{\vec{t},\ell}{}\neq\emptyset~ \right\}\;,\\
    \mathcal{G}_i^{\Bspline{}{}} &:= \left\{~\vec{t}~:~\Bspline{\vec{t},\ell}{}\in\Bbasis{\ell}{}\;,~\supp{\Bspline{\vec{t},\ell}{}}\subseteq \bigcup_{\vec{r}\in \mathcal{G}_i^{\pset{}{}}}\closure{\pbox{\vec{r},\ell}{}}~\right\}\;,
\end{align*}
which by the construction of $(\vec{\pp}+\vec{1}$)-boxes are grids.
By lemma \ref{lem:grid-shortest-chain}, the existence of the shortest chain is trivial if $\mathcal{G}_1^{\Bspline{}{}}$ overlaps $\mathcal{G}_2^{\Bspline{}{}}$.
Hence, assume that $\cap_{i=1}^2\mathcal{G}_i^{\Bspline{}{}} = \emptyset$, and thus $\cap_{i=1}^2\mathcal{G}_i^{\pset{}{}} = \emptyset$.
However, given \eqref{eq:line-intersection}, we must have that $s_2^n = s_1^n + \pp^n + 1$ and $\exists \vec{t_1}\in\mathcal{G}_1^{\pset{}{}}$, so that $(t_1^1,\dots,t_1^n+1)\in\mathcal{G}_2^{\pset{}{}}$.

We exploit this by defining two new grids, $\widetilde{\mathcal{G}}_i^{\Bspline{}{}} := \{ \widetilde{\vec{t}} \in \NN^{n-1} : \exists t^n \in \NN\text{~s.t.~}(\widetilde{\vec{t}}, t^n) \in \mathcal{G}_i^{\Bspline{}{}} \}$.
These grids overlap and, by Lemma \ref{lem:grid-shortest-chain}, there exists a shortest chain $\{\widetilde{\vec{r}}_0,\cdots,\widetilde{\vec{r}}_{\tilde{k}}\}$ between the pair $\widetilde{\vec{s}}_1 = (s_1^1,\dots s_1^{n-1}) \in \widetilde{\mathcal{G}}_1^{\Bspline{}{}}$ and $\widetilde{\vec{s}}_2=(s_2^1,\dots s_2^{n-1})\in \widetilde{\mathcal{G}}_2^{\Bspline{}{}}$.
Let the index $\alpha$ be so that
$\widetilde{\vec{r}}_i \in \widetilde{\mathcal{G}}_1^{\Bspline{}{}}$ if $i\leq \alpha$ and $\widetilde{\vec{r}}_i\in \widetilde{\mathcal{G}}_2^{\Bspline{}{}}$ if $i\geq \alpha$.
The shortest chain $\{\vec{r}_i~:~i = 0,\dots,\tilde{k}+\pp^n+1\}$ between $\vec{s}_0$ and $\vec{s}_1$ is then easily constructed as
\begin{equation}
    \vec{r}_i = \begin{cases}
        (\widetilde{\vec{r}}_i,s_1^n)\;, & 0 \leq i\leq \alpha\;,\\
        (\widetilde{\vec{r}}_\alpha,s_1^n+i - \alpha)\;, &  \alpha < i \leq \alpha + \pp^n + 1\;,\\
        (\widetilde{\vec{r}}_{i-\pp^n -1 },s_2^n)\;, & \alpha + \pp^n + 1 < i \leq \tilde{k}+\pp^n+1\;.
    \end{cases}
\end{equation}
\end{proof}

\begin{lemma}[Assumption 3b from \cite{shepherd_locally-verifiable_2024}]\label{lem:ass-3b}
Consider a subset $\mathcal{A} \subset \Bbasis{\ell}{\act,\ell+1}$, and a level-$\ell$ or level-$(\ell+1)$ B-spline $\phi$ of spline degree $\tilde{\vec{\pp}}$ where $\tilde{\pp}^j \in \{\pp^j, \pp^j-1\}$ so that $\supp{\phi}\subseteq \Omega_{\ell+1}$ and $\supp{\phi}\subset \cup_{\Bspline{}{}\in\mathcal{A}} \closure{\supp{\Bspline{}{}}}$. 
Then, there exists a $\Bspline{\vec{r},\ell}{}\in \Bbasis{\ell}{\act,\ell+1}$ so that $\supp{\phi}\subset \supp{\Bspline{\vec{r},\ell}{}}$.
Moreover, $\supp{\Bspline{\vec{r},\ell}{}}$ is contained in the smallest axis-aligned bounding box containing $\cup_{\Bspline{}{}\in\mathcal{A}} \supp{\Bspline{}{}}$.
\end{lemma}
\begin{proof}
    Let $P = \times_{j=1}^n (\widehat{\xi}_{L^j,\ell}^j,\widehat{\xi}_{U^j,\ell}^j)$ be the smallest axis-aligned bounding box containing $\cup_{\Bspline{}{}\in\mathcal{A}} \supp{\Bspline{}{}}$.
    Assume w.l.o.g. that $\supp{\phi}\cap \supp{\Bspline{}{}}\neq \emptyset$ for each $\Bspline{}{}\in\mathcal{A}$.
    Define:
    \begin{equation*}
        R(\phi) := \left\{~\vec{t} : \Bspline{\vec{t},\ell}{}\in\Bbasis{\ell}{\act,\ell+1}\text{~and~}\supp{\phi} \subset \supp{\Bspline{\vec{t},\ell}{}} ~\right\}\;.
    \end{equation*}
    For $\vec{r}\in R(\phi)$ and $\supp{\Bspline{\vec{r},\ell}{}}= \times_{j=1}^n (\widehat{\xi}_{e^j,\ell}^j,\widehat{\xi}_{f^j,\ell}^j)$:
    \begin{equation*}
        \supp{\Bspline{\vec{r},\ell}{}} \subset P \Leftrightarrow L^j \leq e^j, f^j\leq U^j\;.
    \end{equation*}
    Assume w.l.o.g. that this condition is violated for $I\subseteq \{1,\dots,n\}$ as:
    \begin{align*}
        j \notin I \Rightarrow L^j \leq e^j,f^j\leq U^j\;,\quad
        \text{and~}
        j \in I \Rightarrow e^j < L^j\;.
    \end{align*}
    We define $\tilde{\vec{r}}$ such that $\supp{\Bspline{\widetilde{\vec{r}},\ell}{}}= \times_{j=1}^n (\widehat{\xi}_{\tilde{e}^j,\ell}^j,\widehat{\xi}_{\tilde{f}^j,\ell}^j)$ where:
    \begin{equation*}
        j \notin I \Rightarrow \tilde{e}^j := e^j\;,\quad
        \text{and~}
        j \in I \Rightarrow 
        \tilde{e}^j := L^j\;.
    \end{equation*}
    Clearly $\supp{\Bspline{\tilde{\vec{r}},\ell}{}} \subset P$ and, moreover, $\supp{\phi}\subset\supp{\Bspline{\tilde{\vec{r}},\ell}{}}$ since:
    \begin{equation*}
        \supp{\phi} \subset \supp{\Bspline{\vec{r},\ell}{}}\cap  P\subset \supp{\Bspline{\tilde{\vec{r}},\ell}{}}\cap P\;.
    \end{equation*}
    So we only need to verify that $\tilde{\vec{r}} \in R(\phi)$.
    We will do so by showing that for any $\pset{\widetilde{\vec{s}},\ell}{}\in\meshPbox{\ell}{\vec{p+1}}$ such that $\pbox{\widetilde{\vec{s}},\ell}{}\cap\supp{\Bspline{\widetilde{\vec{r}},\ell}{}}\neq \emptyset$, we have that $\pset{\widetilde{\vec{s}},\ell}{}\subset \Omega_{\ell+1}$.

    Consider one such $\pset{\widetilde{\vec{s}},\ell}{}\in\meshPbox{\ell}{\vec{p+1}}$.
    If $\supp{\phi} \cap \pbox{\widetilde{\vec{s}},\ell}{} \neq \emptyset$, then we must have $\pset{\widetilde{\vec{s}},\ell}{}\subset \Omega_{\ell+1}$.
    Therefore, let $\supp{\phi} \cap \pbox{\widetilde{\vec{s}},\ell}{} = \emptyset$.
    Since $\supp{\Bspline{\widetilde{\vec{r}},\ell}{}}$ is obtained by translating $\supp{\Bspline{\vec{r},\ell}{}}$ along dimensions $j\in I$ by at most $\pp^j$ elements, there exists $\pset{\vec{s},\ell}{}\in\meshPbox{\ell}{\vec{p+1}}$ such that:
    \begin{equation*}
        \begin{split}
            &\pbox{\vec{s},\ell}{} \subset \Omega_{\ell+1}\;,\quad
            \pbox{\vec{s},\ell}{}\cap\supp{\Bspline{\vec{r},\ell}{}}\neq \emptyset\;,\\
            &j \in I \Rightarrow s^j \in \{\tilde{s}^j - 1,\tilde{s}^j\}\;,\quad j\notin I \Rightarrow s^j = \tilde{s}^j\;,\\
            &\exists B \in \mathcal{A} \text{~such that~} \supp{B}\cap\pbox{\vec{s},\ell}{}\neq \emptyset\;.
        \end{split}
    \end{equation*}
    Let:
    \begin{equation*}
        \supp{B} = \times_{j=1}^n (\widehat{\xi}_{\alpha^j,\ell}^j,\widehat{\xi}_{\beta^j,\ell}^j)\;,\quad
        \pbox{\widetilde{\vec{s}},\ell}{} = \times_{j=1}^n (\widehat{\xi}_{C^j,\ell}^j,\widehat{\xi}_{D^j,\ell}^j)\;.
    \end{equation*}
    Then, for $j \in I$, the following inequalities hold:
    \begin{equation*}
        \begin{split}
            \alpha^j < C^j & \qquad (\text{because~}\supp{B}\cap\supp{\phi} \neq \emptyset ~\wedge~ \supp{\phi}\cap\pbox{\widetilde{\vec{s}},\ell}{} = \emptyset)\;,\\
            \tilde{e}^j \leq \alpha^j < \tilde{f}^j \leq \beta^j & \qquad (\text{because~}\supp{B}\cap\supp{\phi} \neq \emptyset ~\wedge~\supp{B}, \supp{\Bspline{\tilde{\vec{r}},\ell}{}} \subset P)\;.
        \end{split}
    \end{equation*}
    Consequently, we have the following for $j \in I$:
    \begin{equation*}
        (\widehat{\xi}_{\alpha^j,\ell}^j,\widehat{\xi}_{\beta^j,\ell}^j)\cap (\widehat{\xi}_{C^j,\ell}^j,\widehat{\xi}_{D^j,\ell}^j)\neq \emptyset\;.
    \end{equation*}
    On the other hand, the above also holds true for $j \notin I$ since $\supp{B}\cap\pbox{\vec{s},\ell}{}\neq \emptyset$ and $\vec{s}$ and $\widetilde{\vec{s}}$ only differ along dimensions in $I$.
    Thus, we have that $\supp{B}\cap\pbox{\widetilde{\vec{s}},\ell}{}\neq \emptyset$ as desired.
\end{proof}

{
\begin{remark}
    Given a domain $\check{\Omega} \subset \RR^n$ and a sufficiently smooth map, $\vec{F} : \Omega \rightarrow \check{\Omega}$, the THB-spline complex defined on $\Omega$ can be used to define a structure-preserving discretization of the de Rham complex on $\check{\Omega}$; the argument follows from \cite{buffa2011isogeometric}, and we only recall it here for the specific two- and three-dimensional cases in \eqref{eq:cont-derham} and \eqref{eq:cont-derham-3d}.

    Define $\check{\mathbb{V}}_h^k := \{ \phi : \iota^k(\phi) \in \mathbb{V}_h^k \}$, $k = 0, 1, \dots, n$, as the $k$-th space in the desired discrete de Rham complex on $\check{\Omega}$, where $\iota^k$ are appropriate pullback operators.
    In particular, $\iota^0(\phi) := \phi \circ \vec{F}$, $\iota^n(\phi) := \det(D\vec{F})(\phi \circ \vec{F})$, and:
    \begin{itemize}
        \item $n = 2$, complex \eqref{eq:cont-derham}: 
        \begin{equation*}
		\begin{aligned}
			\iota^{1}(\phi) := \det(D\vec{F})(D\vec{F})^{-1}(\phi \circ \vec{F})\;;
		\end{aligned}
	\end{equation*}
        \item $n = 3$, complex \eqref{eq:cont-derham-3d}:
        \begin{equation*}
		\begin{aligned}
			\iota^1(\phi) := (D\vec{F})^{\mathrm{T}}(\phi \circ \vec{F})\;,\;\;
			\iota^2(\phi) := \det(D\vec{F})(D\vec{F})^{-1}(\phi \circ \vec{F})\;.
		\end{aligned}
	\end{equation*}
    \end{itemize}
    Then, as in \cite{buffa2014isogeometric}, the following is an exact discrete de Rham complex:
    \begin{equation*}
        \begin{tikzcd}
            \check{\mathbb{V}}_h^0 \arrow{r}{\check{d}^0} & \check{\mathbb{V}}_h^1 \arrow{r}{\check{d}^1}  & \dots \arrow{r}{\check{d}^{n-1}} & \check{\mathbb{V}}_h^n\;,
        \end{tikzcd}
    \end{equation*}
    where $\check{d}^k$ is the $k$-th differential operator in the sequences \eqref{eq:cont-derham} and \eqref{eq:cont-derham-3d}.
\end{remark}
}

\section{Numerical results}
\label{sec:numerical-results}
We have implemented admissible refinement and coarsening with $\vec{\qq}$-boxes via an open-source package \cite{dijkstra_nutils-pbox_2024} that builds on top of Nutils \cite{van_zwieten_nutils_2022}, adapting the algorithms from \cite{carraturo_suitably_2019}, taking into account \eqref{eq:multilevel-support-pbox}.
This GitLab repository also contains the code used to generate the numerical results.
We validate the optimal convergence rates of Theorem \ref{thm:LocalTHBsplineProjEst} for the \Bezier{} projector; we will refer to it as the $\vec{\pp}$-box \Bezier{} projector.
Secondly, we introduce an adaptive refinement approach for $\vec{\qq}$-boxes.
With this, we compare the $\vec{\pp}$-box \Bezier{} projector to other local THB-spline projectors \cite{giust_local_2020,lanini_characterization_2024}.
Even more, we apply the approach to adaptive approximation of a Poisson problem \cite{bracco_refinement_2018} where the solution contains a singularity.
Lastly, we use the ($\vec{\pp}+\vec{1}$)-boxes to solve the time-dependent incompressible Navier-Stokes equations with a structure-preserving formulation.
Note that all meshes presented are element meshes.

\subsection{Optimal convergence rates}
\begin{figure}[h]
    \centering
    \begin{tikzpicture}
        \node at (0,0) {\includegraphics[width = 6cm]{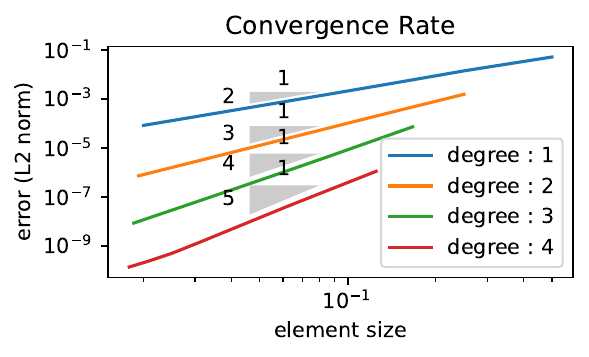}};
        \def\cubeWidth{1.5}
        \def\cubWidthHalf{\cubeWidth/2}
        \def\cubex{4.5}
        \def\cubey{-1}
        \coordinate (O) at (\cubex,\cubey,0);
        \coordinate (B) at (\cubex+\cubeWidth,\cubey+\cubeWidth,-\cubeWidth);
        \draw (O) -- ++ (\cubeWidth,0,0) -- ++ (0,\cubeWidth,0) -- ++(-\cubeWidth,0,0) -- cycle;
        \draw (B) -- ++ (-\cubeWidth,0,0) -- ++ (0,-\cubeWidth,0) -- ++(+\cubeWidth,0,0) -- cycle;
        \draw (O) -- ++ (0,0,-\cubeWidth) -- ++ (0,\cubeWidth,0) -- ++(0,0,\cubeWidth) -- cycle;
        \draw (B) -- ++ (0,0,\cubeWidth) -- ++ (0,-\cubeWidth,0) -- ++(0,0,-\cubeWidth) -- cycle;
        
        \coordinate (O) at (\cubex,\cubey,0);
        \coordinate (B) at (\cubex+\cubWidthHalf,\cubey+\cubWidthHalf,-\cubWidthHalf);
        \draw[dashed] (O) -- ++ (\cubWidthHalf,0,0) -- ++ (0,\cubWidthHalf,0) -- ++(-\cubWidthHalf,0,0) -- cycle;
        \draw[dashed] (B) -- ++ (-\cubWidthHalf,0,0) -- ++ (0,-\cubWidthHalf,0) -- ++(+\cubWidthHalf,0,0) -- cycle;
        \draw[dashed] (O) -- ++ (0,0,-\cubWidthHalf) -- ++ (0,\cubWidthHalf,0) -- ++(0,0,\cubWidthHalf) -- cycle;
        \draw[dashed] (B) -- ++ (0,0,\cubWidthHalf) -- ++ (0,-\cubWidthHalf,0) -- ++(0,0,-\cubWidthHalf) -- cycle;
    \end{tikzpicture}
    \caption{Convergence rate for various spline degrees $\pp$ and the theoretical convergence rates according to Theorem \ref{thm:LocalTHBsplineProjEst}. The same refinement domain is used for all cases for a fair comparison.}
    \label{fig:result-convergence-rate}
\end{figure}

The optimal convergence rate is checked over the unit cube in three dimensions, with $\Omega_1 = [0,1]^3$ and $\Omega_2 = [0,1/2]^3$.
Over this mesh, the target function 
\begin{eqnarray}
    g(\vec{x}) = \prod_{i=1}^3 \sin(\pi x^i)\;,
\end{eqnarray}
is projected for various mesh element sizes $h^i$.
We only consider those $h^i$ so that $\pp h^i$ perfectly divides $1/2$, so that $\Omega_2$ can be represented by $\vec{\pp}$-boxes of level 1 (see Assumption \ref{ass:p-box-refinement-domain}). 
The maximum $L^2$ element error can be seen in Figure \ref{fig:result-convergence-rate}, which perfectly agrees with Theorem \ref{thm:LocalTHBsplineProjEst}.

\subsection{Adaptive refinement: Projection}
\label{sec:NR-adaptive-projection}
For comparing the $\vec{\pp}$-box B\'ezier projector to the local projectors of \cite{lanini_characterization_2024,giust_local_2020}, we adaptively project the following function
\begin{eqnarray}\label{eq:numeric-example-projection-target}
    f(\mathbf{x}) = 1 - \tanh\left( \frac{\|\mathbf{x}\| - 0.3}{0.05\sqrt{2}} \right)\,,\quad \mathbf{x} \in \Omega_1:= [-1,1]^n\,,
\end{eqnarray}
till the $L^\infty$ projection error is below the tolerance $tol = 10^{-4}$, where $n = 2, 3$.

\begin{figure}[h]
    \centering
    \begin{subfigure}{0.3\textwidth}
        \centering
        \includegraphics[width = \textwidth]{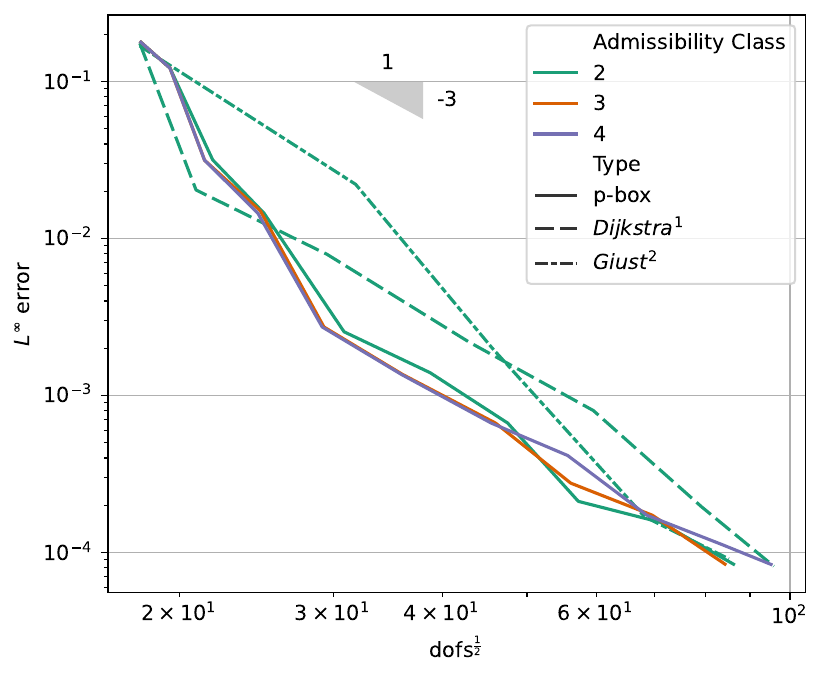}\\
        \includegraphics[width = 0.77\textwidth]{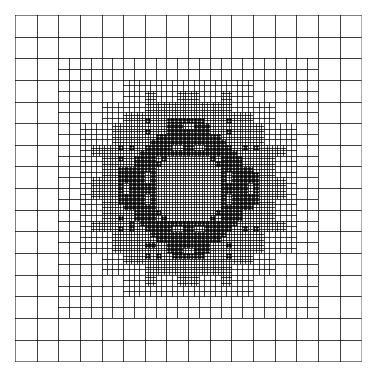}
        \caption{$n = 2; p=2$}\label{fig:result-adaptive-refinement-p2}
    \end{subfigure}
    \begin{subfigure}{0.3\textwidth}
        \centering
        \includegraphics[width = \textwidth]{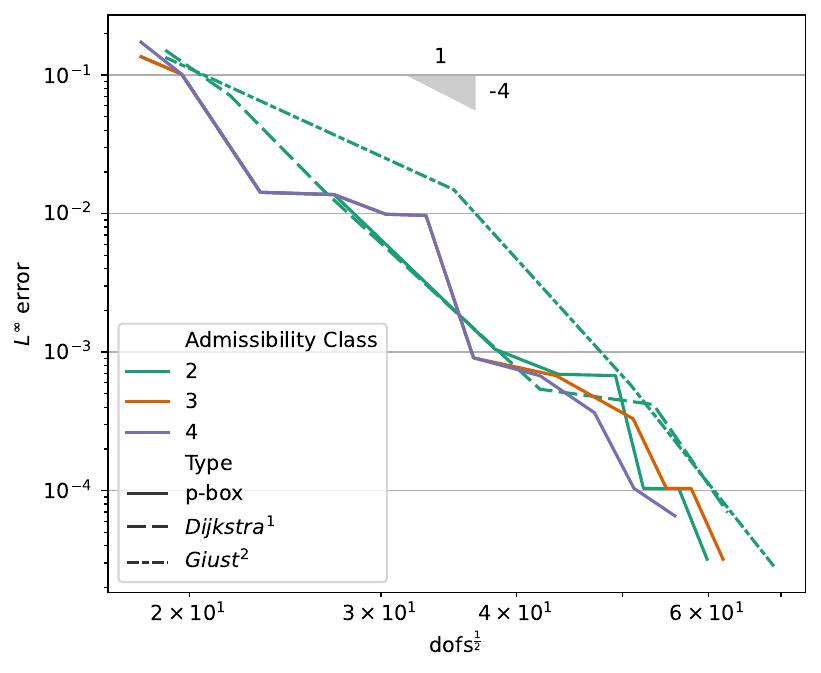}\\
        \includegraphics[width = 0.77\textwidth]{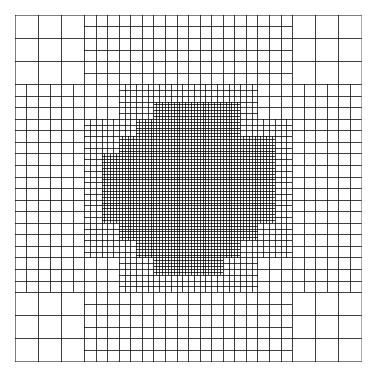}
        \caption{$n = 2; p = 3$}\label{fig:result-adaptive-refinement-p3}
    \end{subfigure}
    \begin{subfigure}{0.3\textwidth}
        \centering
        \includegraphics[width = \textwidth]{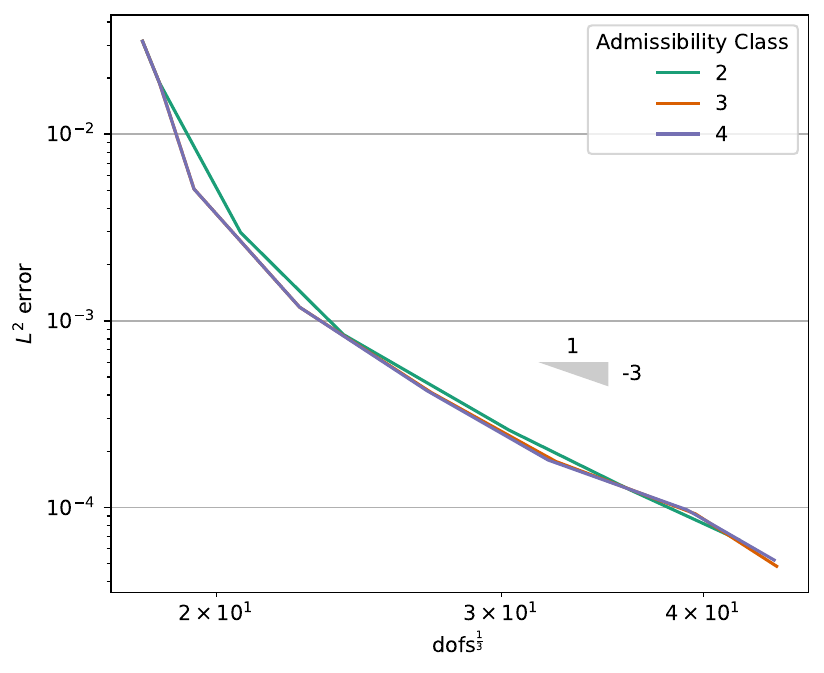}\\
        \includegraphics[width = \textwidth]{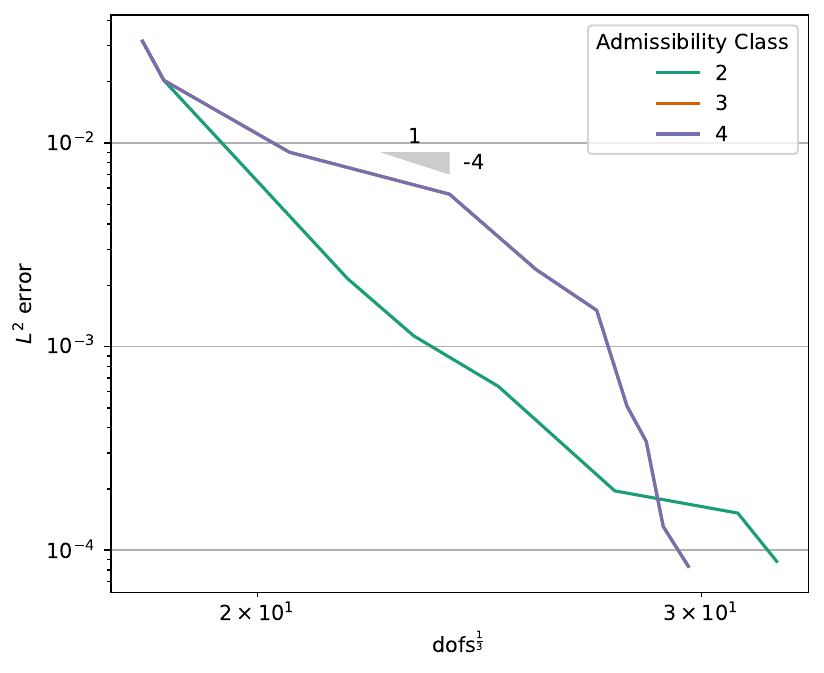}
        \caption{$n=3; p=2$ and $p=3$}\label{fig:result-adaptive-refinement-d3}
    \end{subfigure}
    \caption{
    Error reduction during adaptive approximation of \eqref{eq:numeric-example-projection-target}.
    The $\vec{\pp}$-box B\'ezier projector (solid lines) for admissibility classes $c=2,3,4$ is compared to the projectors of \cite{lanini_characterization_2024} (dashed lines) and \cite{giust_local_2020} (dashed-dotted lines) for admissibility class $c=2$.
    Results are shown for $\vec{\pp}=(2,2)$ (top left) or $\vec{\pp}=(3,3)$ (top center) and dimension $n=2$.
    The final meshes of admissibility class $c=2$ for spline degrees $\vec{\pp}=(2,2)$ (bottom left) and $\vec{\pp}=(3,3)$ (bottom center) are also presented.
    The results obtained for $n=3$ and $\vec{\pp}=(2,2,2)$ (top right) and $\vec{\pp}=(3,3,3)$ (bottom right) are also shown; we omit the 3D meshes since they are hard to visualize.
    Note that in the $n=3$ plots, the plots for $c=3$ and $c=4$ are overlapping.
    }
    \label{fig:result-adaptive-refinement}
\end{figure}

For $n = 2$, this is compared to the results from \cite{lanini_characterization_2024,giust_local_2020}.
In \cite{lanini_characterization_2024,giust_local_2020}, this benchmark problem is solved for $\pp=2,3$, and initial meshes with $16\times 16$ elements, and with a D\"orfler constant $\theta = 0.5$.
For degree $\pp=2$, $\vec{\pp}$-boxes can perfectly reproduce this initial mesh, but for degree $\pp=3$, we take an initial mesh consisting of $5$ $\vec{\pp}$-boxes ($15$ mesh elements).
The resulting $\vec{\pp}$-box meshes for admissibility class 2 are also shown in Figures \ref{fig:result-adaptive-refinement-p2} and \ref{fig:result-adaptive-refinement-p3}.
In Figures \ref{fig:result-adaptive-refinement-p2} and \ref{fig:result-adaptive-refinement-p3}, little difference between $\vec{\pp}$-boxes and the methods of \cite{lanini_characterization_2024,giust_local_2020} can be observed by measuring in the $L^\infty$ norm (for comparison reasons).
This is noteworthy because $\vec{\pp}$-boxes refine multiple elements simultaneously, which could potentially make them less effective at capturing 
finer details. 
However, there is no evidence to suggest that this is the case. Instead, $\vec{\pp}$-boxes provide a straightforward refinement strategy.

\subsection{Adaptive refinement: Poisson}\label{sec:NR-poisson-l-shape}
\begin{figure}[htb]
    \centering
    \centering
    \def\widthFigures{0.3}
    \def\FigureSize{0.9\textwidth}
    \begin{subfigure}[b]{\widthFigures\textwidth}
        \centering
        \includegraphics[width=\FigureSize]{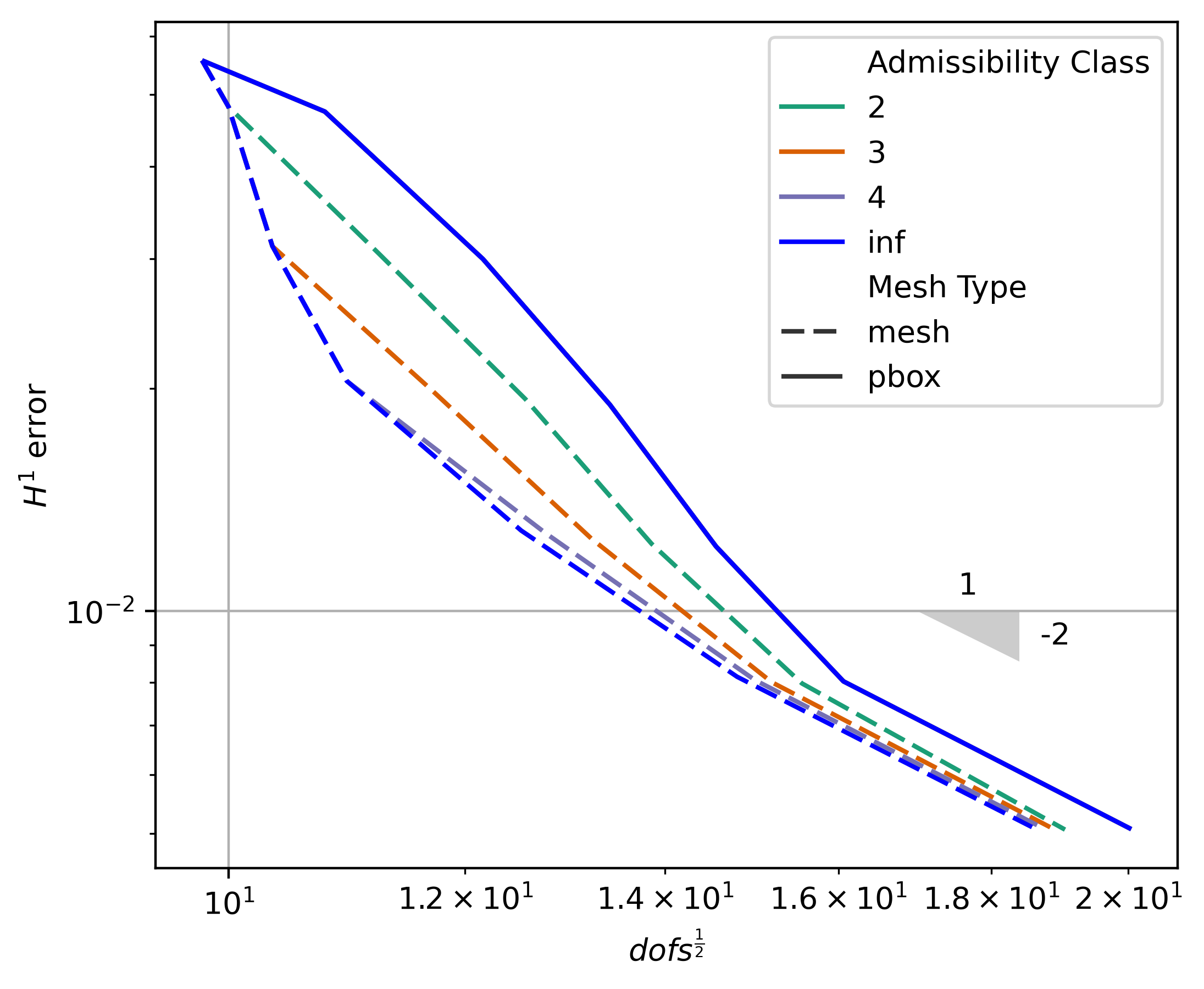}
        \caption{convergence rate for $\pp = 2$}
    \end{subfigure}
    \begin{subfigure}[b]{\widthFigures\textwidth}
        \centering
        \includegraphics[width=\FigureSize]{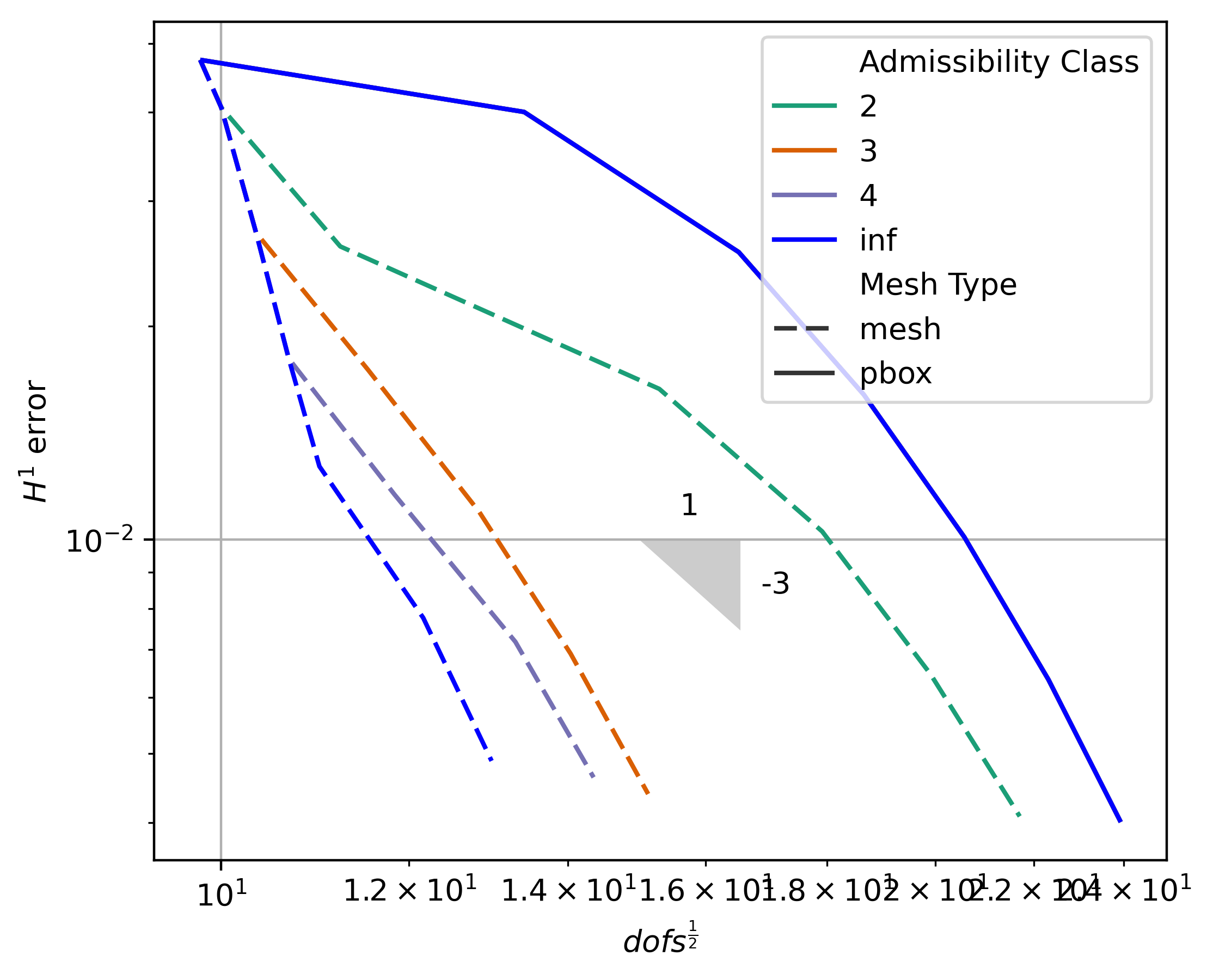}
        \caption{convergence rate for $\pp = 3$}
    \end{subfigure}
    \begin{subfigure}[b]{\widthFigures\textwidth}
        \centering
        \includegraphics[width=\FigureSize]{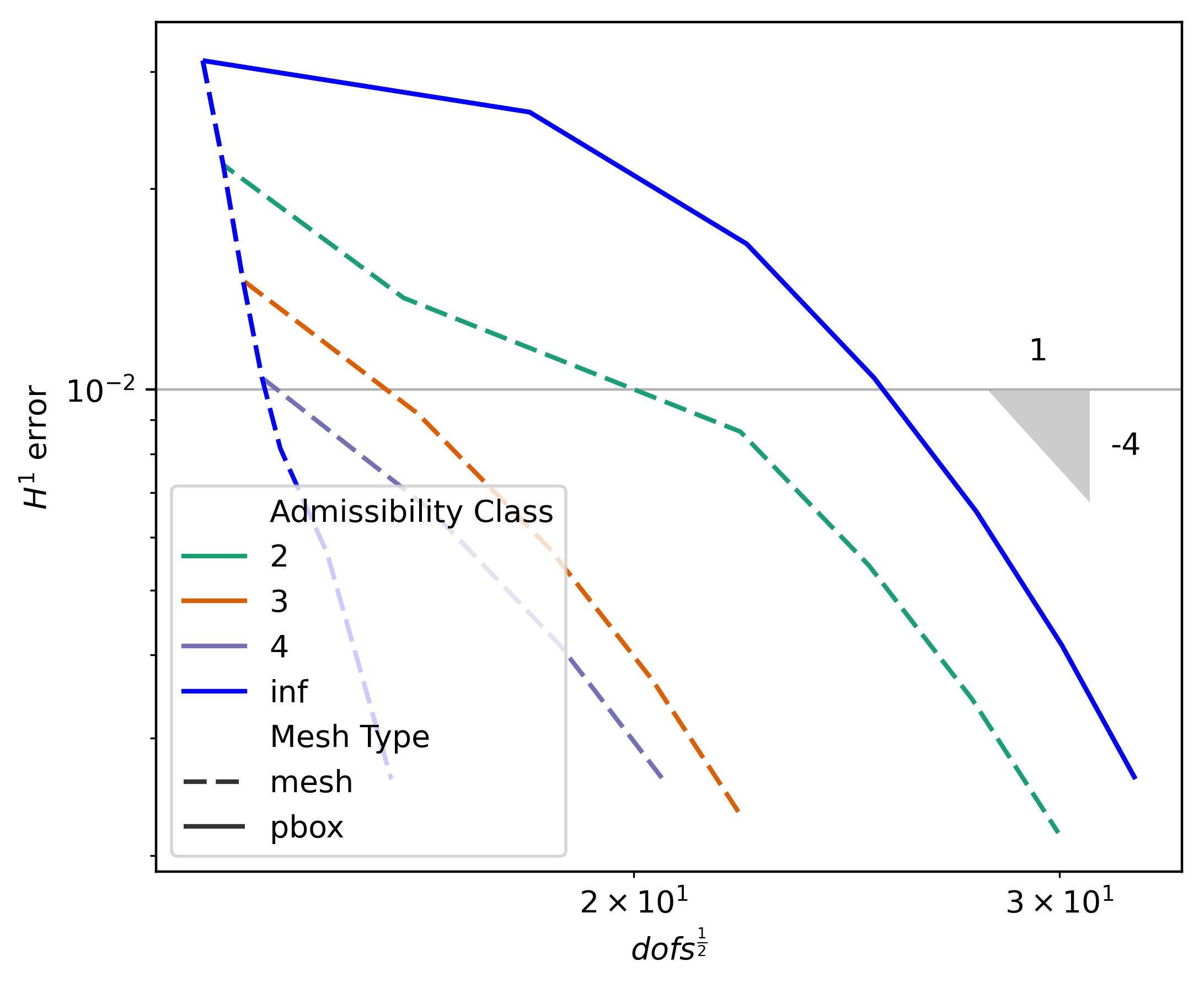}
        \caption{convergence rate for $\pp = 4$}
    \end{subfigure}
    \begin{subfigure}[b]{\widthFigures\textwidth}
        \centering
        \includegraphics[width=0.8\textwidth]{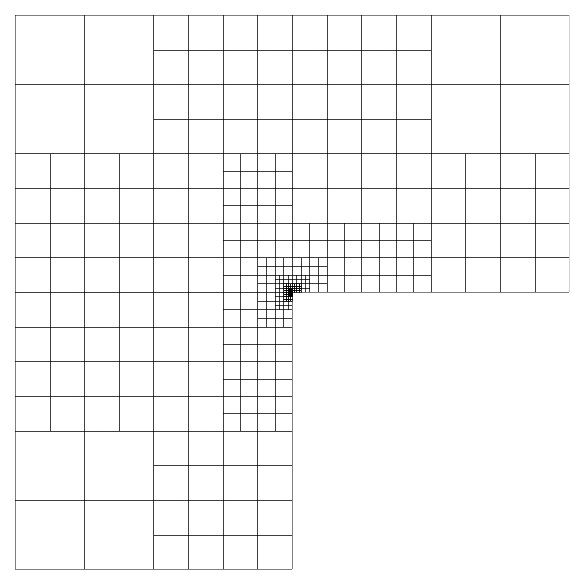}
        \caption{mesh for $\pp = 2, c=2$}
    \end{subfigure}
    \begin{subfigure}[b]{\widthFigures\textwidth}
        \centering
        \includegraphics[width=0.8\textwidth]{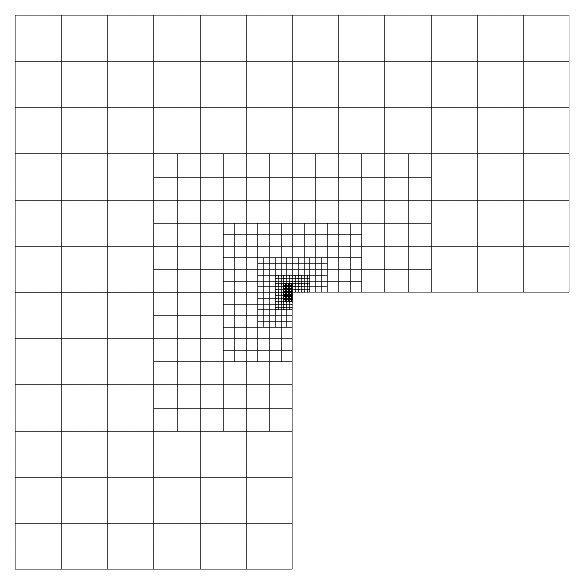}
        \caption{mesh for $\pp = 3, c=2$}
    \end{subfigure}
    \begin{subfigure}[b]{\widthFigures\textwidth}
        \centering
        \includegraphics[width=0.8\textwidth]{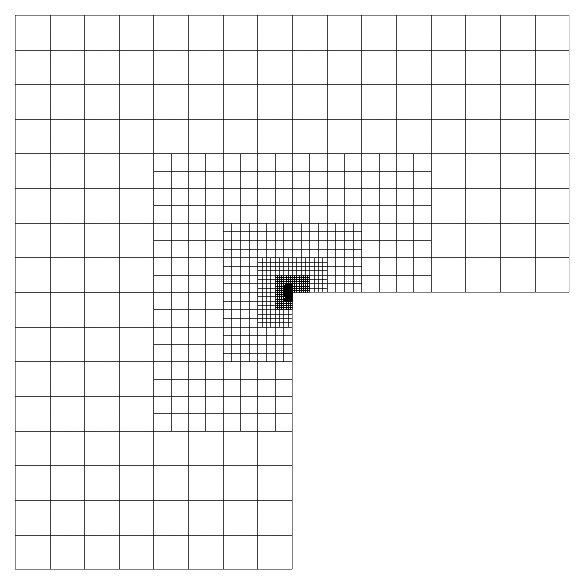}
        \caption{mesh for $\pp = 4, c=2$}
    \end{subfigure}
    \caption{Convergence rates of benchmark 2 for degrees $\vec{\pp}=(2,2), (3,3), (4,4)$.{ The problem is solved for admissibility class $c=2,3,4,\infty$. The $\vec{\pp}$-box mesh results for admissibility classes $c=2,3,4,\infty$ are also shown, but the plots overlap, so only one is visible.}}
    \label{fig:benchmark-2-results}
\end{figure}

The following example is inspired by \cite{bracco_refinement_2018}, in which a Poisson problem is solved where the solution has a singularity.
For this problem, lower admissibility classes were observed to use more degrees of freedom for a given level of accuracy.
We refer the reader to \cite{bracco_refinement_2018} for an in-depth discussion of this method and how it is affected by admissibility classes.
We have implemented the adaptive refinement method for the following problem, where $\Omega = (-1,1)^2 \backslash [0,1]\times[-1,0]$ is built using three patches,
\begin{subequations}
\begin{align}
    \Delta u &= f\;,&\text{ on }\Omega\;,\\
    u&= g\;,&\text{ on }\partial \Omega\;.
\end{align}
\end{subequations}
We discretise the solution $u_h$ with THB-splines that are $C^0$ continuous across the patches, and on each patch, the mesh is given by a $\vec{\pp}$-box mesh.
To investigate the convergence rate, we use the following manufactured solution
\begin{equation}\label{eq:benchmark-2-solution}
    u(r,\theta) := r^{{2}/{3}} \sin\left({{2}\theta/{3}}\right)
\end{equation}
and $g$ is chosen as its restriction to $\partial \Omega$.
We use the following elementwise error estimator,
\begin{eqnarray}
    \epsilon_{\Omega^{\vec{e}_{\ell}}}^2(u_h) = h_{\vec{e}_{\ell}}^2\int_{F(\Omega^{\vec{e}_{\ell}})} |f-\Delta u_h|^2 dV + h_{\vec{e}_{\ell}}\sum_{i\neq j}\int_{\Gamma_{ij}\cap \partial F(\Omega^{\vec{e}_{\ell}})}|\nabla u_h\cdot \vec{n}|^2 dS\;,
\end{eqnarray}
where $h_{\vec{e}_{\ell}}$ is the diameter of element $\Omega^{\vec{e}_{\ell}}$, $F$ maps $\Omega^{\vec{e}_{\ell}}$ to the unit square, and $\Gamma_{ij} := \partial \Omega_i \cap \partial \Omega_j$ is the boundary between the two patches.
For a $\vec{\pp}$-box $\pbox{}{}$, the estimator is given by,
\begin{equation}
    \epsilon_{\pbox{}{}}^2(u_h) := \sum_{\mathclap{\Omega^{\vec{e}_{\ell}}\in\mesh{}{}~:~\Omega^{\vec{e}_{\ell}}\subset\pbox{}{}}}~\epsilon_{\Omega^{\vec{e}_{\ell}}}^2(u_h)\;.
\end{equation}
Lastly, we use D\"orfler marking with constant $\theta = 0.9$.
The convergence rates can be seen in Figure \ref{fig:benchmark-2-results}.
We note two observations.
Firstly, as observed in \cite{bracco_refinement_2018}, imposing an admissibility class worsens the method's pre-asymptotic behaviour.
However, in the case of $\vec{\pp}$-boxes, the various admissibility classes generate identical meshes of admissibility class $c=2$.
This can be explained by Proposition \ref{prop:sufficient-condition-admiss-class-2} and the observation that the singularity is located at a patch corner.
The closest $\vec{\pp}$-box to the corner will always be an active regular $\vec{\pp}$-box; as a result, all well-behaved $\vec{\pp}$-boxes are border $\vec{\pp}$-boxes so that the mesh is of admissibility class $c=2$.
While this seems to present a downside of $\vec{\pp}$-box refinement, we believe it more effectively highlights the advantage of THB-splines at the borders and corners, which $\vec{\pp}$-box refinement is unable to exploit. See Remark \ref{rem:benchmark-2-discussion}.

\begin{figure}[htb]
    \centering
    \begin{subfigure}{0.3\textwidth}
        \centering
        \includegraphics[width=\textwidth]{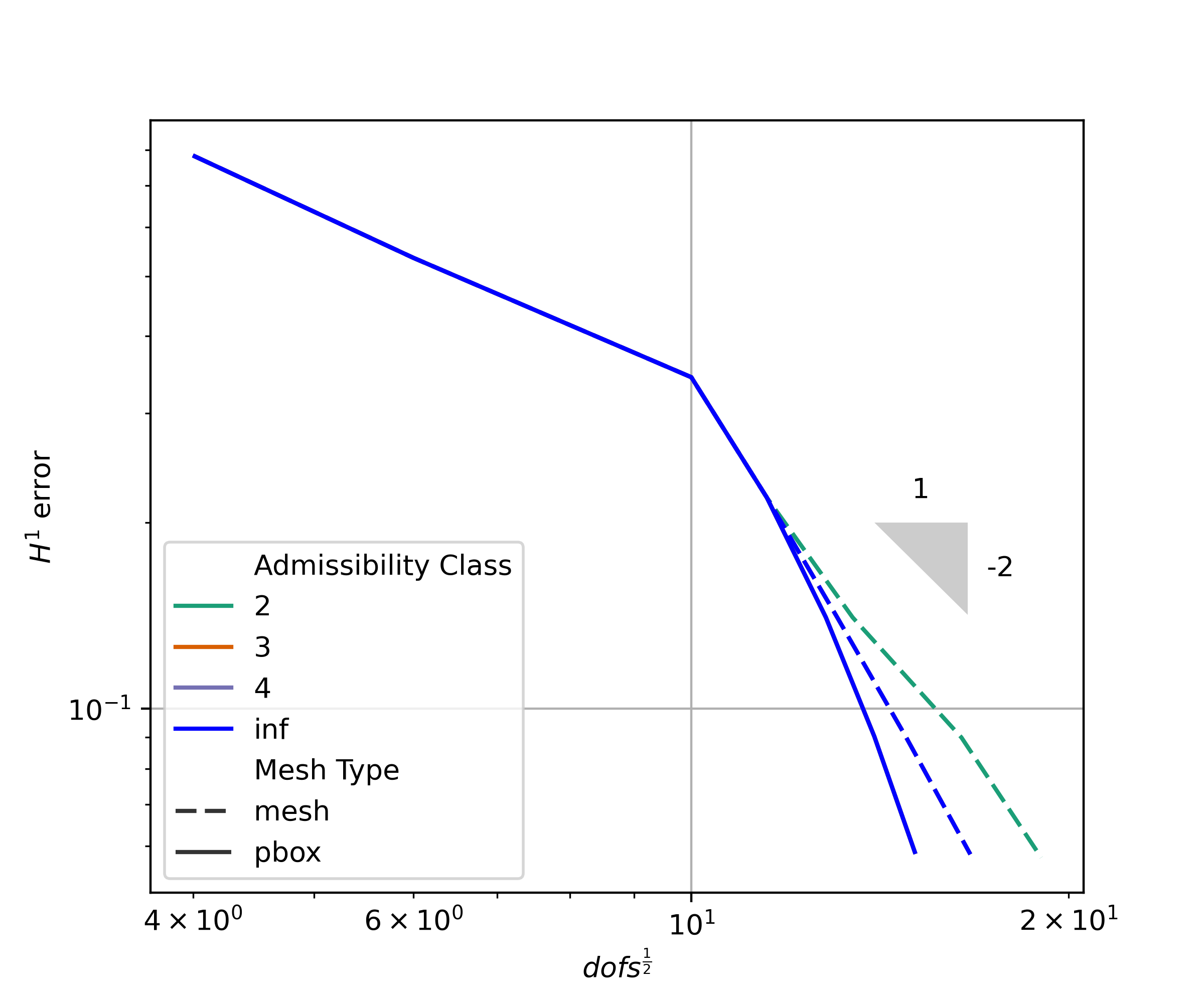}
        \caption{convergence rate for $\pp=2$}
        \label{fig:benchmark-2-remark-results-convergence-p2}
    \end{subfigure}
    \hfill
    \begin{subfigure}{0.3\textwidth}
        \centering
        \includegraphics[width=\textwidth]{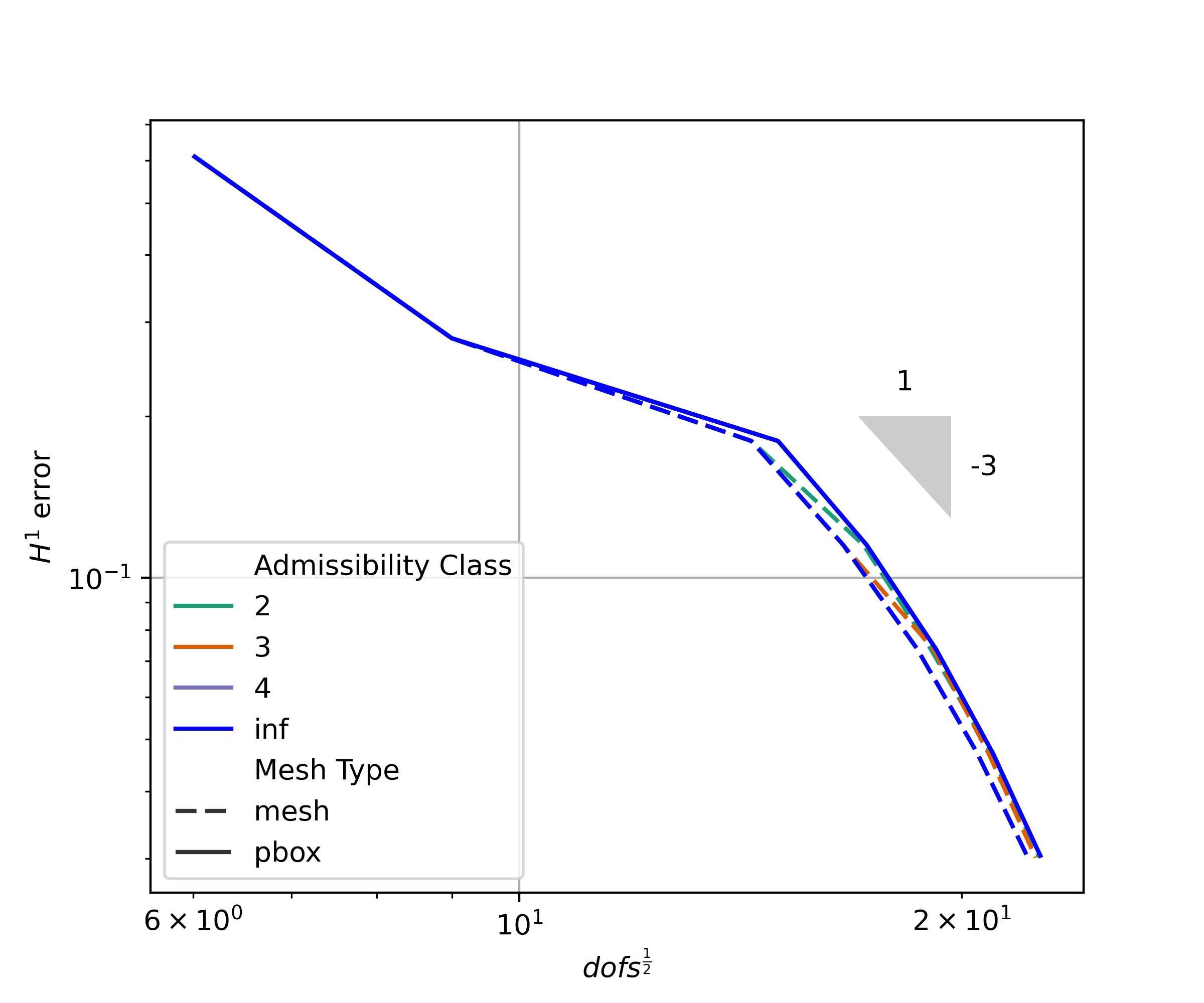}
        \caption{convergence rate for $\pp=3$}
    \end{subfigure}
    \hfill
    \begin{subfigure}{0.3\textwidth}
        \centering
        \includegraphics[width=\textwidth]{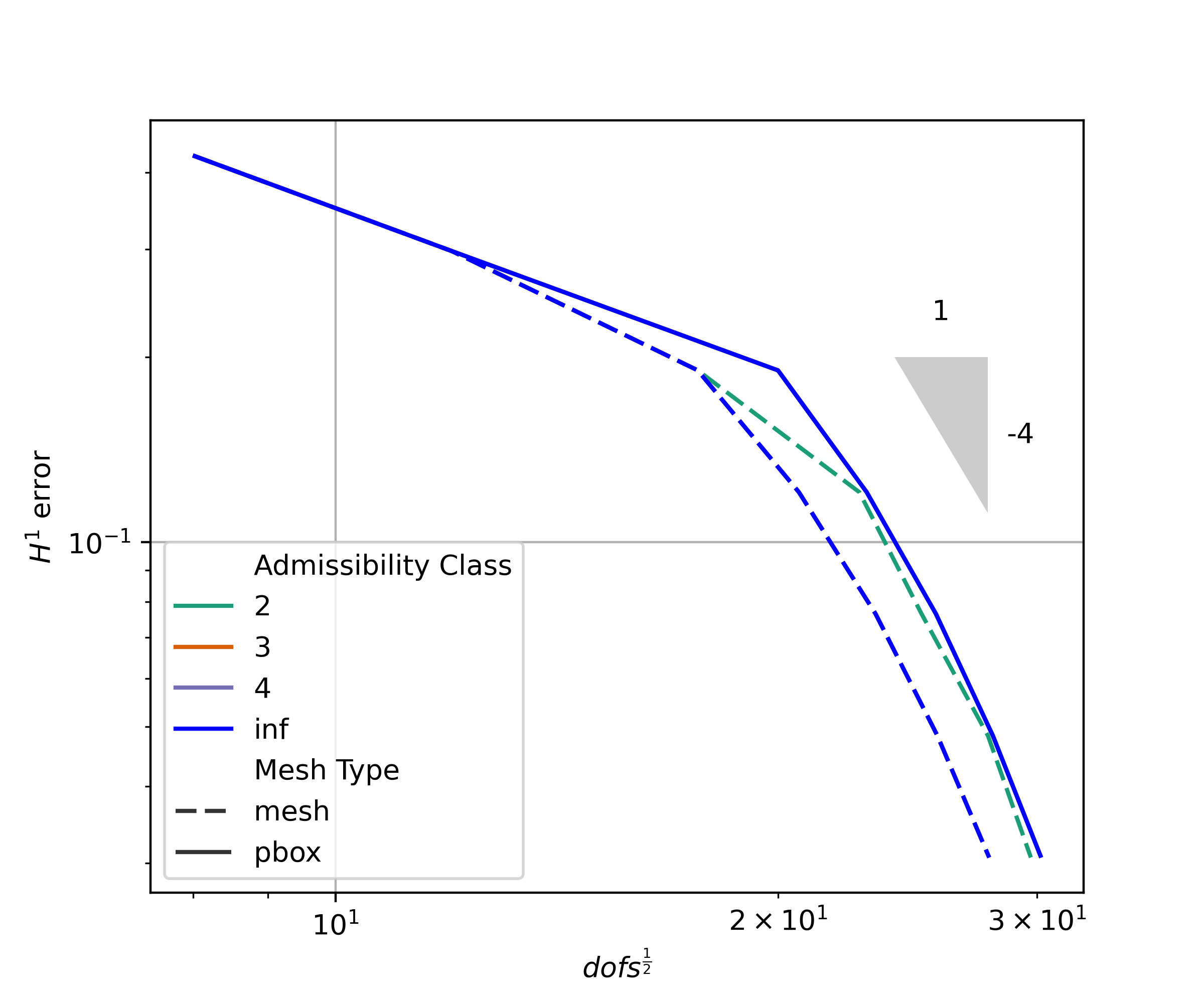}
        \caption{convergence rate for $\pp=4$}
    \end{subfigure}
    \begin{subfigure}{0.3\textwidth}
        \centering
        \includegraphics[width=\textwidth]{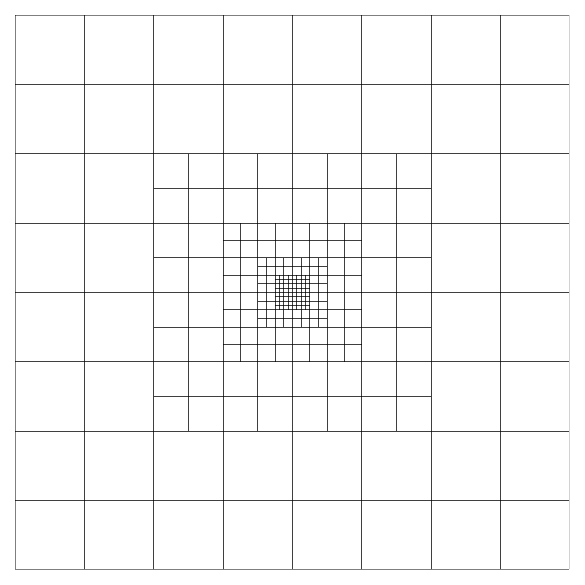}
        \caption{mesh for $\pp=2,c=2$}
    \end{subfigure}
    \hfill
    \begin{subfigure}{0.3\textwidth}
        \centering
        \includegraphics[width=\textwidth]{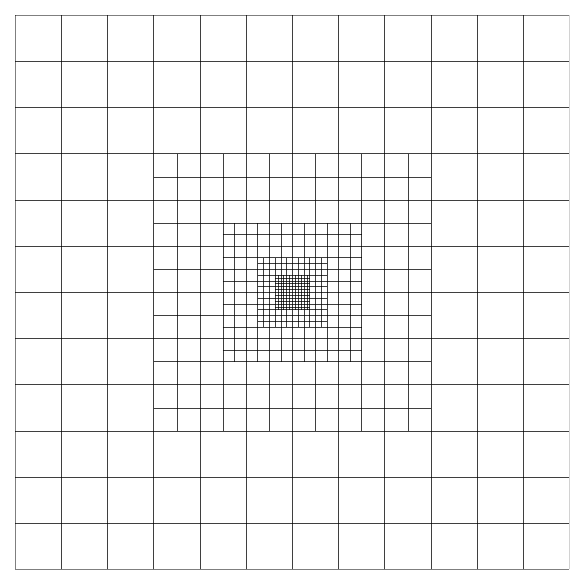}
        \caption{mesh for $\pp=3,c=2$}
    \end{subfigure}
    \hfill
    \begin{subfigure}{0.3\textwidth}
        \centering
        \includegraphics[width=\textwidth]{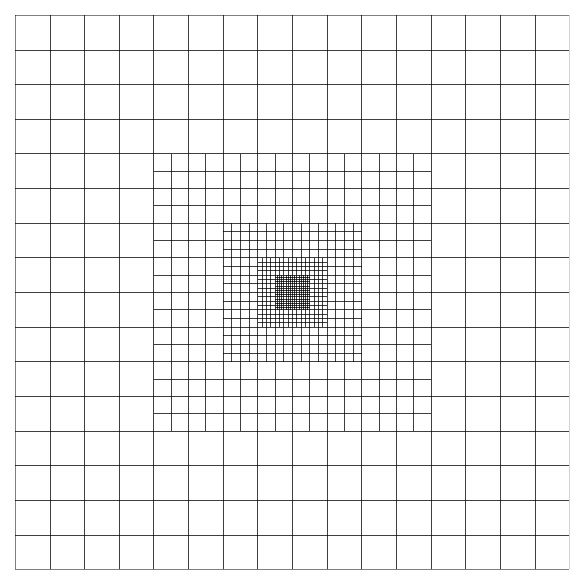}
        \caption{mesh for $\pp=4,c=2$}
    \end{subfigure}
    \caption{The Poisson problem of Section \ref{sec:NR-poisson-l-shape} solved on a square domain $\Omega = [-1,1]^2$, with the manufactured solution $u(r,\theta) = r^{\frac{2}{3}}$. For the element mesh, B-spline supports are refined instead of single mesh elements of Section \ref{sec:NR-poisson-l-shape}. As a result, the error is calculated over spline supports; see \cite{vuong_hierarchical_2011}.
    This problem is solved with admissibility classes $c=2,3,4,\infty$.
    For both the $\vec{\pp}$-box meshes and element meshes, the convergence lines coincide or are incredibly close.
    Notably, for $\pp=2$, the $\vec{\pp}$-box mesh outperforms the element mesh, 
    and for degrees $\pp=3$ and $\pp=4$, the gap between $\vec{\pp}$-box and element meshes of Figure \ref{fig:benchmark-2-results} is reasonable.}
    \label{fig:benchmark-2-remark-results}
\end{figure}

\begin{remark}\label{rem:benchmark-2-discussion}
    In numerical experiments involving a singularity placed in the middle of the domain and performing adaptive projection, the discrepancy observed in Figure \ref{fig:benchmark-2-results} almost vanishes.
    This is because in this case, both $\vec{\pp}$-box refinement and mesh-element refinement require the refinement of B-spline supports, i.e., multiple elements are refined simultaneously.
    The performance of $\vec{\pp}$-box refinement in this setup is only slightly worse, and interestingly, mesh-element refinements for different admissibility classes yield similar results (i.e., no advantage of using very high admissibility classes).
    This suggests that the behaviour observed in Figure \ref{fig:benchmark-2-results} is because the corner singularity ends up being a best-case scenario for mesh-element refinement, as it allows the refinement of a single element, which coincides with the support of corner B-splines.
    {For completeness, the results of this experiment are shown in Figure \ref{fig:benchmark-2-remark-results}. 
    Notably, in Figure \ref{fig:benchmark-2-remark-results-convergence-p2}, for $\pp=2$, the $\vec{\pp}$-box mesh outperforms mesh-element refinement. 
    However, the resulting meshes appear over-refined, suggesting a too large D\"orfler marking parameter for mesh-element refinement (for $\pp=2$).
    To keep the discussion focused, we choose to forego tests that require finding optimal D\"orfler parameters for the two different adaptivity approaches.
    }
\end{remark}

\subsection{Stability of the THB-spline de Rham complex}\label{sec:NR-stability}
{
Stable discretizations of PDEs such as Maxwell's equations and the Stokes equations on contractible domains require discrete de Rham complexes that are not only exact but also stable.
By Theorem \ref{thm:exact-discrete-complex}, any THB-spline complex resulting from a $(\vec{\pp}+\vec{1})$-box refinement is exact.
While a proof of the stability of such complexes is an open problem in isogeometric analysis and beyond the scope of this work, for completeness, we investigate the stability of the THB-spline complexes for different refinement configurations through numerical experiments.
Similar experiments have also been performed in \cite{evans2020hierarchical}.

\subsubsection{Refinement configurations tested}
As in \cite{evans2020hierarchical}, we will test the stability of the complex for $\vec{\pp} = (4,4)$ and the $(\vec{\pp}+\vec{1})$-box meshes given in Figure \ref{fig:spuriousmode-check-meshes}.
For both patterns, the initial mesh consists of $4\times 4$ $(\vec{\pp}+\vec{1})$-boxes, so that the initial mesh is a tensor-product mesh and generates a stable complex by \cite{buffa2011isogeometric}.
From this initial mesh, multiple meshes are generated through subsequent refinements.
The first refinement pattern, shown in Figure \ref{fig:spuriousmode-check-meshes-corner-to-corner}, is corner-to-corner refinement, which is graded towards three points: two opposite corners and the center of the domain.
For $\ell>1$, each subsequent mesh is constructed by considering the previous level-$(\ell-1)$ mesh and refining the four level-$(\ell-1)$ $(\vec{\pp}+\vec{1})$-boxes to the bottom-left and/or the top-right with respect to the three points.
The result is an artificial mesh in which the admissibility class $c$ equals the total number of refinement levels $L$; e.g., $c=5$ for the mesh of Figure \ref{fig:spuriousmode-check-meshes-corner-to-corner}.
The second refinement pattern, shown in Figure \ref{fig:spuriousmode-check-meshes-diagonal}, is refinement along the diagonal.
Here, for $\ell>1$, each subsequent refinement is performed by 
refining the level-$(\ell-1)$ $(\vec{\pp}+\vec{1})$-boxes that make up the diagonal and the first lower and upper diagonals.
However, as diagonal refinement generates meshes with more elements, simulations were only possible for three levels of refinement.
The generalized eigenvalue problems on these meshes are solved using the MultiParEig add-on in MATLAB.

\begin{figure}[htb]
    \centering
    \begin{subfigure}{0.45\textwidth}
        \centering
        \includegraphics[width=0.6\textwidth]{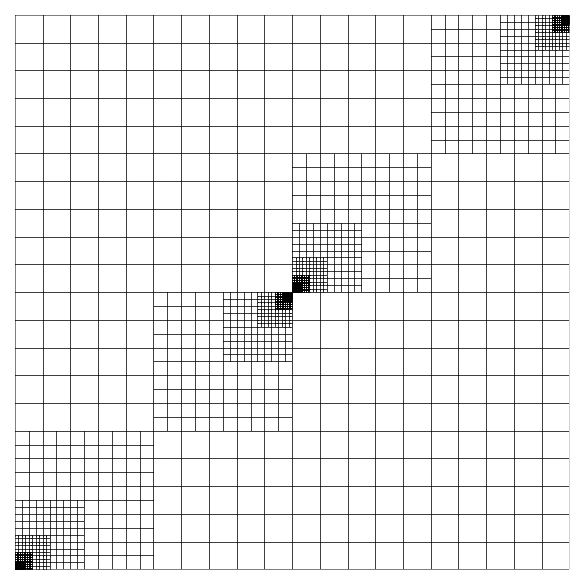}
        \caption{corner-to-corner refinement}
        \label{fig:spuriousmode-check-meshes-corner-to-corner}
    \end{subfigure}
    \hfill
    \begin{subfigure}{0.45\textwidth}
        \centering
        \includegraphics[width=0.6\textwidth]{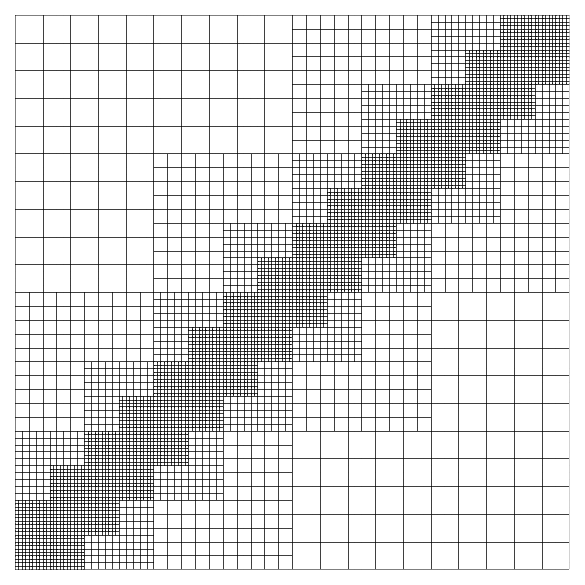}
        \caption{diagonal refinement}
        \label{fig:spuriousmode-check-meshes-diagonal}
    \end{subfigure}
    \caption{The two types of refinement patterns considered in Section \ref{sec:NR-stability}. Figure \ref{fig:spuriousmode-check-meshes-corner-to-corner} shows corner-to-corner refinement, where we investigate both refinement into the corner of the domain and adjacent box refinement, which affects the admissibility class.  It is important to note that this refinement pattern is highly artificial and serves as a worst-case scenario. In Figure \ref{fig:spuriousmode-check-meshes-diagonal}, we consider diagonal refinement, which is a more representative refinement pattern. }
    \label{fig:spuriousmode-check-meshes}
\end{figure}

\subsubsection{Maxwell eigenvalue problem:}
We start by analysing the Maxwell eigenvalue problem on $\Omega =(0,\pi)^2$.
For this problem, the appropriate de Rham complex is given as:
\begin{equation}
    \begin{tikzcd}
        H_0^1(\Omega) \arrow{r}{\mathrm{grad}} & H_0(\mathrm{curl};\Omega) \arrow{r}{\mathrm{curl}} & L_0^2(\Omega)\;,
    \end{tikzcd}
\end{equation}
where $H_0(\mathrm{curl};\Omega) := \{ \vec{v} \in H(\mathrm{curl};\Omega) : \vec{v}\times \vec{n} = 0\text{ on }\partial \Omega \}$.
The corresponding discrete spaces are given by
\begin{equation}
    \begin{split}
        \mathbb{V}_{h,0}^0 &:= \mathbb{T}_{\DomainHierarchy,[p^0,p^1]} \cap H_0^1(\Omega)\;,\\
        \mathbb{V}_{h,0}^{1,*} &:= \mathbb{T}_{\DomainHierarchy,[p^0-1,p^1]} \times \mathbb{T}_{\DomainHierarchy,[p^0,p^1-1]} \cap H_0(\mathrm{curl};\Omega)\;,\\
        \mathbb{V}_{h,0}^{2} &:= \mathbb{T}_{\DomainHierarchy,[p^0-1,p^1-1]} \cap L_0^2(\Omega)\;.
    \end{split}
\end{equation}
Then, the discrete problem formulations reads: find $\vec{u}_h \in \mathbb{V}_{h,0}^{1,*}$ and $\lambda_h \in\RR$, such that
\begin{equation}
    \left( \mathrm{curl} \vec{u}_h,\mathrm{curl} \vec{v}_h \right) = \lambda_h \left( \vec{u}_h,\vec{v}_h\right)\;,\quad \forall \vec{v}_h\in \mathbb{V}_{h,0}^{1,*}\;.
\end{equation}
The analytical non-zero eigenvalues for the continuous formulation of the same problem are given by $\lambda = m_0^2 + m_1^2$ for $m_0,m_1\in\mathbb{N}$.
In addition, the number of discrete zero eigenvalues should equal $\mathrm{dim}(\mathbb{V}_{h,0}^0)$.

\begin{figure}[htb]
    \centering
    \begin{subfigure}{0.3\textwidth}
        \centering
        \includegraphics[width=\textwidth]{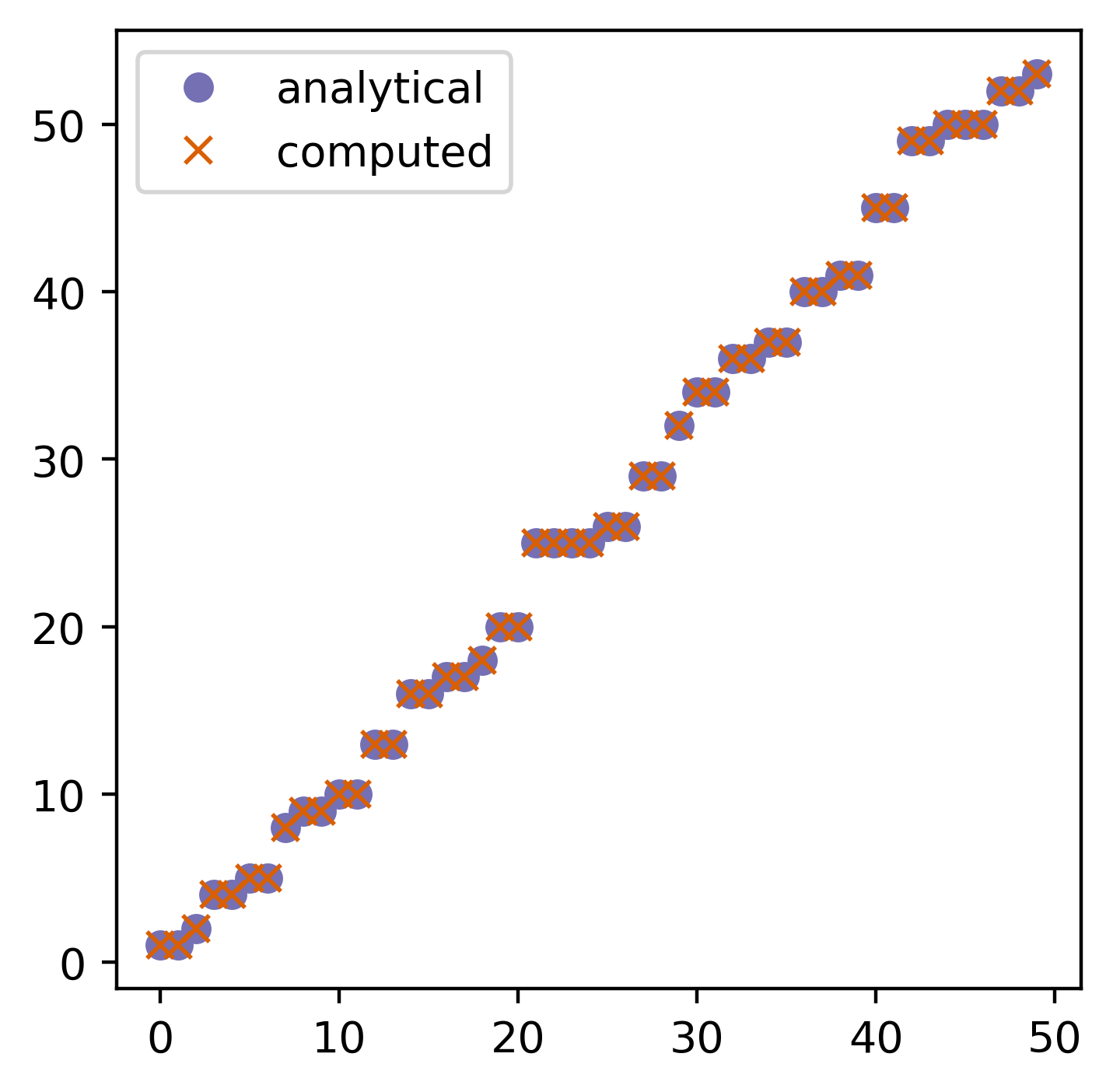}
        \caption{corner-to-corner $\ell=2$}
    \end{subfigure}
    \hfill
    \begin{subfigure}{0.3\textwidth}
        \centering
        \includegraphics[width=\textwidth]{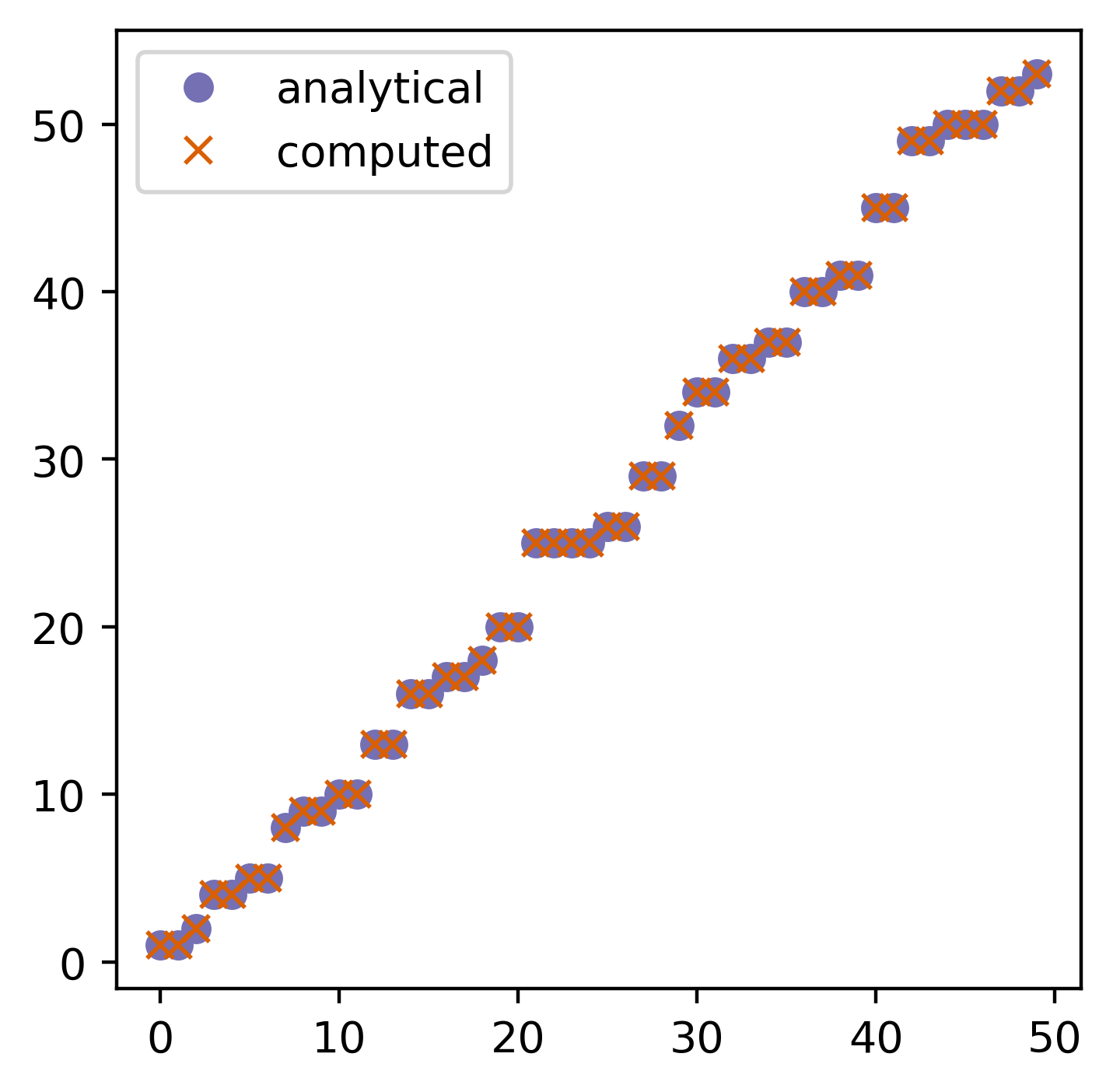}
        \caption{corner-to-corner $\ell=3$}
    \end{subfigure}
    \hfill
    \begin{subfigure}{0.3\textwidth}
        \centering
        \includegraphics[width=\textwidth]{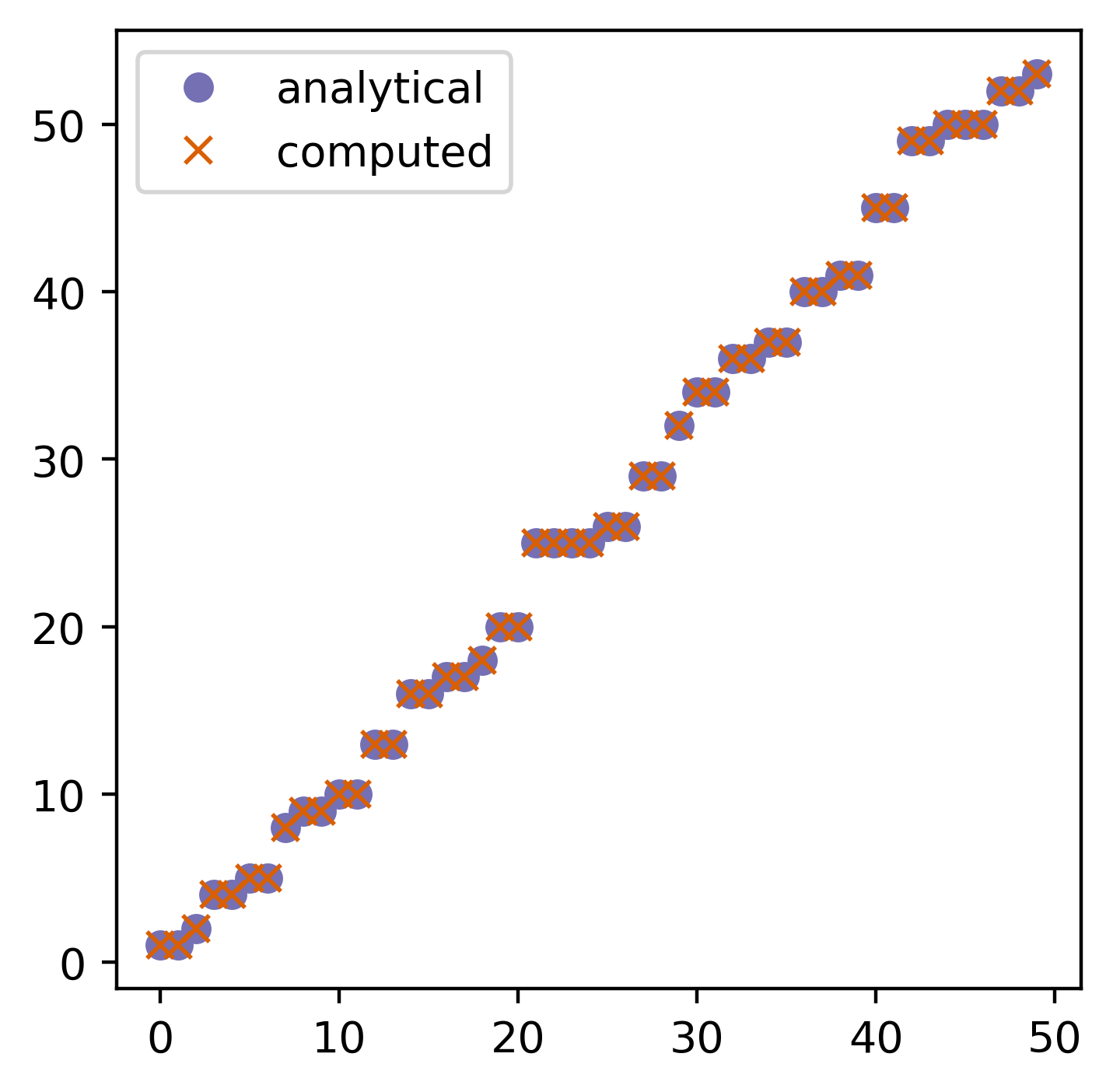}
        \caption{corner-to-corner $\ell=4$}
    \end{subfigure}
    \begin{subfigure}{0.3\textwidth}
        \centering
        \includegraphics[width=\textwidth]{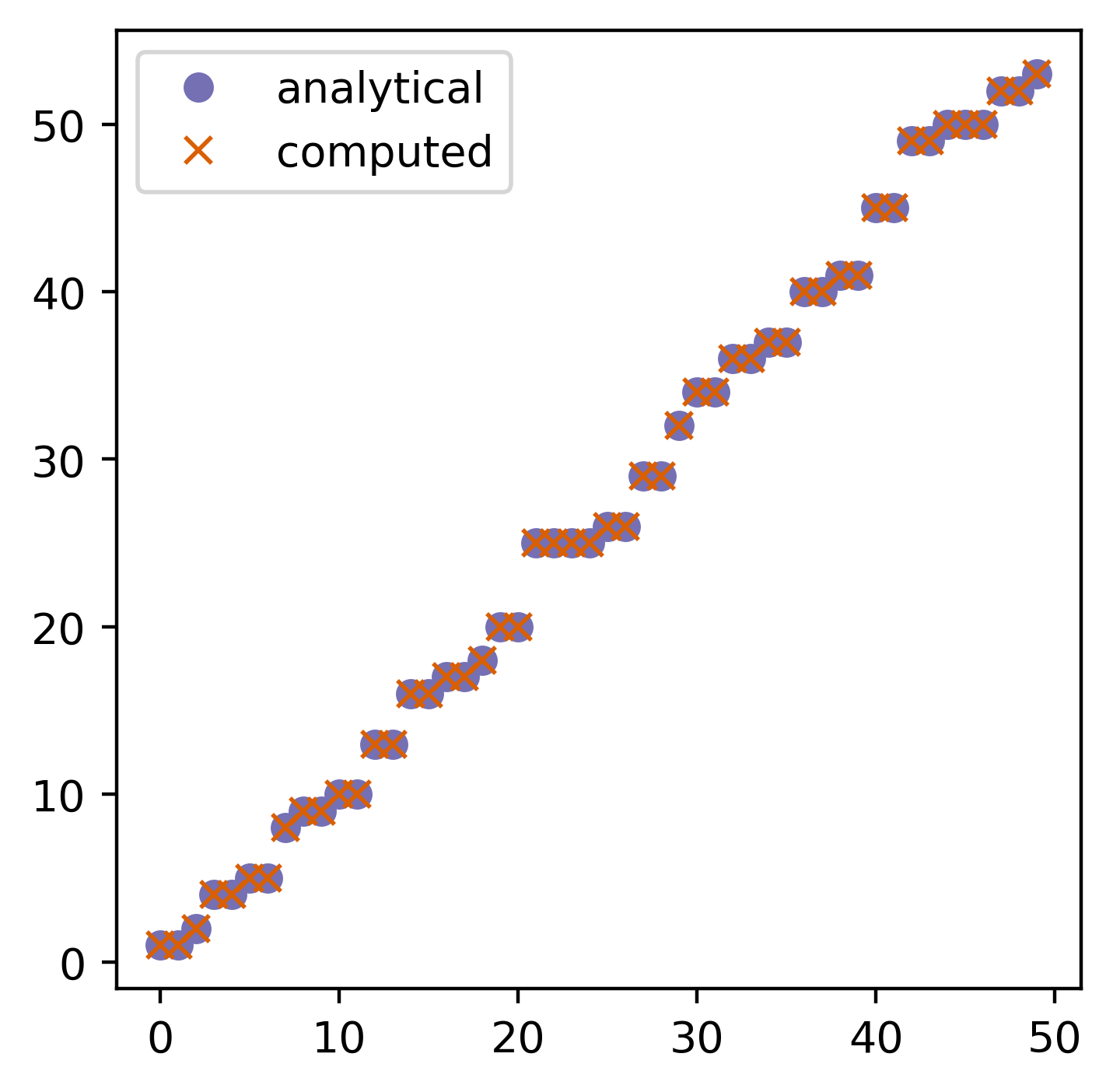}
        \caption{corner-to-corner $\ell=5$}
    \end{subfigure}
    \hfill
    \begin{subfigure}{0.3\textwidth}
        \centering
        \includegraphics[width=\textwidth]{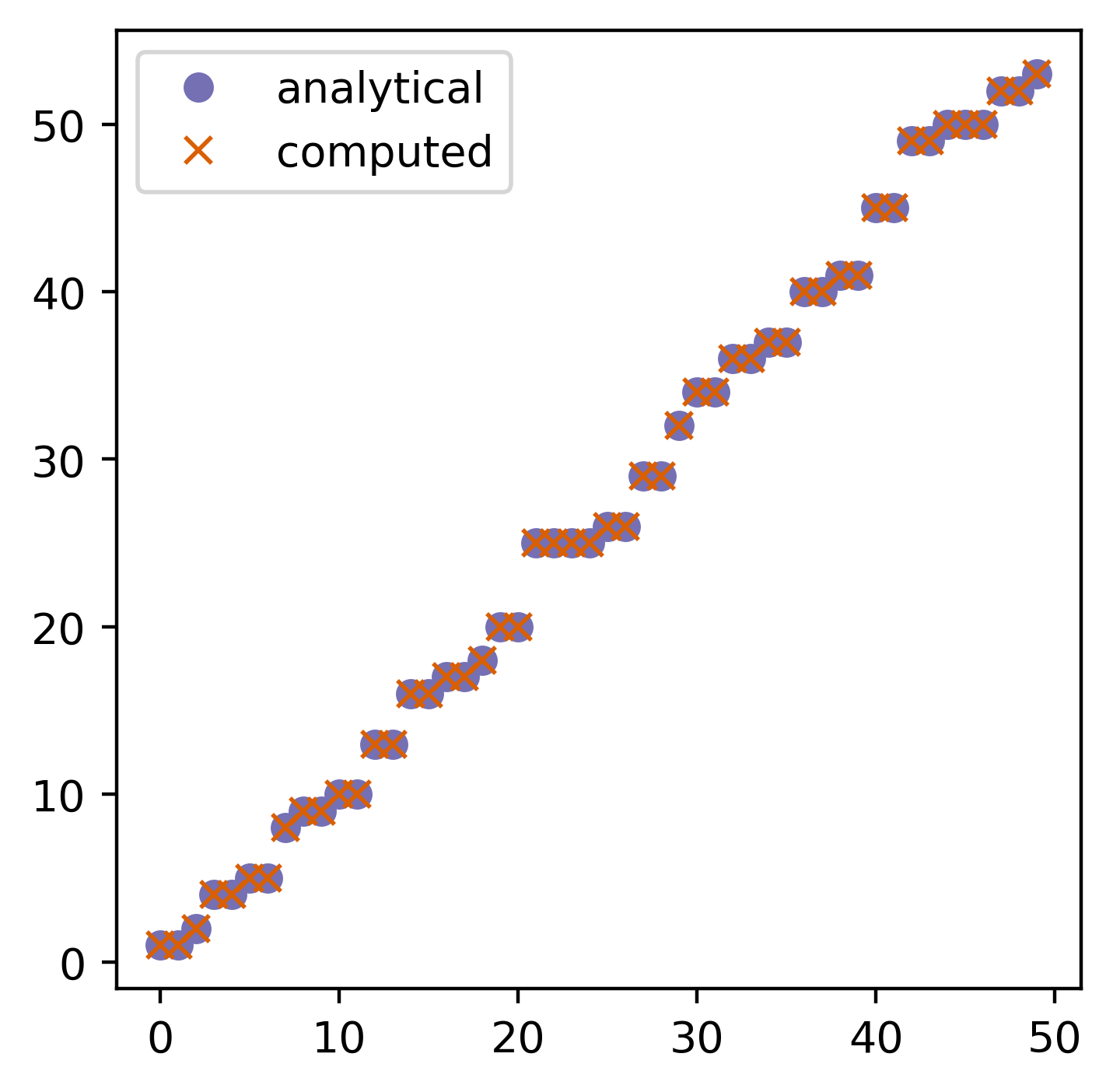}
        \caption{diagonal $\ell=2$}
    \end{subfigure}
    \hfill
    \begin{subfigure}{0.3\textwidth}
        \centering
        \includegraphics[width=\textwidth]{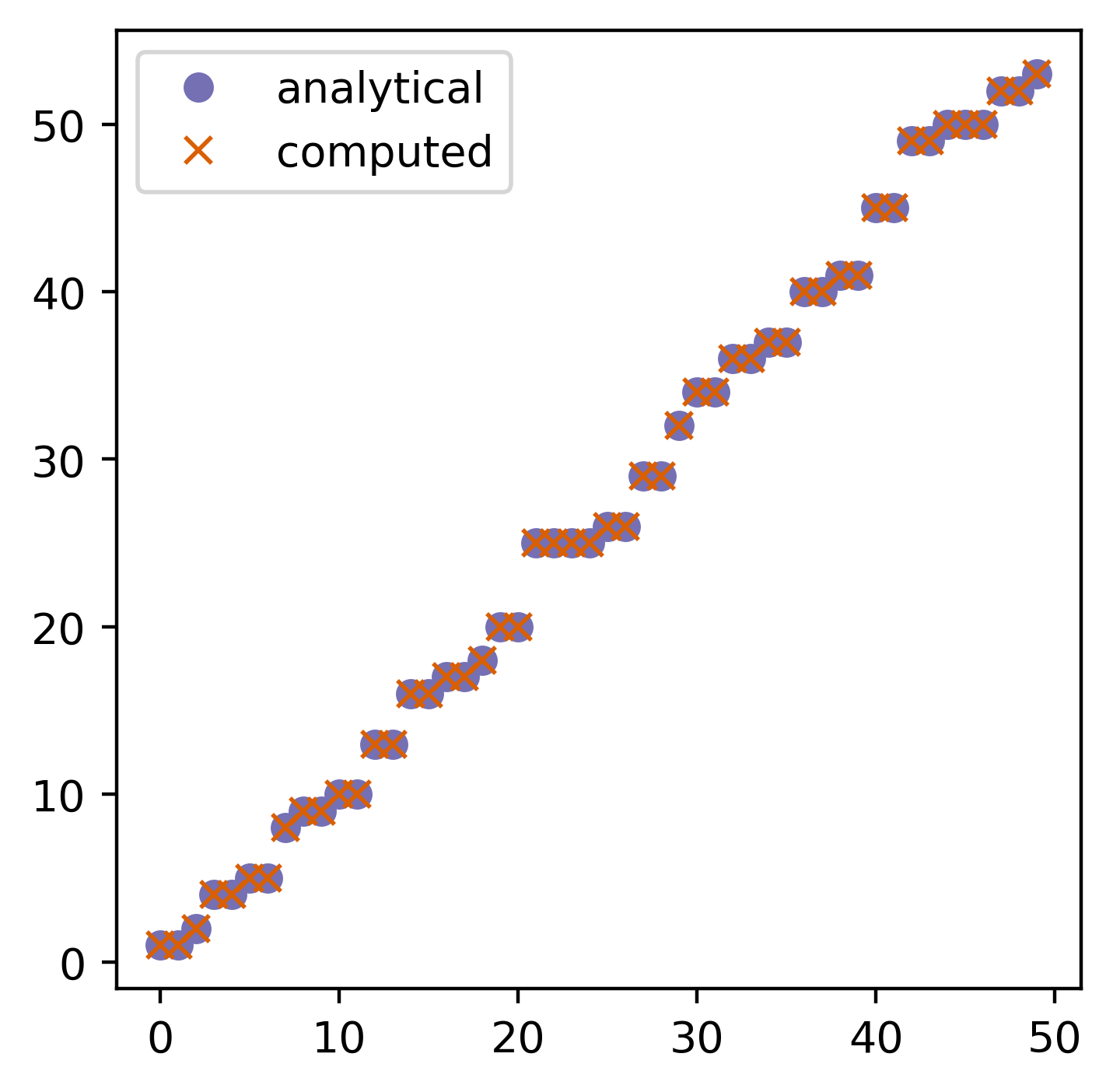}
        \caption{diagonal $\ell=3$}
    \end{subfigure}
    \caption{The first 50 non-zero computed eigenvalues for the meshes of Figure \ref{fig:spuriousmode-check-meshes}, where $\ell$ denotes the number of refinement levels considered. The analytical eigenvalues are also shown, and the spectra all look correct.}
    \label{fig:spuriousmode-check-spectra}
\end{figure}

In all cases tested, we observed no spurious zero eigenvalues.
The first 50 non-zero eigenvalues for the different cases are shown in Figure \ref{fig:spuriousmode-check-spectra}.
We observed no spurious modes in any of the cases.

\subsubsection{Inf-sup condition}
We also investigate the discrete inf-sup condition \cite{chapelleInfsupTest1993}, which plays an essential role in Stokes flow and the Navier-Stokes equations.
For this, we consider the domain $\Omega=[0,2]^2$ and the spaces $\mathbb{V}_{h}^1 := \mathbb{T}_{\DomainHierarchy,[\pp^0,\pp^1-1]}\times\mathbb{T}_{\DomainHierarchy,[p^0-1,p^1]}$ and $\mathbb{V}_h^2 := \Tspace{\DomainHierarchy,[p^0-1,p^1-1]}{}$.
We compute the inf-sup constant $\beta_h$ defined as
\begin{equation}
    \inf_{q_h\in\mathbb{V}_{h}^2} \sup_{\vec{v}_{h} \in\mathbb{V}_{h}^{1}}\frac{\int_{\Omega} q_h \mathrm{div}\vec{v}_h \mathrm{d}V}{\Vert q_h\Vert_{L^2(\Omega)}\Vert \vec{v}_h \Vert_{H(\mathrm{div};\Omega)}} =: \beta_h\;,
\end{equation}
and investigate whether $\beta_h$ goes to zero as the mesh is refined.
For the same meshes considered in Figure \ref{fig:spuriousmode-check-meshes}, Table \ref{tab:spurious-inf-sup-data} shows the computed $\beta_h$.
Here, we note that for $\ell=1$, the spaces are B-spline spaces and the complex is known to be stable (see \cite{buffa2011isogeometric}). 
For all refinement patterns, we observe that the calculated constant is very close to the B-spline constant on the initial tensor-product mesh and does not deteriorate with increasing refinement level.
}

\begin{table}[htb]
    \centering
    \caption{The calculated inf-sup constant for various meshes of Figure \ref{fig:spuriousmode-check-meshes} using THB-splines of degree $\vec{\pp}=(4,4)$.}
    \label{tab:spurious-inf-sup-data}
    \begin{tabular}{r|rr}\hline
\multicolumn{1}{c|}{Refinement level} & \multicolumn{2}{c}{inf-sup constant}                                \\
\multicolumn{1}{c|}{}                 & \multicolumn{1}{c}{corner-to-corner} & \multicolumn{1}{c}{diagonal} \\ \hline \hline
$\ell=1$                             & 0.911867527519551                           & 0.911867527519557                      \\
$\ell=2$                             & 0.911867527519544                             & 0.911867527519528                     \\
$\ell=3$                             & 0.911867527519541                             & 0.911867527519521                      \\
$\ell=4$                             & 0.911867527519544                             & -                            \\
$\ell=5$                             &0.911867527519547                              & - \\ \hline                           
\end{tabular}
\end{table}

\subsection{Adaptive refinement: Structure-preserving discretisations of fluid flow problems}
\begin{figure}[ht]
    \centering
    \begin{subfigure}{0.45\textwidth}
        \includegraphics[width=0.9\textwidth]{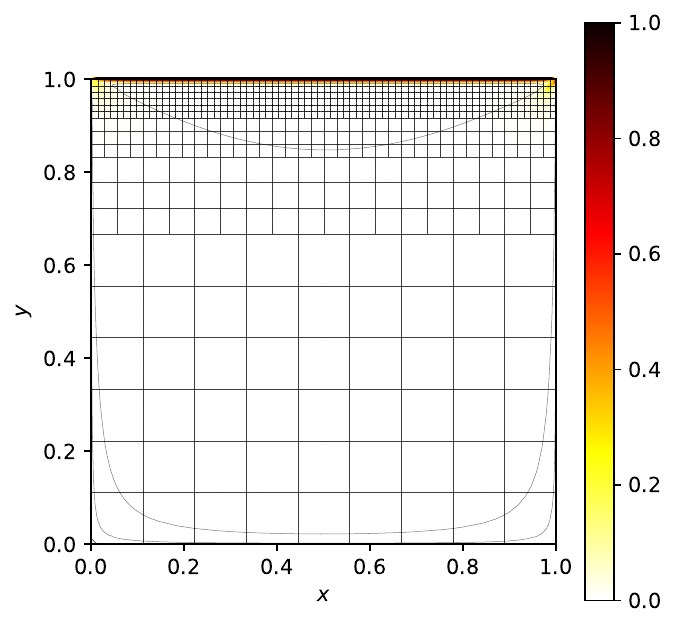}
    \end{subfigure}
    \hfill
    \begin{subfigure}{0.45\textwidth}
        \includegraphics[width=0.9\textwidth]{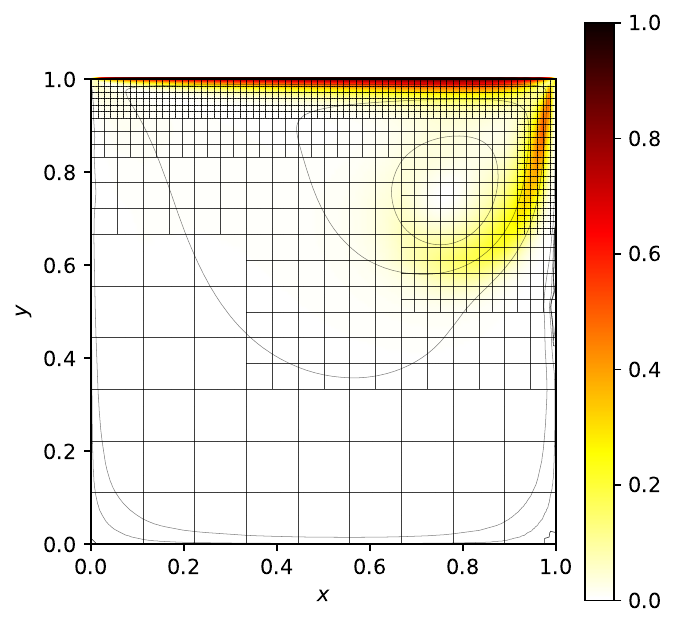}
    \end{subfigure}\\
    \begin{subfigure}{0.45\textwidth}
        \includegraphics[width=0.9\textwidth]{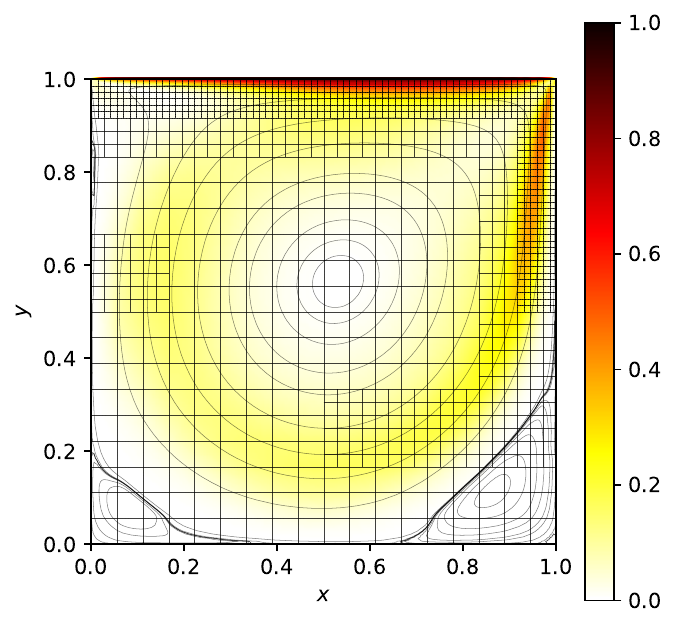}
    \end{subfigure}
    \hfill
    \begin{subfigure}{0.45\textwidth}
        \includegraphics[width=0.9\textwidth]{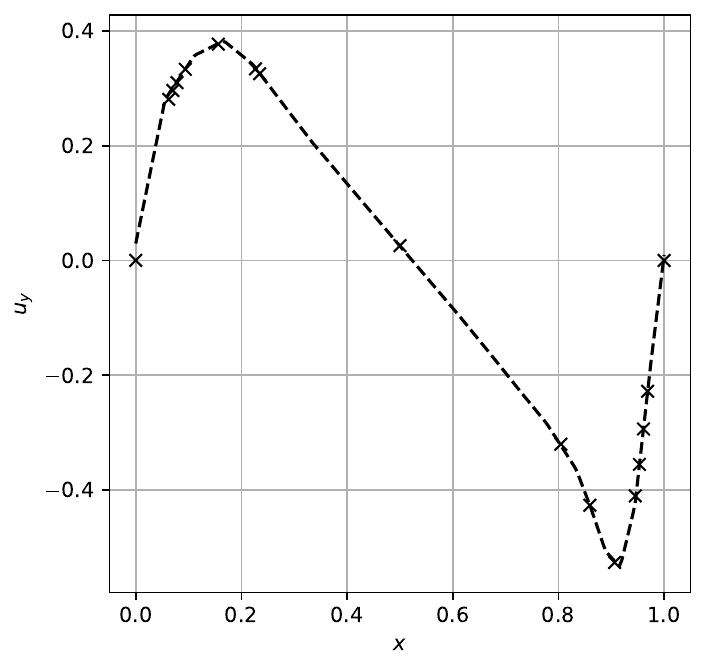}
    \end{subfigure}
    \caption{The numerical results of Method \eqref{eq:discrete-time-dependent-NS} for degree $2$. The mesh has been overlayed on the plots.
    {The colors show the magnitude of the velocity field, and the contourlines show the value of the streamfunction, where the contour height has been taken from \cite[Table 7]{botella_benchmark_1998}). The results} are shown at time $t=0.0$ (top left), at time $t=3.0$ (top right), and at time $t=80.0$ (bottom left).
    The bottom right figure compares the vertical velocity (dashed line) at $t=80.0$ on the horizontal center-line with benchmark results from \cite[Table 10]{botella_benchmark_1998}.}
    \label{fig:results_navier_stokes}
\end{figure}

We validate the use of $(\vec{\pp}+\vec{1})$-boxes for structure-preserving methods by solving the 2D vorticity-velocity-pressure formulation of the time-dependent incompressible Navier-Stokes equations on $\Omega\times [0,T]$, where the mesh is allowed to refine and coarsen.
We choose the domain $\Omega = [0,1]^2$ where the outward facing normal is denoted by $\vec{n}$. Denote $L^2(\Omega)$ as the space of square-integrable functions, from which we define the Hilbert spaces of the de Rham complex \eqref{eq:cont-derham}:
\begin{align}
    H^1(\Omega) &:= \left\{~\sigma\in L^2(\Omega)~:~\mathrm{rot}(\sigma)\in \left[L^2(\Omega)\right]^2~\right\}\;,\\
    H(\mathrm{div};\Omega) &:= \left\{~\vec{v}\in L^2(\Omega)~:~\mathrm{div}(\vec{v})\in L^2(\Omega) ~\right\}\;,
\end{align}
where $\mathrm{rot}(\sigma) := (\partial_y \sigma, -\partial_x \sigma)$ and $\mathrm{div}(\vec{v}) := \partial_x v_x + \partial_y v_y$. In addition, define
\[H_0(\mathrm{div};\Omega) := \left\{~\vec{v}\in H(\mathrm{div};\Omega)~:~\vec{n}\cdot\vec{v} = 0 \text{ on }\partial\Omega ~\right\}\;.\]
The weak form of time-dependent vorticity-velocity-pressure formulation of the two-dimensional lid-driven cavity problem is given by finding $\vec{u}\in H_{0}(\mathrm{div};\Omega),\omega\in H^1(\Omega),p\in L^2(\Omega) $, such that for all $\vec{v}\in H_{0}(\mathrm{div};\Omega),\sigma\in H^1(\Omega),q\in L^2(\Omega) $:
\begin{subequations}
    \begin{align}
        (\partial_t \vec{u},\vec{v})_{\Omega} + (\omega\times\vec{u},\vec{v})_{\Omega} + \mathrm{Re}^{-1}(\mathrm{rot}\left(\omega\right),\vec{v})_{\Omega} -(p,\mathrm{div}(\vec{v}))_{\Omega} &= 0\;,\\
        (\mathrm{div}\left(u\right),q)_{\Omega} &= 0\;,\\
        (\omega,\sigma)_{\Omega} - (\vec{u},\mathrm{rot}\left(\sigma\right))_{\Omega} + (\vec{u}^{\parallel}\times \vec{n},\sigma)_{\partial \Omega} &=0\;.
    \end{align}
\end{subequations}
With initial conditions $\vec{u}_0 = \vec{0}$ and the boundary conditions
\begin{align}
\vec{u}^{\parallel}(\vec{x}) &:= [1,0]^T\;,&\quad & \vec{x}\in\Gamma_{\mathrm{top}}\subset \partial \Omega\;\;,\\
    \vec{u}^{\parallel}(\vec{x}) &:= \vec{0}\;,&\quad & \vec{x}\in \partial \Omega \backslash \Gamma_{\mathrm{top}}\;.
\end{align}
Where $\Gamma_{\mathrm{top}}\subset \partial \Omega$ is the lid of $\Omega$.
We use the discrete spaces \eqref{eq:discrete-spaces-derham-2d}, which we supplement with the boundary conditions:
\begin{equation}
    \mathbb{V}_{h,\vec{0}}^1 := \{~\vec{v}\in \mathbb{V}_h^1~:~ \vec{v}\cdot\vec{n} = {0}~ \}\;.
\end{equation}
Lastly, we temporally discretize the equations using the Crank-Nicolson method, with a time step of $\Delta T$. Then, for time-step $k$, at time $t^k := \Delta T\cdot k$, we discretize the quantities at each full time-step as $[\cdot]_h^k$ and define the half time-step quantities:
\begin{align}
    \left[\cdot\right]_h^{k+\frac{1}{2}} := \frac{1}{2}\left(\left[\cdot\right]_h^{k+1} + \left[\cdot\right]_h^{k} \right)\;.
\end{align}
Define the initial spaces as $\omega_h^0 = 0,\vec{u}_h^0 = \vec{0},p_h^0 = 0$.
We can now introduce our discrete structure-preserving formulation for time step $ k\geq 0$ where the quantities $\omega_h^k\in \mathbb{V}_h^0,\vec{u}_h^k\in\mathbb{V}_{h,\vec{0}}^1,p_h^k\in\mathbb{V}_h^2$ are known. 
Find $\omega_h^{k+1}\in \mathbb{V}_h^0,\vec{u}_h^{k+1}\in \mathbb{V}_{h,\vec{0}}^1,p_h^{k+1}\in\mathbb{V}_h^2$, such that:
\begin{subequations}\label{eq:discrete-time-dependent-NS}
    \begin{align}
        \nonumber \left(\frac{\vec{u}_h^{k+1}-\vec{u}_h^k}{\Delta T},\vec{v}_h\right)_{\Omega} + (\omega_h^{k+\frac{1}{2}}\times \vec{u}_h^{k+\frac{1}{2}},\vec{v}_h)_{\Omega} & & &  \\ 
         + \mathrm{Re}^{-1}(\mathrm{rot}(\omega_h^{k+\frac{1}{2}}),\vec{v}_h)_{\Omega} - (p_h^{k+\frac{1}{2}},\mathrm{div}(\vec{v}_h))_{\Omega} &= 0\;,&  &\forall \vec{v}_h\in\mathbb{V}_{h,\vec{0}}^1\;,\\
        (q_h,\mathrm{div}(\vec{u}_h^{k+\frac{1}{2}})_{\Omega} &=0\;, & &\forall q_h\in\mathbb{V}_h^2\;,\\
        (\omega_h^{k+\frac{1}{2}},\sigma_h)_{\Omega} - (\vec{u}_h^{k+\frac{1}{2}},\mathrm{rot}(\sigma_h))_{\Omega} + (\vec{u}^{\parallel}\times \vec{n},\sigma_h)_{\partial\Omega} &= 0\;, & &\forall \sigma_h \in \mathbb{V}_h^0\;.
    \end{align}
\end{subequations}
For refinement/coarsening, we introduce the following residual-based error indicator at time-step $k$ for a ($\vec{\pp}+\vec{1}$)-box $\pbox{}{}$:
\begin{equation}\label{eq:discrete-time-dependent-ns-estimator}
    \eta^2_k(\pset{}{}) := \sum_{e\in\pset{}{}} \left\Vert \frac{\vec{u}_h^{k+1}-\vec{u}_h^k}{\Delta T} + \omega_h^{k+\frac{1}{2}}\times\vec{u}_h^{k+\frac{1}{2}} + \mathrm{Re}^{-1}\mathrm{curl}(\omega_h^{k+\frac{1}{2}}) + \mathrm{grad}(p_h^{k+\frac{1}{2}}) \right\Vert_{L^2(\Omega^e)}^2\;.  
\end{equation}
Then, based on D\"orfler marking, after each time step, we refine/coarsen the mesh, leading to algorithm \ref{alg:adapt-ref-coar-NS}.
\begin{algorithm}
\begin{algorithmic}[1]
    \State {given $\DomainHierarchy$ and initial conditions $\omega_h^0,\vec{u}_h^0,p_h^0$}
    \For{$k=0$ to $K := T/\Delta T - 1$}
        \Loop
            \State {$\omega_h^{k+1},\vec{u}_h^{k+1},p_h^{k+1} \leftarrow$ solving \eqref{eq:discrete-time-dependent-NS}}
            \State {$\mathcal{E}_{\mathrm{r}} \leftarrow$ estimate error via \eqref{eq:discrete-time-dependent-ns-estimator}}
            \If {$\max{\mathcal{E}_{\mathrm{r}}}< \epsilon$}
                \State {break out of loop}
            \Else
                \State {refine $\DomainHierarchy$ based on $\mathcal{E}_{\mathrm{r}}$ and D\"orfler marking with constant $\Theta_{\mathrm{r}}$}
            \EndIf
        \EndLoop
        \State {$\mathcal{E}_{\mathrm{c}} \leftarrow$ estimate error via \eqref{eq:discrete-time-dependent-ns-estimator}}
        \State {coarsen $\DomainHierarchy$ based on $\mathcal{E}_{\mathrm{c}}$ and d\"orfler marking with constant $\Theta_{\mathrm{c}}$}
    \EndFor
\end{algorithmic}
\caption{The adaptive refinement/coarsening algorithm for the incompressible Navier-Stokes equations.}
\label{alg:adapt-ref-coar-NS}
\end{algorithm}
We solve the incompressible Navier-Stokes equations for $\mathrm{Re}=10^3$, with degree $\vec{p}=[2,2]^T$, and refinement and coarsening parameters $\theta_r = 0.75,\theta_c = 0.02$ and start with a grid of $3\times 3$ ($\vec{\pp}+\vec{1}$)-boxes (or an equivalent $9\times 9$ element mesh), and limit refinement up to 4 levels and the target tolerance is chosen to be $\epsilon=10^{-3}$.
In practice, it was observed that coarsening every time step was excessive, and usually, the subsequent adaptive-refinement loop would reconstruct a similar mesh.
Hence, it was chosen to perform coarsening every $10$ time-steps.
In Figure \ref{fig:results_navier_stokes}, the results are shown for time-steps $t=0.0, t=3.0$ and $t=80.0$.
In addition, we compare our results to the benchmark results from \cite{botella_benchmark_1998}, which agree well.
Note that the small discrepancy at the domain borders is a consequence of weak imposition of tangential boundary conditions in the weak form.

\section{Conclusions}
In this work, we introduced a new adaptive refinement strategy for THB-splines designed to simplify the development of advanced IGA methods where certain critical application-specific mesh properties need to be satisfied.
All theoretical results in the paper are valid for THB-splines of arbitrary smoothness, any admissibility class, and in an arbitrary number of dimensions, $n$.
The core idea is the reinterpretation of the (hierarchical) mesh as a macro-element mesh, where each macro-element is $q_1\times \ldots \times q_n$-elements in size and is referred to as a $\vec{q}$-box.
This reinterpretation leads to a straightforward extension of standard adaptive refinement and coarsening algorithms for THB-splines to our setting: instead of refining individual elements, we refine $\vec{q}$-boxes.
By tailoring the block size $\vec{q}$ to the specific application, this approach can circumvent the need for complex application-specific a posteriori checks and corrections.
We focus on two such applications and our main contributions are threefold:
\begin{itemize}
    \item The development of a \Bezier{}-type local projector using $\vec{p}$-box refinement ($\vec{q}=\vec{p}$). This strategy guarantees local linear independence of the basis on each refined block, thereby resolving four key deficiencies of a previous formulation: it is fully constructive, simple to implement, general with respect to mesh admissibility, and applicable to THB-splines of arbitrary smoothness.

    \item The creation of an adaptive framework for structure-preserving methods using ($\vec{p}+1$)-box refinement ($\vec{q}=\vec{p}+\vec{1}$). This is the first adaptive IGA method that guarantees, a priori, the exactness of the THB-spline de Rham complex on the adaptively refined mesh. This result significantly simplifies the application of adaptive IGA to problems in electromagnetism and incompressible fluid dynamics.

    \item Numerical verification of the framework's performance. We demonstrated optimal convergence rates for adaptive approximation with the $\vec{p}$-box projector and successfully applied the ($\vec{p}+1$)-box strategy to an adaptive, structure-preserving simulation of the incompressible Navier-Stokes equations, confirming the practical viability of our approach.
\end{itemize}

Future research will focus on investigating the disparity between $\vec{\qq}$-boxes and mesh-element refinement for border/corner singularities. 
We want to introduce refinement of a single element at the corner, by considering corner macro-elements of $1\times\dots\times 1$-elements in size.
For this, the remaining border macro-elements have to be redefined as macro-elements of intermediate sizes between $q_1 \times  \dots \times q_n$ and $1\times \dots\times 1$, so that the ``shell'' of macro-elements at the border of the domain is one element thick.
This could be achieved by using knot sequences for the definition of $\vec{\qq}$-boxes, rather than breakpoint sequences; however, in this setup, ensuring nested $\vec{\qq}$-box refinement is less straightforward.

\section{Acknowledgements}

CG 
acknowledges the partial support of the Italian Ministry of University and Research (MUR), through the PRIN projects COSMIC
(No. 2022A79M75) and NOTES (No. P2022NC97R), funded by the European Union - Next Generation EU. CG is a member of
INdAM-GNCS, whose support is gratefully acknowledged.
The research of DT is supported by project number 202.150 awarded through the Veni research programme by the Dutch Research Council (NWO).
The research of KD is supported by the Peter Paul Peterich / DIAM Fast Track Scholarship.


\bibliography{references_LatexFixed}

@article{chapelleInfsupTest1993,
  title = {The Inf-Sup Test},
  author = {Chapelle, D. and Bathe, K.J.},
  date = {1993-06},
  journal = {Computers \& Structures},
  shortjournal = {Computers \& Structures},
  volume = {47},
  number = {4--5},
  pages = {537--545},
  year = {1993},
  issn = {00457949},
  doi = {10.1016/0045-7949(93)90340-J},
  urldate = {2025-12-01},
  langid = {english}
}

@article{borden2011isogeometric,
  title={Isogeometric finite element data structures based on B{\'e}zier extraction of NURBS},
  author={Borden, Michael J and Scott, Michael A and Evans, John A and Hughes, Thomas JR},
  journal={International Journal for Numerical Methods in Engineering},
  volume={87},
  number={1-5},
  pages={15--47},
  year={2011},
  publisher={Wiley Online Library}
}

@article{d2018multi,
  title={Multi-level B{\'e}zier extraction for hierarchical local refinement of Isogeometric Analysis},
  author={D’Angella, Davide and Kollmannsberger, Stefan and Rank, Ernst and Reali, Alessandro},
  journal={Computer Methods in Applied Mechanics and Engineering},
  volume={328},
  pages={147--174},
  year={2018},
  publisher={Elsevier}
}

@article{johannessen2015divergence,
  title={Divergence-conforming discretization for Stokes problem on locally refined meshes using LR B-splines},
  author={Johannessen, Kjetil Andr{\'e} and Kumar, Mukesh and Kvamsdal, Trond},
  journal={Computer Methods in Applied Mechanics and Engineering},
  volume={293},
  pages={38--70},
  year={2015},
  publisher={Elsevier}
}

@article{evans2020hierarchical,
  title={Hierarchical B-spline complexes of discrete differential forms},
  author={Evans, John A and Scott, Michael A and Shepherd, Kendrick M and Thomas, Derek C and V{\'a}zquez Hern{\'a}ndez, Rafael},
  journal={IMA Journal of Numerical Analysis},
  volume={40},
  number={1},
  pages={422--473},
  year={2020},
  publisher={Oxford University Press}
}

@article{buffa2014isogeometric,
  title={Isogeometric methods for computational electromagnetics: B-spline and T-spline discretizations},
  author={Buffa, Annalisa and Sangalli, Giancarlo and V{\'a}zquez, Rafael},
  journal={Journal of Computational Physics},
  volume={257},
  pages={1291--1320},
  year={2014},
  publisher={Elsevier}
}

@article{buffa2011isogeometric,
  title={Isogeometric discrete differential forms in three dimensions},
  author={Buffa, Annalisa and Rivas, Judith and Sangalli, Giancarlo and V{\'a}zquez, Rafael},
  journal={SIAM Journal on Numerical Analysis},
  volume={49},
  number={2},
  pages={818--844},
  year={2011},
  publisher={SIAM}
}

@article{evans2013isogeometric,
  title={Isogeometric divergence-conforming B-splines for the unsteady Navier--Stokes equations},
  author={Evans, John A and Hughes, Thomas JR},
  journal={Journal of Computational Physics},
  volume={241},
  pages={141--167},
  year={2013},
  publisher={Elsevier}
}

@article{botella_benchmark_1998,
	title = {Benchmark spectral results on the lid-driven cavity flow},
	volume = {27},
	copyright = {https://www.elsevier.com/tdm/userlicense/1.0/},
	issn = {00457930},
    doi = {10.1016/S0045-7930(98)00002-4},
	language = {en},
	number = {4},
	urldate = {2025-09-08},
	journal = {Computers \& Fluids},
	author = {Botella, O. and Peyret, R.},
	month = may,
	year = {1998},
	pages = {421--433},
}

@book{arnold_finite_2018,
	address = {Philadelphia SE - xi, 120 pages : illustrations ; 26 cm.},
	title = {Finite {Element} {Exterior} {Calculus}},
	isbn = {978-1-61197-554-3},
	language = {English},
	publisher = {Society for Industrial and Applied Mathematics},
	author = {Arnold, D. N.},
	year = {2018},
	doi = {https://doi.org/10.1137/1.9781611975543},
	note = {Publication Title: CBMS-NSF regional conference series in applied mathematics ; 93 TA - TT -},
}

@misc{cabanas_construction_2025,
	title = {Construction of exact refinements for the two-dimensional {HB}/{THB}-spline de {Rham} complex},
	copyright = {Creative Commons Attribution Non Commercial Share Alike 4.0 International},
	doi = {10.48550/ARXIV.2502.19542},
	urldate = {2025-05-01},
	publisher = {arXiv},
	author = {Cabanas, D. C. and Shepherd, K. M. and Toshniwal, D. and Vázquez, R.},
	year = {2025},
	note = {Version Number: 2},
}

@article{shepherd_locally-verifiable_2024,
	title = {Locally-{Verifiable} {Sufficient} {Conditions} for {Exactness} of the {Hierarchical} {B}-spline {Discrete} de {Rham} {Complex} in $\mathbb{R}^n$},
	issn = {1615-3375, 1615-3383},
	doi = {10.1007/s10208-024-09659-6},
	language = {en},
	journal = {FoCM},
	author = {Shepherd, K. and Toshniwal, D.},
	month = dec,
	year = {2024},
}

@article{speleers_effortless_2016,
	title = {Effortless quasi-interpolation in hierarchical spaces},
	volume = {132},
	issn = {0029-599X, 0945-3245},
	doi = {10.1007/s00211-015-0711-z},
	language = {en},
	number = {1},
	urldate = {2025-01-22},
	journal = {Numer. Math.},
	author = {Speleers, H. and Manni, C.},
	year = {2016},
	pages = {155--184},
}

@misc{dijkstra_nutils-pbox_2024,
	title = {nutils-pbox},
	url = {https://gitlab.tudelft.nl/kwdijkstra/nutils-pbox},
	author = {Dijkstra, K.W.},
	year = {2024},
}

@incollection{lanini_characterization_2024,
	address = {Singapore},
	title = {A {Characterization} of {Linear} {Independence} of {THB}-{Splines} in $\mathbb{R}^n$ and {Application} to {Bézier} {Projection}},
	volume = {60},
	isbn = {978-981-9765-07-2 978-981-9765-08-9},
	language = {en},
	urldate = {2025-01-16},
	booktitle = {Approximation {Theory} and {Numerical} {Analysis} {Meet} {Algebra}, {Geometry}, {Topology}},
	publisher = {Springer Nature Singapore},
	author = {Dijkstra, K.W. and Toshniwal, D.},
	editor = {Lanini, M. and Manni, C. and Schenck, H.},
	year = {2024},
	pages = {115--148},
}

@misc{van_zwieten_nutils_2022,
	title = {Nutils},
	copyright = {MIT License},
	url = {https://zenodo.org/doi/10.5281/zenodo.6006701},
	urldate = {2024-07-03},
	publisher = {Zenodo},
	author = {van Zwieten, J.S.B and van Zwieten, G.J. and Hoitinga, W.},
	year = {2022},
}

@article{li_s-splines_2019,
	title = {S-splines: {A} simple surface solution for {IGA} and {CAD}},
	volume = {350},
	issn = {00457825},
	doi = {10.1016/j.cma.2019.03.035},
	journal = {Comput. Methods Appl. Mech. Eng.},
	author = {Li, X. and Sederberg, T.W.},
	year = {2019},
	pages = {664--678},
}

@article{thomas_bezier_2015,
	title = {Bézier projection: {A} unified approach for local projection and quadrature-free refinement and coarsening of {NURBS} and {T}-splines with particular application to isogeometric design and analysis},
	volume = {284},
	issn = {00457825},
	doi = {10.1016/j.cma.2014.07.014},
	journal = {Comput. Methods Appl. Mech. Eng.},
	author = {Thomas, D. C. and Scott, M. A. and Evans, J. A. and Tew, K. and Evans, E. J.},
	year = {2015},
	pages = {55--105},
}

@article{dokken_polynomial_2013,
	title = {Polynomial splines over locally refined box-partitions},
	volume = {30},
	issn = {01678396},
	doi = {10.1016/j.cagd.2012.12.005},
	number = {3},
	journal = {CAGD},
	author = {Dokken, T. and Lyche, T. and Pettersen, K.F.},
	year = {2013},
	pages = {331--356},
}

@article{giannelli_thb-splines_2012,
	title = {{THB}-splines: {The} truncated basis for hierarchical splines},
	volume = {29},
	issn = {01678396},
	shorttitle = {{THB}-splines},
	doi = {10.1016/j.cagd.2012.03.025},
	language = {en},
	number = {7},
	journal = {CAGD},
	author = {Giannelli, C. and Jüttler, B. and Speleers, H.},
	year = {2012},
	pages = {485--498},
}

@incollection{devore_theory_2009,
	address = {Berlin, Heidelberg},
	title = {Theory of adaptive finite element methods: {An} introduction},
	isbn = {978-3-642-03412-1 978-3-642-03413-8},
	shorttitle = {Theory of adaptive finite element methods},
	language = {en},
	urldate = {2024-05-23},
	booktitle = {Multiscale, {Nonlinear} and {Adaptive} {Approximation}},
	publisher = {Springer Berlin Heidelberg},
	author = {Nochetto, R.H. and Siebert, K.G. and Veeser, A.},
	editor = {DeVore, R. and Kunoth, A.},
	year = {2009},
	pages = {409--542},
}

@article{sederberg_t-splines_2003,
	title = {T-{Splines} and {T}-{NURCCs}},
	volume = {22},
	issn = {0730-0301},
	doi = {10.1145/882262.882295},
	number = {3},
	journal = {ACM Trans. Graph.},
	author = {Sederberg, T.W. and Zheng, J. and Bakenov, A. and Nasri, A.},
	year = {2003},
	pages = {477--484},
}

@article{forsey_hierarchical_1988,
	title = {Hierarchical {B}-{Spline} {Refinement}},
	volume = {22},
	doi = {10.1145/378456.378512},
	number = {4},
	journal = {SIGGRAPH Comput. Graph.},
	author = {Forsey, D.R. and Barrels, R.H.},
	year = {1988},
	pages = {205--212},
}

@article{sande_explicit_2020,
	title = {Explicit error estimates for spline approximation of arbitrary smoothness in isogeometric analysis},
	volume = {144},
	issn = {0029-599X, 0945-3245},
	doi = {10.1007/s00211-019-01097-9},
	language = {en},
	number = {4},
	urldate = {2025-01-16},
	journal = {Numer. Math.},
	author = {Sande, E. and Manni, C. and Speleers, H.},
	year = {2020},
	pages = {889--929},
}

@article{patrizi_linear_2020,
	title = {Linear dependence of bivariate {Minimal} {Support} and {Locally} {Refined} {B}-splines over {LR}-meshes},
	volume = {77},
	issn = {01678396},
	doi = {10.1016/j.cagd.2019.101803},
	language = {en},
	urldate = {2023-07-26},
	journal = {CAGD},
	author = {Patrizi, F. and Dokken, T.},
	year = {2020},
	pages = {101803},
}

@article{giust_local_2020,
	title = {Local ({T}){HB}-spline projectors via restricted hierarchical spline fitting},
	volume = {80},
	issn = {0167-8396},
	doi = {https://doi.org/10.1016/j.cagd.2020.101865},
	journal = {CAGD},
	author = {Giust, A. and Jüttler, B. and Mantzaflaris, A.},
	year = {2020},
	pages = {101865},
}

@article{carraturo_suitably_2019,
	title = {Suitably graded {THB}-spline refinement and coarsening: {Towards} an adaptive isogeometric analysis of additive manufacturing processes},
	volume = {348},
	issn = {00457825},
	shorttitle = {Suitably graded {THB}-spline refinement and coarsening},
	doi = {10.1016/j.cma.2019.01.044},
	language = {en},
	urldate = {2024-11-26},
	journal = {Comput. Methods Appl. Mech. Eng.},
	author = {Carraturo, M. and Giannelli, C. and Reali, A. and Vázquez, R.},
	year = {2019},
	pages = {660--679},
}

@article{bracco_refinement_2018,
	title = {Refinement {Algorithms} for {Adaptive} {Isogeometric} {Methods} with {Hierarchical} {Splines}},
	volume = {7},
	copyright = {https://creativecommons.org/licenses/by/4.0/},
	issn = {2075-1680},
	doi = {10.3390/axioms7030043},
	language = {en},
	number = {3},
	urldate = {2024-12-02},
	journal = {Axioms},
	author = {Bracco, C. and Giannelli, C. and Vázquez, R.},
	year = {2018},
	pages = {43},
}

@article{gantner_adaptive_2017,
	title = {Adaptive {IGAFEM} with optimal convergence rates: {Hierarchical} {B}-splines},
	volume = {27},
	issn = {0218-2025, 1793-6314},
	shorttitle = {Adaptive {IGAFEM} with optimal convergence rates},
	doi = {10.1142/S0218202517500543},
	language = {en},
	number = {14},
	urldate = {2025-01-09},
	journal = {Math. Models Methods Appl. Sci.},
	author = {Gantner, G. and Haberlik, D. and Praetorius, D.},
	year = {2017},
	pages = {2631--2674},
}

@article{buffa_adaptive_2017,
	title = {Adaptive isogeometric methods with hierarchical splines: {Optimality} and convergence rates},
	volume = {27},
	issn = {0218-2025, 1793-6314},
	shorttitle = {Adaptive isogeometric methods with hierarchical splines},
	doi = {10.1142/S0218202517500580},
	language = {en},
	number = {14},
	urldate = {2025-01-08},
	journal = {Math. Models Methods Appl. Sci.},
	author = {Buffa, A. and Giannelli, C.},
	year = {2017},
	pages = {2781--2802},
}

@article{buffa_complexity_2016,
	title = {Complexity of hierarchical refinement for a class of admissible mesh configurations},
	volume = {47},
	issn = {01678396},
	doi = {10.1016/j.cagd.2016.04.003},
	language = {en},
	urldate = {2025-01-08},
	journal = {CAGD},
	author = {Buffa, A. and Giannelli, C. and Morgenstern, P. and Peterseim, D.},
	year = {2016},
	pages = {83--92},
}

@article{buffa_adaptive_2016,
	title = {Adaptive isogeometric methods with hierarchical splines: {Error} estimator and convergence},
	volume = {26},
	issn = {0218-2025, 1793-6314},
	shorttitle = {Adaptive isogeometric methods with hierarchical splines},
	doi = {10.1142/S0218202516500019},
	language = {en},
	number = {01},
	urldate = {2025-01-16},
	journal = {Math. Models Methods Appl. Sci.},
	author = {Buffa, A. and Giannelli, C.},
	year = {2016},
	pages = {1--25},
}

@article{evans_hierarchical_2015,
	title = {Hierarchical {T}-splines: {Analysis}-suitability, {Bézier} extraction, and application as an adaptive basis for isogeometric analysis},
	volume = {284},
	issn = {00457825},
	shorttitle = {Hierarchical {T}-splines},
	doi = {10.1016/j.cma.2014.05.019},
	language = {en},
	urldate = {2025-01-09},
	journal = {Comput. Methods Appl. Mech. Eng.},
	author = {Evans, E.J. and Scott, M.A. and Li, X. and Thomas, D.C.},
	year = {2015},
	pages = {1--20},
}

@article{kuru_goal-adaptive_2014,
	title = {Goal-adaptive {Isogeometric} {Analysis} with hierarchical splines},
	volume = {270},
	issn = {00457825},
	doi = {10.1016/j.cma.2013.11.026},
	language = {en},
	urldate = {2025-01-09},
	journal = {Comput. Methods Appl. Mech. Eng.},
	author = {Kuru, G. and Verhoosel, C.V. and Van Der Zee, K.G. and Van Brummelen, E.H.},
	year = {2014},
	pages = {270--292},
}

@article{johannessen_isogeometric_2014,
	title = {Isogeometric analysis using {LR} {B}-splines},
	volume = {269},
	issn = {00457825},
	doi = {10.1016/j.cma.2013.09.014},
	language = {en},
	urldate = {2023-03-01},
	journal = {Comput. Methods Appl. Mech. Eng.},
	author = {Johannessen, K.A. and Kvamsdal, T. and Dokken, T.},
	year = {2014},
	pages = {471--514},
}

@article{giannelli_bases_2013,
	title = {Bases and dimensions of bivariate hierarchical tensor-product splines},
	volume = {239},
	copyright = {https://www.elsevier.com/tdm/userlicense/1.0/},
	issn = {03770427},
	doi = {10.1016/j.cam.2012.09.031},
	language = {en},
	urldate = {2024-10-22},
	journal = {J. Comput. Appl. Math.},
	author = {Giannelli, C. and Jüttler, B.},
	year = {2013},
	pages = {162--178},
}

@article{bressan_properties_2013,
	title = {Some properties of {LR}-splines},
	volume = {30},
	issn = {01678396},
	doi = {10.1016/j.cagd.2013.06.004},
	language = {en},
	number = {8},
	urldate = {2023-07-26},
	journal = {CAGD},
	author = {Bressan, A.},
	year = {2013},
	pages = {778--794},
}

@article{dorfel_adaptive_2010,
	title = {Adaptive isogeometric analysis by local h-refinement with {T}-splines},
	volume = {199},
	copyright = {https://www.elsevier.com/tdm/userlicense/1.0/},
	issn = {00457825},
	doi = {10.1016/j.cma.2008.07.012},
	language = {en},
	number = {5-8},
	urldate = {2025-01-09},
	journal = {Comput. Methods Appl. Mech. Eng.},
	author = {Dörfel, M.R. and Jüttler, B. and Simeon, B.},
	year = {2010},
	pages = {264--275},
}

@article{buffa_linear_2010,
	title = {Linear independence of the {T}-spline blending functions associated with some particular {T}-meshes},
	volume = {199},
	issn = {00457825},
	doi = {10.1016/j.cma.2009.12.004},
	language = {en},
	number = {23-24},
	urldate = {2023-07-26},
	journal = {Comput. Methods Appl. Mech. Eng.},
	author = {Buffa, A. and Cho, D. and Sangalli, G.},
	year = {2010},
	pages = {1437--1445},
}

@article{bazilevs_isogeometric_2010,
	title = {Isogeometric analysis using {T}-splines},
	volume = {199},
	copyright = {https://www.elsevier.com/tdm/userlicense/1.0/},
	issn = {00457825},
	doi = {10.1016/j.cma.2009.02.036},
	language = {en},
	number = {5-8},
	urldate = {2025-01-09},
	journal = {Comput. Methods Appl. Mech. Eng.},
	author = {Bazilevs, Y. and Calo, V.M. and Cottrell, J.A. and Evans, J.A. and Hughes, T.J.R. and Lipton, S. and Scott, M.A. and Sederberg, T.W.},
	year = {2010},
	pages = {229--263},
}

@article{hughes_isogeometric_2005,
	title = {Isogeometric analysis: {CAD}, finite elements, {NURBS}, exact geometry and mesh refinement},
	volume = {194},
	issn = {00457825},
	shorttitle = {Isogeometric analysis},
	doi = {10.1016/j.cma.2004.10.008},
	language = {en},
	number = {39-41},
	urldate = {2023-03-01},
	journal = {Comput. Methods Appl. Mech. Eng.},
	author = {Hughes, T.J.R. and Cottrell, J.A. and Bazilevs, Y.},
	year = {2005},
	pages = {4135--4195},
}

@article{beirao_da_veiga_estimates_2011,
	title = {Some estimates for h-p-k-refinement in {Isogeometric} {Analysis}},
	volume = {118},
	issn = {0029599X},
	doi = {10.1007/s00211-010-0338-z},
	number = {2},
	journal = {Numer. Math.},
	author = {Beirão da Veiga, L. and Buffa, A. and Rivas, J. and Sangalli, G.},
	year = {2011},
	pages = {271--305},
}

@article{vuong_hierarchical_2011,
	title = {A hierarchical approach to adaptive local refinement in isogeometric analysis},
	volume = {200},
	issn = {00457825},
	doi = {10.1016/j.cma.2011.09.004},
	number = {49-52},
	journal = {Comput. Methods Appl. Mech. Eng.},
	author = {Vuong, A. V. and Giannelli, C. and Jüttler, B. and Simeon, B.},
	year = {2011},
	pages = {3554--3567},
}

@phdthesis{kraft_adaptive_1998,
	address = {Stuttgard},
	title = {Adaptive und linear unabh\"angige {Multilevel} {B}-{Splines} und ihre {Anwendungen}},
	school = {Univ. Stuttgard},
	author = {Kraft, R.},
	year = {1998},
	note = {Publication Title: PhD Thesis},
}
\bibliographystyle{elsarticle-num}

\appendix
\section{Generic \textbf{\textit{q}}-box results}\label{sec:AppendixA}

The proof of Propositions \ref{prop:regular-p-box-spline-level} is straightforward.
\begin{proof}[Proof: Proposition \ref{prop:regular-p-box-spline-level}]
    By the definition of the THB-splines, only level $k\leq \ell$ THB-splines have to be considered.
    Clearly, a level $k < \ell$ THB-spline $\Tspline{j}{}$ can be written as a linear combination of $\Bbasis{\ell}{\nact,\ell}$ B-splines, due to it being truncated to level $k$,
    \begin{eqnarray}
        \Tspline{j}{} \in \lspan{\Bbasis{\ell}{\nact,\ell}}\;.
    \end{eqnarray}
    As the supports of $\Bspline{}{}\in\Bbasis{\ell}{\nact,\ell}$ are $\pp^i+1$ elements wide, of which at least one element cannot be contained in $\Omega_\ell$, these supports can only cover level-$\ell$ border $\vec{\qq}$-boxes.
    But then, no spline of level $k<\ell$ can have support on level-$\ell$ regular $\vec{\qq}$-boxes.
\end{proof}

The proof of Proposition \ref{prop:sufficient-condition-admiss-class-2} is the direct result of the following two Lemmas.
\begin{lemma}\label{lem:iff1}
    Given any $\vec{\qq}$-box mesh consisting of active border $\vec{\qq}$-boxes and active regular $\vec{\qq}$-boxes where $\qq^i\geq \pp^i$. Then all the active $\vec{\qq}$-boxes are of admissibility class $c=2$.
\end{lemma}
\begin{proof}
    Consider an active level-$\ell$ $\vec{\qq}$-box $\pbox{\vec{r},\ell}{}$ for which we want to show that no splines of level $\ell-2$ or lower have support.
    It is sufficient to only consider level-$(\ell-2)$ B-splines, as due to the truncation operator, any lower level spline becomes a linear combination of level-$(\ell-2)$ B-splines.
    Truncating a level-$(\ell-2)$ B-spline to level $\ell-1$, a linear combination of $\Bbasis{\ell-1}{\nact,\ell-1}$ is retrieved.
    \begin{eqnarray}
        \trunc{\ell-1}{\Bbasis{\ell-2}{}} \subset \lspan{\Bbasis{\ell-1}{\nact,\ell-1}}\;.
    \end{eqnarray}
    Crucially, the splines of $\Bbasis{\ell-1}{\nact,\ell-1}$ have supports of at most $\pp^i+1$ elements in direction $i$.
    Secondly, these splines must have support outside of $\Omega_{\ell-1}$, and as such, can only penetrate $\Omega_{\ell-1}$ at most $\pp^i$ elements in each direction.
    Hence, within $\Omega_{\ell-1}$, these splines can only be supported by level-$(\ell-1)$ border $\vec{\qq}$-boxes.
    But as these level-$(\ell-1)$ border $\vec{\qq}$-boxes are active, these truncated level-$(\ell-2)$ B-splines have no support on any active level-$\ell$ $\vec{\qq}$-box, as desired.
\end{proof}

\begin{lemma}\label{lem:iff2}
    Any well-behaved $\vec{\qq}$-box of admissibility class $2$ with $\qq^i\geq\pp^i$ is an active border $\vec{\qq}$-box.
\end{lemma}
\begin{proof}
    Without loss of generality, assume that this well-behaved $\vec{\qq}$-box $\pbox{\vec{r},\ell}{}$ is of level $\ell>1$. Then, by definition, it neighbours some level-$m<\ell$ (potentially not active) element $\Omega^e$. 
    Hence, a B-spline $\Bspline{\vec{j},m}{}\in \Bbasis{m}{}$ exsists that has support on both $\pbox{\vec{r},\ell}{}$ and $\Omega^e$. 
    
    In addition, for any active $\vec{\qq}$-box $\pbox{\vec{s},k}{}$ of level $k$, by Assumption \ref{ass:p-box-refinement-domain}, $\Omega_k$ contains at least $\prod_{i=1}^d 2 \pp$ level-$k$ mesh elements, of which $\prod_{i=1}^d \pp$ elements are contained in $\pbox{\vec{s},k}{}$. As a result, a level-$k$ B-spline $\Bspline{\vec{i},k}{}\in\Bbasis{k}{}$ exists with support on $\pbox{\vec{s},k}{}$, and $\supp{\Bspline{\vec{i},k}{}}\subset \Omega_k$. Hence, any active $\vec{\qq}$-box of level $k$ supports THB-splines of level $k$.
    
    But then, since $\pbox{\vec{r},\ell}{}$ only supports splines of two levels, these levels must be $\ell-1$ and $\ell$, showing that $\pbox{\vec{r},\ell}{}$ is an active border $\vec{\qq}$-box, as it cannot contain any higher level-$k>\ell$ active border $\vec{\qq}$-boxes.
\end{proof}

\begin{proof}[Proof: Proposition \ref{prop:sufficient-condition-admiss-class-2}]
    This is a direct consequence of Proposition \ref{prop:regular-p-box-spline-level}, Lemma \ref{lem:iff1} and Lemma \ref{lem:iff2}.
\end{proof}

\section{Proof of Lemma \ref{lem:non-overloaded-border-p-box}}\label{sec:appendixB}
For Lemma \ref{lem:non-overloaded-border-p-box}, we only need to consider active border $\vec{\pp}$-boxes and regular $\vec{\pp}$-boxes.
By Proposition \ref{prop:regular-p-box-spline-level}, regular $\vec{\pp}$-boxes are never overloaded w.r.t. $\Tbasis{\DomainHierarchy}{}$, leaving the subject of active border $\vec{\pp}$-boxes for this section.

Before proceeding with the details of the proof, we would like to highlight the intuition behind the steps.
For the active border $\vec{\pp}$-box $\pbox{\vec{r},\ell}{}$, we want to show that the set 
$$\mathcal{D}:= \{~\Tspline{}{}\in \Tbasis{\DomainHierarchy}{}~:~\pbox{\vec{r},\ell}{}\cap \supp{\Tspline{}{}}\neq \emptyset~\}$$
is linearly independent on $\pbox{}{}$. For this, we consider an ordering of the mesh-elements $(e_1,e_2,e_3,\dots)$ for $e_i\in \pset{}{}$, and define the sets
$$ \mathcal{D}_i:= \{~\Tspline{}{}\in\mathcal{D}~:~\Omega^{e_i}\in \supp{\Tspline{}{}},~\Omega^{e_j}\notin \supp{\Tspline{}{}},~\forall j=1,\dots,i-1~\}\;. $$
Then, $\mathcal{D} = \cup_{e_i\in\pset{}{}} \mathcal{D}_i$ and if each $\mathcal{D}_i$ is linearly independent on $\Omega^{e_i}$, so must $\mathcal{D}$ be on $\pbox{}{}$.
The primary benefit of this approach is that on each element, all splines are polynomials, which simplifies the proof of linear independence. 
However, finding the correct ordering is not trivial, and depends on the active border $\vec{\pp}$-box $\pbox{\vec{r},\ell}{}$ and its relation to the boundary of the refinement domain $\Omega_\ell$.

In Section \ref{sec:refinement-border-characterization}, this relation between border $\vec{\pp}$-boxes of level $\ell$ and refinement domain $\Omega_\ell$ is studied.
Section \ref{sec:equivalence-b-splines} uses this relation to show that B-splines of different levels span the same polynomial spaces over mesh-elements.
Section \ref{sec:well-behaved-border-elements} builds on top of this relation by identifying the starting element $e_1$, which we call a well-behaved border element, and constructing an ordering from it. 
Which is used in the proof of Theorem \ref{thm:non-overloaded-well-behaved-p-box}.

\begin{figure}[htb]
    \centering
    \begin{subfigure}{.3\textwidth}
        \centering
        \begin{tikzpicture}
            \definecolor{meshcolor}{RGB}{127,127,127}
            \draw (1,1) circle (1.5mm);
            \draw[xstep = 1cm, ystep = 1cm,meshcolor] (-0.25,-0.25) grid (3.25,3.25);
            \draw[xstep = 0.5cm, ystep = 0.5cm,meshcolor] (1,-0.25) grid (3.25,3.25);
            \draw[thick,meshcolor] (1,1) rectangle (2.5,2.5);
            \foreach \i in {1.25,1.75,2.25}{
                \node at (1.25,\i) {$\times$};
            }
            \foreach \i in {1.25,1.75}{
                \draw[->,shorten <= 0.1cm, shorten >= 0.05cm] (\i,1.25) -- ++ (0.5,0);
                \draw[->,shorten <= 0.1cm, shorten >= 0.05cm] (\i,1.75) -- ++ (0.5,0);
                \draw[->,shorten <= 0.1cm, shorten >= 0.05cm] (\i,2.25) -- ++ (0.5,0);
            }
            \draw[line width = 0.5mm] (1,1) -- (1,2.5);
        \end{tikzpicture}
        \caption{Edge $\pset{\vec{r}_1,\ell}{}$}
        \label{fig:example-border-non-overloaded-side}
    \end{subfigure}
    \hfill
    \begin{subfigure}{.3\textwidth}
        \centering
        \begin{tikzpicture}
            \definecolor{meshcolor}{RGB}{127,127,127}
            \draw (1,1) circle (1.5mm);
            \draw[xstep = 1cm, ystep = 1cm,meshcolor] (-0.25,-0.25) grid (3.25,3.25);
            \draw[xstep = 0.5cm, ystep = 0.5cm,meshcolor] (1,-0.25) grid (3.25,3.25);
            \draw[xstep = 0.5cm, ystep = 0.5cm,meshcolor] (-0.25,1) grid (3.25,3.25);
            \draw[thick,meshcolor] (1,1) rectangle (2.5,2.5);
            \node at (1.25,1.25) {$\times$};
            \foreach \i in {1.25,1.75}{
                \draw[->,shorten <= 0.1cm, shorten >= 0.05cm] (\i,1.25) -- ++ (0.5,0);
                \draw[->,shorten <= 0.1cm, shorten >= 0.05cm] (\i,1.75) -- ++ (0.5,0);
                \draw[->,shorten <= 0.1cm, shorten >= 0.05cm] (\i,2.25) -- ++ (0.5,0);
                \draw[->,shorten <= 0.1cm, shorten >= 0.05cm] (1.25,\i) -- ++ (0,0.5);
                \draw[->,shorten <= 0.1cm, shorten >= 0.05cm] (1.75,\i) -- ++ (0,0.5);
                \draw[->,shorten <= 0.1cm, shorten >= 0.05cm] (2.25,\i) -- ++ (0,0.5);
            }
            \fill[] (1,1) circle (0.07);
        \end{tikzpicture}
        \caption{Concave corner $\pset{\vec{r}_2,\ell}{}$}
        \label{fig:example-border-non-overloaded-corner-1}
    \end{subfigure}
    \hfill
    \begin{subfigure}{.3\textwidth}
        \centering
        \begin{tikzpicture}
            \definecolor{meshcolor}{RGB}{127,127,127}
            \draw (1,1) circle (1.5mm);
            \draw[xstep = 1cm, ystep = 1cm,meshcolor] (-0.25,-0.25) grid (3.25,3.25);
            \draw[xstep = 0.5cm, ystep = 0.5cm,meshcolor] (1,1) grid (3.25,3.25);
            \draw[thick,meshcolor] (1,1) rectangle (2.5,2.5);
            \node at (1.25,1.25) {$\times$};
            \foreach \i in {1.25,1.75}{
                \draw[->,shorten <= 0.1cm, shorten >= 0.05cm] (\i,1.25) -- ++ (0.5,0);
                \draw[->,shorten <= 0.1cm, shorten >= 0.05cm] (\i,1.75) -- ++ (0.5,0);
                \draw[->,shorten <= 0.1cm, shorten >= 0.05cm] (\i,2.25) -- ++ (0.5,0);
                \draw[->,shorten <= 0.1cm, shorten >= 0.05cm] (1.25,\i) -- ++ (0,0.5);
                \draw[->,shorten <= 0.1cm, shorten >= 0.05cm] (1.75,\i) -- ++ (0,0.5);
                \draw[->,shorten <= 0.1cm, shorten >= 0.05cm] (2.25,\i) -- ++ (0,0.5);
            }
            \draw[line width = 0.5mm] (1,2.5) -- (1,1) -- (2.5,1);
        \end{tikzpicture}
        \caption{Convex corner $\pset{\vec{r}_3,\ell}{}$}
        \label{fig:example-border-non-overloaded-corner-2}
    \end{subfigure}
    \caption{In two dimensions, there are three distinct border $\vec{\pp}$-boxes $\pbox{\vec{r},\ell}{}$ based on $\partial \pbox{\vec{r},\ell}{} \cap \left( \partial\Omega_\ell \backslash \partial \Omega\right)$, shown with a black line or a filled circle.
    Please note that each boundary includes a unique corner vertex, which is indicated by an empty circle.
    Here $\partial \pbox{\vec{r},\ell}{} \cap \left( \partial\Omega_\ell \backslash \partial \Omega\right)$ is an edge (a), a concave corner (b) or a convex corner (c), and the well-behaved border elements are marked with a $\times$.
    For each well-behaved border element $\vec{e}_\ell$, the mesh-elements of $\widehat{E}(\vec{e}_\ell)$ are given by translation vectors starting from $\vec{e}_\ell$, showing an ordering of the mesh-elements.}
    \label{fig:example-border-non-overloaded}
\end{figure}

\subsection{Refinement border characterization}\label{sec:refinement-border-characterization}
For any given mesh element or $\vec{\pp}$-box, its boundary can be decomposed into vertices, edges, faces, etc.
Note that intersections of multiple hyperplanes reproduce each of these.
As a result, each possible boundary consisting of combinations of vertices, edges, faces, etc, can be reduced to appropriate unions of hyperplane intersections.
In this section, we exploit the refinement pattern of pboxes/THB-splines to introduce an element/pbox boundary description in terms of hyperplanes.  

Consider level $\ell>1$, and an active level-$\ell$ (macro-)element $\Omega^{\kappa} \subset \Omega_\ell$.
E.g., a mesh element for $\kappa = \vec{e}_\ell \in \mesh{\DomainHierarchy}{}$ or a $\vec{\pp}$-box for $\kappa = \pset{\vec{r},\ell}{}\in\meshPbox{\DomainHierarchy}{\vec{\pp}}$.
{By construction, $\Omega^{\kappa} = \bigtimes_{i=1}^n \left(\widehat{\xi}_{k^i,\ell}^i,\widehat{\xi}_{k^i+t^i,\ell}^i \right)$ for indices $k^i,t^i$.}
This (macro-)element is obtained by bisecting a level-$(\ell-1)$ (macro-)element, we will denote it with $\Omega^{\star \kappa}$ {$ = \bigtimes_{i=1}^n\left(\widehat{\xi}_{K^{i},\ell-1}^i,\widehat{\xi}_{K^{i}+T^i,\ell-1}^i \right)$ for indices $K^{i},T^{i}$.}
Moreover, {for each $i$, either $\xi_{k^i,\ell}^i = \xi_{K^{i},\ell-1}^i$ or $\xi_{k^{i}+t^i,\ell}^i = \xi_{K^{i}+T^i,\ell-1}^i$. Hence, $\kappa$ and $\star \kappa$ share a corner vertex/extreme point, that is, a level-$(\ell -1)$ breakpoint, which we denote as $\vec{v}$.
}
Since $\Omega^{\kappa}$ is active, we must have $\Omega^{\kappa} \subset \Omega^{\star \kappa} \subseteq \Omega_\ell$ so that we can decompose $\partial \Omega^\kappa \cap \partial \Omega_\ell \backslash \partial \Omega$ into subsets defined as the intersection of one or multiple hyperplanes containing the corner vertex $\vec{v}$; see Figure \ref{fig:example-border-non-overloaded}.
\begin{definition}\label{def:shared-boundary-set}
    Given $\Omega^{\kappa}$, let $\vec{v}$ be the corner vertex and breakpoint $\Omega^\kappa$ shares with $\Omega^{\star \kappa}$.
    For a set of vectors $\mathfrak{n} := \{\vec{n}_1, \vec{n}_2, \dots\}$, $H(\kappa, \mathfrak{n})$ is called a shared-boundary set and it is defined as
    \begin{equation}
        H(\kappa, \mathfrak{n}) := \left\{ P_{\vec{v}, \vec{n}} \cap \left( \partial \Omega^\kappa \backslash \partial \Omega \right) : \vec{n}\in \mathfrak{n} \right\}\,,
    \end{equation}
    where $P_{\vec{v}, \vec{n}}$ are intersections of hyperplanes passing through $\vec{v}$,
    \begin{equation}\label{eq:hyperplane-definition}
        \begin{split}
            P_i &:= \mathbb{R}^{i-1}\times\left\{ v_{i} \right\}\times\mathbb{R}^{n-i}\;,\\
            P_{\vec{v},\vec{n}} &:= \bigcap\limits_{i:n_i\neq 0} P_i\;.
        \end{split}
    \end{equation}
\end{definition}
When for a level-$\ell$ (macro-)element $\kappa$, $H(\kappa,\mathfrak{n})=\partial \Omega^{\kappa} \cap \left( \partial \Omega_\ell \backslash \partial \Omega \right)$ and $\mathfrak{n}$ has the smallest cardinality, we will denote it as $\mathfrak{n}^{\mathrm{b}}$ or $\mathfrak{n}^{\mathrm{b}}_\kappa$ when $\kappa$ is not clear from context.
For the $\vec{\pp}$-boxes of Figure \ref{fig:example-border-non-overloaded},
\begin{equation}
    \text{a)}\,\mathfrak{n}^{\mathrm{b}}_{\pset{\vec{r}_1,\ell}{}}=\left\{ \begin{bmatrix}
        1\\ 
        0
    \end{bmatrix} \right\}\;,\quad
    \text{b)}\,\mathfrak{n}^{\mathrm{b}}_{\pset{\vec{r}_2,\ell}{}}=\left\{ \begin{bmatrix}
        1\\ 
        1
    \end{bmatrix} \right\}\;,\quad
    \text{c)}\,\mathfrak{n}^{\mathrm{b}}_{\pset{\vec{r}_3,\ell}{}}=\left\{ \begin{bmatrix}
        1\\ 
        0
    \end{bmatrix}, \begin{bmatrix}
        0\\ 
        1
    \end{bmatrix} \right\}\;.
\end{equation}

\subsection{Equivalence of B-spline subspaces of different levels }
\label{sec:equivalence-b-splines}
On a given mesh element ${\vec{e}_\ell} \in \mesh{\ell}{}$, the supported B-splines of level-$(\ell-1)$ and level-$\ell$ are both linearly independent, count the same number of basis functions, and span the same space.
This same relation between splines of different levels also holds when imposing boundary conditions on element $\Omega^e$, such that there are fewer supported splines.
For example, on boundaries described by shared-boundary sets such as $H(\vec{e}_\ell,\mathfrak{n}^{\mathrm{b}})$.

The boundary conditions we wish to impose are inspired by the definition of $\Bbasis{\ell}{\act,k}$. 
Let $\Bspline{\vec{i},\ell}{}\in \Bbasis{\ell}{\act,k}$, so that $\supp{\Bspline{\vec{i},\ell}{}}\subseteq \Omega_k$. 
Then, for point $\vec{x}\in\partial \Omega_k$ with axis aligned normal $\vec{n} = e_i$,
\begin{equation}\label{eq:boundary-condition}
    \Bspline{\vec{i},\ell}{}(\vec{x}) = 0\;, \nabla_{\vec{n}} \Bspline{\vec{i},\ell}{}(\vec{x}) = 0\;, \dots , \nabla_{\vec{n}}^{\pp^i - \kk^i}  \Bspline{\vec{i},\ell}{}(\vec{x}) = 0\;,
\end{equation}
where $\nabla_{\vec{n}}$ is the partial derivative along the vector $\vec{n}$.
Observe that the B-splines that satisfy \eqref{eq:boundary-condition} vanish at a facet of $\vec{e}_\ell$ with normal $\vec{n}=e_i$, or even a hyperplane with normal $\vec{n}=e_i$ going through $\vec{x}$.
We can naturally extend \eqref{eq:boundary-condition} to $k$-dimensional sets (e.g. vertices and edges) for $k<n-1$ by considering them as intersections of hyperplanes.

Intersecting hyperplanes to construct lower-dimensional boundary sets was the motivation behind the definition of shared-boundary sets.
Given a shared-boundary set $H(\vec{e}_\ell,\mathfrak{n})$ of element $\Omega^{\vec{e}_\ell}$ with corner vertex $\vec{v}$ and a supported B-spline $\Bspline{}{}$. Call $\Bspline{}{}$ vanishing $H(\vec{e}_\ell,\mathfrak{n})$, if for each $\vec{n}\in \mathfrak{n}$, there is at least one $i$ such that $n^i\neq 0$ and:
\begin{align}\label{eq:vanishing-condition}
    \nabla_{\vec{e}^i}^0 \Bspline{}{}(\vec{x}) = 0\;,\dots, \nabla_{\vec{e}^i}^{\pp^i-\kk^i} \Bspline{}{}(\vec{x}) = 0\;,\quad \forall \vec{x}\in P_{\vec{v},\vec{n}}\;.
\end{align}
Since B-splines are piecewise-polynomial tensor-product functions, this definition places restrictions on the polynomial components of the spline $\Bspline{}{}$. 
In Figure \ref{fig:example-vanishing-polynomials}, an example of polynomials is given for shared-boundary sets of hyperplanes.
For more general shared-boundary sets, the resulting polynomials are combinations of those described in Figure \ref{fig:example-vanishing-polynomials}.
\begin{figure}[htb]
    \centering
    \begin{subfigure}{0.45\textwidth}
        \centering
        \begin{tikzpicture}
            \draw (0,0) rectangle (1,1);
            \draw[line width = 1mm, red] (0,0) -- (0,1);
            \node at (0.5,-0.3) {$\{x,x^2,x^3\}\otimes \{1,y,y^2\}$};
        \end{tikzpicture}
    \end{subfigure}
    \hfill
    \begin{subfigure}{0.45\textwidth}
        \centering
        \begin{tikzpicture}
            \draw (0,0) rectangle (1,1);
            \draw[line width = 1mm, red] (0,0) -- (1,0);
            \node at (0.5,-0.3) {$\{1,x,x^2,x^3\}\otimes \{y^2\}$};
        \end{tikzpicture}
    \end{subfigure}
    \caption{A polynomial example of vanishing over a boundary set. Over the element $\Omega^e = [0,1]\times [0,1]$, two different shared-boundary sets are drawn, indicated by a thick red line. For the polynomial degrees $\pp = [3,2]^T$ and continuity $\kk = [1,2]^T$, the tensor product decomposition of the vanishing polynomials is given. }
    \label{fig:example-vanishing-polynomials}
\end{figure}

For vanishing B-splines, we have the following result.
\begin{lemma}\label{lem:B-spline-polynomomial-change-of-basis}
    Given an element $\vec{e}_\ell \in\mesh{\ell}{}$ of a B-spline space of degree $\vec{\pp}$ and internal knot multiplicity $\vec{\kk}$ and a shared-boundary set $H(\vec{e}_\ell, \mathfrak{n})$.
    The subset of supported level-$(\ell-1)$ B-spline basis functions that vanish over $H\left(\vec{e}_\ell,\mathfrak{n}\right)$ is a change of basis to the subset of supported level-$\ell$ B-spline basis functions that vanish over $H\left(\vec{e}_\ell,\mathfrak{n}\right)$.
\end{lemma}
\begin{proof}
    It is sufficient to show that B-splines of different levels are changes of basis when a single hyperplane defines the shared-boundary set.
    For general shared-boundary sets, the vanishing splines are combinations of splines that vanish over different hyperplanes, and do not depend on the spline level.

    For a single hyperplane, the normal is axis-aligned, and the tensor product B-splines are restricted in one dimension. 
    Without loss of generality, we consider the one-dimensional case.
    Then, $\Omega^e := (\xi_{k,\ell},\xi_{k+\kk,\ell}) = (\xi_{K,\ell-1},\frac{1}{2}\xi_{K+\kk,\ell-1})$ for the lowest indices $k, K$ and we pick $H(e,[+1])=\xi_{k,\ell}$. 
    Note, the element is written in terms of knot-indices instead of break-points, so that the right mesh element boundary is $\kk$ knots over.
    Then, the subsets of vanishing B-splines and polynomials are:
    \begin{align}
        \left\{\Bspline{k,\ell}{},\dots,\Bspline{k+\kk-1,\ell}{}\right\} &\subset \Bbasis{\ell}{}\;,\\
        \left\{\Bspline{K,\ell-1}{},\dots,\Bspline{K+\kk-1,\ell-1}{}\right\} &\subset \Bbasis{\ell-1}{}\;,\\
        \left\{ \left(x-\xi_{k,\ell}\right)^{\pp-\kk+1},\dots,\left(x-\xi_{k,\ell}\right)^\pp \right\}&\subset \mathcal{P}_{\pp}(\Omega^e)\;.
    \end{align}
    Each of these function sets spans the same space with the same number of basis functions. Ending our proof.
\end{proof}
We end this section with the following relation to vanishing splines and $\Bbasis{\ell}{\act,\ell}$.
\begin{lemma}
\label{lem:pbox-active-vanishing-equivalence}
    For an active level-$\ell$ $\vec{\pp}$-box $\pset{\vec{r},\ell}{}$, the set of B-splines of $\Bbasis{\ell}{\act,\ell}$ with support on $\pset{\vec{r},\ell}{}$ and the set of level-$\ell$ B-splines with support on $\pset{\vec{r},\ell}{}$ that vanish over $H(\pset{\vec{r},\ell}{},\mathfrak{n}^{\mathrm{b}})$ are a basis for the same space.
\end{lemma}
\begin{proof}
By construction, every B-spline of $\Bbasis{\ell}{\act,\ell}$ vanishes over $H(\pset{\vec{r},\ell}{},\mathfrak{n}^{\mathrm{b}})$.
For the converse, fix the level-$\ell$ B-spline $\Bspline{}{}$ with support on $\pbox{\vec{r},\ell}{}$ and which vanishes over $H(\pset{\vec{r},\ell}{},\mathfrak{n}^{\mathrm{b}})$.
Without loss of generality, for each $\vec{n}\in \mathfrak{n}^{\mathrm{b}}$ we assume that all entries are non-negative.
We note that the support of a level-$\ell$ B-spline is at most $\pp^i+1$ level-$\ell$ elements wide, so that it can only be supported by $\pset{\vec{r},\ell}{}$ and its direct neighbours.
Then, let $\pset{\vec{k},\ell-1}{} = \star\pset{\vec{r},\ell}{}$ so that $\Bspline{}{}$ can only have support on $\pset{\tilde{\vec{k}},\ell-1}{}$ for $\tilde{k}^i = k^i, k^i -1$ for all $i$, provided $\pset{\tilde{\vec{k}},\ell-1}{}$ exists. If for all $\pset{\tilde{\vec{k}},\ell-1}{}$, $\pbox{\tilde{\vec{k}},\ell-1}{}\subset \Omega_\ell$ we are done. 
Else, consider a $\pset{\tilde{\vec{k}},\ell-1}{}$ so that $\pbox{\tilde{\vec{k}},\ell-1}{}\cap\Omega^\ell=\emptyset$. 
Then, $\emptyset \neq \partial \pbox{\tilde{\vec{k}},\ell-1}{} \cap \partial \pbox{\vec{r},\ell}{} \subset H(\pset{\vec{r},\ell}{},\mathfrak{n}^{\mathrm{b}})$ and $\Bspline{}{}$ will vanish over $\partial \pbox{\tilde{\vec{k}},\ell-1}{} \cap \partial \pbox{\vec{r},\ell}{}$, showing that $\Bspline{}{}$ is not supported on $\pset{\tilde{\vec{k}},\ell-1}{}$. Showing that $\supp{\Bspline{}{}}\subset \Omega_\ell$.
\end{proof}

\subsection{Well-behaved border (macro-)elements}\label{sec:well-behaved-border-elements}
The following definition introduces well-behaved border elements and the associated well-behaved border macro-elements.
\begin{definition}\label{def:well-behaved-border-macro-element}
    Given an active border $\vec{\pp}$-box $\pset{\vec{r},\ell}{}$, call an element ${\vec{e}_{\ell}}\in \pset{\vec{r},\ell}{}$ a well-behaved border element if,
    \begin{equation}
        \partial \Omega^{\vec{e}_{\ell}} \cap P_{\vec{v},\vec{n}} \neq \emptyset,\quad \forall P_{\vec{v},\vec{n}} \in H\left(\pset{\vec{r},\ell}{},\mathfrak{n}^{\mathrm{b}}\right)\;.
    \end{equation}
    For every element $\tilde{\vec{e}}_\ell \in \pset{\vec{r},\ell}{}$ that is not a well-behaved border element, assign it to the set $\widehat{E}({\vec{e}_\ell})\subset \mesh{\DomainHierarchy}{}$ if ${\vec{e}_\ell}\in \pset{\vec{r},\ell}{}$ is the closest (in the $l^1$ metric) well-behaved border element to $\tilde{\vec{e}}_\ell$. By convention, $\widehat{E}({\vec{e}_\ell})$ also contains the well-behaved border element $\vec{e}_\ell$.
\end{definition}
Note that the closest well-behaved border element in the above definition is unique.
Indeed, for border $\vec{\pp}$-box $\pset{\vec{r},\ell}{}$, let ${\vec{e}_\ell}, {\vec{e}_\ell+\vec{t}} \in\pset{\vec{r},\ell}{}$ be well-behaved border elements, then we must have
$$
    \partial \Omega^{{\vec{e}_{\ell}}} \cap P_{\vec{v},\vec{n}} \neq \emptyset \neq \partial \Omega^{{\vec{e}_{\ell}}+\vec{t}} \cap P_{\vec{v},\vec{n}}\;,\quad \forall P_{\vec{v},\vec{n}} \in H\left({\pset{\vec{r},\ell}{}}{},{\mathfrak{n}^{\mathrm{b}}}\right)\;,
$$ 
which can only occur if $\vec{t} \perp \vec{n}$ for all $\vec{n} \in \mathfrak{n}^{\mathrm{b}}$.
This observation allows for a unique, closest, well-behaved border element. 
For these macro-elements, we have the following elementwise linear independence result.

\begin{proposition}\label{prop:lin-ind-propagation}
Let $\pset{\vec{r},\ell}{}$  be an active border $\vec{\pp}$-box and consider a well-behaved border element ${\vec{e}_{\ell}}\in \pset{\vec{r},\ell}{}$.
Then, for any $\vec{t}_1 = (t_1^1,\dots,t_1^n)$ such that ${{\vec{e}_{\ell}}+\vec{t}_1} \in \widehat{E}({\vec{e}_{\ell}})$, the following set of THB-splines on $\Omega^{{\vec{e}_{\ell}}+\vec{t}_1}$ is linearly independent,
\begin{equation}\label{eq:prop-lin-ind-propagation-THB-splines}
	\begin{split}
		\bigg\{
			\Tspline{j}{}\,:\,
			&\supp{\Tspline{j}{}} \cap \Omega^{{\vec{e}_{\ell}}+\vec{t}_1} \neq \emptyset,
			\text{ and, }\\
			&\supp{\Tspline{j}{}} \cap \Omega^{{\vec{e}_{\ell}}+\vec{t}} = \emptyset,
			\text{ where }
			{{\vec{e}_{\ell}}+\vec{t}} \in \widehat{E}_{\vec{e}_{\ell}},\,
			\vec{t}_1 \neq \vec{t}\,\land\,
			|t_1^i| \geq |t^i| \, \forall i\,
			\bigg\}\,.
	\end{split}
\end{equation}
\end{proposition}
The proof of Proposition \ref{prop:lin-ind-propagation} for well-behaved elements is a direct result of Lemma \ref{lem:B-spline-polynomomial-change-of-basis} and \ref{lem:pbox-active-vanishing-equivalence}. However, for general mesh elements, we require some intermediate results.
First, we have that for most elements, Proposition \ref{prop:lin-ind-propagation} is trivial, due to the following Lemma from \cite{lanini_characterization_2024}.
\begin{lemma}[{\cite[Lemma\,1.3.1]{lanini_characterization_2024}}]\label{lem:non-decreasing-trunc}
Let ${\vec{e}_{\ell}}$ be a well-behaved border element with macro-element $\widehat{E}(\vec{e}_\ell)$.
In addition, let ${\vec{e}_\ell + \vec{t}_1}\in\widehat{E}({\vec{e}_\ell})$ and $\Omega^{\vec{e}_\ell +\vec{t}}\subset \Omega^{\star(\vec{e}_\ell+\vec{t}_1)}$ with $\vec{t}_1\neq\vec{t}$, $|t^i|\leq |t_1^i|$ for all $i$.
Then, for any B-spline $\Bspline{}{}\in \Bbasis{\ell-1}{\nact,\ell}$, we have that if $\Omega^{\vec{e}_\ell+\vec{t}_1}\subset \supp{\trunc{\ell}{\Bspline{}{}}}$, then $\Omega^{\vec{e}_\ell+\vec{t}}\subset \supp{\trunc{\ell}{\Bspline{}{}}}$.
\end{lemma}
And as a consequence.
\begin{corollary}\label{cor:trivial-odd-translation-vector}
    Consider a well-behaved border element ${\vec{e}_{\ell}}$ with macro-element $\widehat{E}({\vec{e}_\ell})$. Then, for the element ${\vec{e}_\ell + \vec{t}_1}\in\widehat{E}({\vec{e}_\ell})$ with any $t^i_1$ odd, then the set of splines given in \eqref{eq:prop-lin-ind-propagation-THB-splines} of Proposition \ref{prop:lin-ind-propagation} consists solely out of level-$\ell$ B-splines.
\end{corollary}
\begin{proof}
For Proposition \ref{prop:lin-ind-propagation}, we only consider the splines that have no support on $\Omega^{\vec{e}_\ell + \vec{t}}$ for $\vec{t}_1\neq \vec{t}$, $|t^i|\leq |t_1^i|$ for all $i$.
In addition, since the mesh is assumed to be of admissibility class $2$, the only possible splines are either level-$\ell$ B-splines or truncated level-$(\ell-1)$ B-splines.
Pick $\vec{t}$ as the vector with only even entries, by subtracting one from the odd entries of $\vec{t}_1$; see Figure \ref{fig:example-lemma-non-decreasing-trunc}.
Then, $\Omega^{\vec{e}_\ell+\vec{t}}\subset \Omega^{\star ( \vec{e}_\ell + \vec{t}_1)}$ and by Lemma \ref{lem:non-decreasing-trunc} each truncated level-$(\ell-1)$ B-spline is supported on $\Omega^{\vec{e}_\ell + \vec{t}}$.
\end{proof}
For the remaining elements, the difficulty is in incorporating the fact that the THB-splines cannot have support on prior elements (with prior as in Proposition \ref{prop:lin-ind-propagation}).  
Within the tools we have introduced, this can be achieved by introducing a new shared-boundary set. 
The relevant splines vanish over the hyperplanes separating these elements, as shown in Figure \ref{fig:example-lemma-non-decreasing-trunc}.
These hyperplanes pass through the common vertex of ${{\vec{e}_{\ell}}+\vec{t}_1}$ and ${\star({\vec{e}_{\ell}}+\vec{t}_1)}$, so that they are described by $H\left(\vec{e}_\ell + \vec{t}_1,\mathfrak{n}^{\mathrm{t}}_{\vec{t}}\right)$:
\begin{eqnarray}
    \mathfrak{n}^{\mathrm{t}}_{\vec{t}} := \left\{ \, t^i\hat{\vec{i}}\, :\,t^i\neq 0 \,\right\}\;,
\end{eqnarray}
where $\hat{\vec{i}}$ is the dimension-$i$ axis-aligned unit vector.
Then, the relevant B-splines with support on ${\vec{e}_\ell + \vec{t}_1}$ and without support on any $\vec{e}_\ell + \vec{t}\in \widehat{E}({\vec{e}_\ell})$, $\vec{t}\neq \vec{t}_1$, $|t^i| \leq |t^i_1|$ for all $i$, vanish over $H\left(\vec{e}_\ell + \vec{t}_1,\mathfrak{n}^{\mathrm{t}}_{\vec{t}_1}\right)$.
\begin{figure}[htb]
\centering
\begin{tikzpicture}
\definecolor{color1}{RGB}{102,194,165}
\definecolor{color2}{RGB}{252,141,98}
\definecolor{color3}{RGB}{141,160,203}

\fill[color3,opacity=0.5] (1,1) rectangle (1.5,1.5);
\draw[pattern={north west lines},pattern color=color1] (1,1) rectangle (1.5,1.5);
\draw[pattern={north west lines},pattern color=color1] (1,2) rectangle (1.5,2.5);
\draw[pattern={north west lines},pattern color=color1] (2,1) rectangle (2.5,1.5);
\draw[pattern={north west lines},pattern color=color1] (2,2) rectangle (2.5,2.5);

\draw[step = 0.5, lightgray] (1,1) grid (3,3);
\draw (0,0) grid (3,3);
\draw[color2,thick] (1,1) rectangle (2.5,2.5);
\end{tikzpicture}
\caption{A well-behaved border element (highlighted in purple) and the outline of the accompanying projection element shaded/drawn in orange. The four elements in the projection element (indicated with green sideways lines) are the exception to Lemma \ref{lem:non-decreasing-trunc}.}
\label{fig:example-lemma-non-decreasing-trunc}
\end{figure}
We can now state the proof of Proposition \ref{prop:lin-ind-propagation}.
\begin{proof}[Proof of Proposition \ref{prop:lin-ind-propagation}]
By Corollary \ref{cor:trivial-odd-translation-vector}, the claim is immediate if at least one $t_1^i$ is odd.
Then, let $\vec{t}_1$ be the translation vector where each $t_1^i$ is even.
In the following, we use $\vec{t}$ to denote a translate vector such that ${{\vec{e}_{\ell}}+\vec{t}} \in \widehat{E}({\vec{e}_{\ell}})$, $\vec{t}_1 \neq \vec{t}$, $|t^i|\leq |t_1^i| $ for all $i$.
For level $k$, define the set of level-$k$ B-splines, without support on any such ${\vec{e}_\ell + \vec{t}}$,
\begin{align}\label{eq:proof-lin-ind-start-bsplines}
\begin{split}
\Bbasis{k}{\#} &:= \big\{\Bspline{k}{}\in\Bbasis{k}{}:\Omega^{\vec{e}_\ell+\vec{t}_1} \subset \supp{{\Bspline{k}{}}},\\ &\quad\quad\quad  \,\Omega^{\vec{e}_\ell+\vec{t}} \nsubset \supp{\Bspline{k}{}},\,\forall \vec{e}_\ell + \vec{t}\in\widehat{E}({\vec{e}_\ell}),\,\vec{t}\neq\vec{t}_1 ,\,|t^i|\leq |t_1^i|\,\forall i\, \big\}\;,
\end{split}\\
    \Bbasis{k}{\mathrm{t}} &:= \left\{\Bspline{k}{\#}\in\Bbasis{k}{}:\Omega^{\vec{e}_\ell+\vec{t}_1} \subset \supp{\trunc{\ell}{\Bspline{k}{}}}\, \right\}\;.
\end{align}
Note that $\trunc{\ell}{\Bspline{}{}}$ is a linear combination of $\Bbasis{\ell}{\nact,\ell}$, so that the splines of $\Bbasis{\ell-1}{\mathrm{t}}$ do not vanish over $H(\vec{e}_{\ell}+\vec{t}_1,\mathfrak{n}^{\mathrm{b}})$, and are thus linearly independent to the level-$(\ell-1)$ B-splines that do vanish over $H(\vec{e}_{\ell}+\vec{t}_1,\mathfrak{n}^{\mathrm{b}}\cup\mathfrak{n}^{\mathrm{t}}_{\vec{t}_1})$.
By Lemma \ref{lem:B-spline-polynomomial-change-of-basis}, both the level-$(\ell-1)$ and level-$\ell$ B-splines that vanish over $H(\vec{e}_\ell+\vec{t}_1,\mathfrak{n}^{\mathrm{b}}\cup \mathfrak{n}^{\mathrm{t}}_{\vec{t}})$ are a basis of the same polynomial space over ${\vec{e}_\ell+\vec{t}_1}$.
In addition, the level-$\ell$ B-splines that vanish over $H(\vec{e}_\ell+\vec{t}_1,\mathfrak{n}^{\mathrm{b}}\cup \mathfrak{n}^{\mathrm{t}}_{\vec{t}_1})$, vanish over $H(\pset{\vec{r},\ell}{},\mathfrak{n}^{\mathrm{b}})$ and by Lemma \ref{lem:pbox-active-vanishing-equivalence} belong to $\Bbasis{\ell}{\act,\ell}$.
Hence, over the element ${\vec{e}_\ell +\vec{t}_1}$, the set
\begin{equation}
\label{eq:proof-hb-spline-definition}
    \Bbasis{\ell-1}{\mathrm{t}} ~\cup~\left( \Bbasis{\ell}{\#}\cap \Bbasis{\ell}{\act,\ell}\right) 
\end{equation}
is linearly independent.
To construct the THB-splines, the splines of \eqref{eq:proof-hb-spline-definition} should be truncated.
By the assumption of admissibility class $c=2$, the truncation operation \eqref{eq:def-trunc-operator} is equivalent to
\begin{equation}
    \trunc{\ell}{\Bspline{\vec{i},\ell-1}{}} \;=\; \trunc{k}{\sum_{\vec{j}}\Coeff{\vec{i},\vec{j}}{}\Bspline{\vec{j},\ell}{}} \;=\; \Bspline{\vec{i},\ell-1}{} - \;\sum_{\mathclap{\vec{j}:\Bspline{\vec{j},\ell}{}\in\Bbasis{\ell}{\#}\cap \Bspline{\ell}{\act,\ell} }}\; \Coeff{\vec{i},\vec{j}}{} \Bspline{\vec{j},\ell}{}\;,\forall \Bspline{\vec{i},\ell-1}{}\in \Bbasis{\ell-1}{\mathrm{t}} \;.
\end{equation}
Here we use Lemma \ref{lem:B-spline-polynomomial-change-of-basis} to note that the level-$(\ell-1)$ B-splines of $\Bbasis{\ell-1}{\#}$ are linear combinations of $\Bbasis{\ell}{\#}$ on ${\vec{e}_\ell+\vec{t}_1}$.
Hence, $\trunc{\ell}{\cdot}$ is a change of basis for the splines of \eqref{eq:proof-hb-spline-definition}, and
\begin{eqnarray}\label{eq:proof-thb-spline-definition}
    \trunc{\ell}{\Bbasis{\ell-1}{\mathrm{t}}}\;\;\cup\;\; \left( \Bbasis{\ell}{\#} \cap \Bbasis{\ell}{\act,\ell}\right)\;.
\end{eqnarray}
is linearly independent on ${{\vec{e}_{\ell}}+\vec{t}_1}$, and are precisely the THB-splines of \eqref{eq:prop-lin-ind-propagation-THB-splines}.
\end{proof}
We can now prove that active border $\vec{\pp}$-boxes are not overloaded w.r.t. $\Tbasis{\DomainHierarchy}{}$, which is the non-trivial step of proving Lemma \ref{lem:non-overloaded-border-p-box}.
\begin{proof}[Proof of Lemma \ref{lem:non-overloaded-border-p-box}]
    Consider any regular $\vec{\pp}$-box $\pset{\vec{r},\ell}{}$, which by Proposition \ref{prop:regular-p-box-spline-level} only supports B-splines and is thus trivially not overloaded w.r.t. $\Tbasis{\DomainHierarchy}{}$.
    If $\pset{\vec{r},\ell}{}$ is an active border $\vec{\pp}$-box, it can be partitioned into smaller macro-elements $\widehat{E}({\vec{e}_{\ell}})$ for the well-behaved border elements ${\vec{e}_{\ell}}\in\pset{\vec{r},\ell}{}$.
    For such a macro-element $\widehat{E}({\vec{e}_{\ell}})$, consider any linear combination of THB-splines that sums to zero on this macro-element,
    \begin{equation}\label{eq:proof-thm-non-overloaded-projection-element-sum}
     0 = \sum_{j}\Coeff{j}{} \Tspline{j}{}\;.
    \end{equation}
    Starting from the linear independence on the well-behaved border element $\Omega^{\vec{e}_{\ell}}$, which is obtained by considering $\vec{t}=(0,\dots,0)$ in proposition \ref{prop:lin-ind-propagation}, we can conclude that some of the $\Coeff{j}{}$ in \eqref{eq:proof-thm-non-overloaded-projection-element-sum} are zero.
    Consequently, repeated application of Proposition \ref{prop:lin-ind-propagation} reveals that all the remaining $\Coeff{j}{}$ ought to be zero since every element in $\widehat{E}({\vec{e}_{\ell}})$ can be obtained by translations of ${\vec{e}_{\ell}}$.
    A similar argument shows that the THB-splines supported on $\pbox{\vec{r},\ell}{}$ are linearly independent on $\pbox{\vec{r},\ell}{}$.
\end{proof}

\end{document}